\documentclass[]{siamart0516}
 \usepackage{amsmath, amssymb}
 \usepackage{multirow}
 \usepackage{graphicx}
 \usepackage{gensymb}
 \usepackage{url}
 \usepackage{color}
 \usepackage{tikz}
 \usepackage{algorithmic}
 \usetikzlibrary{arrows, decorations.markings, shapes}
 \usepackage[ansinew]{inputenc}
 \usepackage[T1]{fontenc}
 \tikzstyle{realitybox} = [rectangle, rounded corners, minimum width=3cm, minimum height=1cm,text centered, draw=black, fill=white!30]
 \tikzstyle{modelbox} = [rectangle, rounded corners, minimum width=3cm, minimum height=1cm,text centered, draw=black, fill=gray!30]

 \newcommand{\walls}{{\Gamma_{\mathrm{walls}}}}

 \newcommand{\laplace}{\Delta}
 \newcommand{\del}{\nabla}
 \newcommand{\grad}{\del}
 \renewcommand{\div}{\del\cdot}

 \newcommand{\set}[1]{\{ #1 \}}
 \renewcommand{\brack}[1]{\langle #1 \rangle}
 \newcommand{\norm}[1]{\lVert #1 \rVert}
 \newcommand{\abs}[1]{\lvert #1 \rvert}

\ifpdf
\hypersetup{
pdftitle={Variational data assimilation for transient blood flow simulations},
pdfauthor={S. W. Funke, M. Nordaas, {\O}. Evju, M. S. Aln{\ae}s, K.-A. Mardal}
}
\fi

\begin{document}
 \title{Variational data assimilation for transient blood flow simulations}
 \author{S. W. Funke\thanks{Center for Biomedical Computing, Simula Research Laboratory, Norway (\email{simon@simula.no})} \and
         M. Nordaas\thanks{Center for Biomedical Computing, Simula Research Laboratory, Norway} \and
         {\O}. Evju\thanks{Center for Biomedical Computing, Simula Research Laboratory, Norway} \and
         M. S. Aln{\AE}s\thanks{Center for Biomedical Computing, Simula Research Laboratory, Norway} \and
         K.-A. Mardal\thanks{Department of Mathematics, University of Oslo, Norway}}
 \maketitle

\begin{abstract}
Several cardiovascular diseases are caused from localised abnormal blood flow
such as in the case of stenosis or aneurysms. Prevailing theories propose
that the development is caused by abnormal wall-shear stress in focused
areas.  Computational fluid mechanics have arisen as a promising tool for
a more precise and quantitative analysis, in particular because the
anatomy is often readily available even by standard imaging techniques
such as magnetic resolution and computed tomography angiography.  However,
computational fluid mechanics rely on accurate boundary conditions which
is difficult to obtain.  In this paper we address the problem of
recovering high resolution information from noisy, low-resolution
measurements of blood flow using variational data assimilation (also known
as 4DVar).  We show that accurate flow reconstruction is obtained with
proper regularisation even in the presence of significant noise and for a
range of regularisation parameters spanning orders of magnitude.
Numerical experiments are performed in both 2D and 3D and with pulsatile
flow relevant for physiological flow in cerebral aneurysms.
\end{abstract}



 \section{Introduction}

Detailed insight of blood flow has the potential to assist clinical decisions, for
example when evaluating the risk of rupture of an aneurysm \cite{cebral2011quantitative,Xiang2011,takao2012hemodynamic}.
Different
non-invasive measurement techniques for blood flow exist today, such as
ultrasound or phase-contrast magnetic resonance imaging. Unfortunately, the spatial and
temporal resolution of these techniques are still too coarse to unveil
potentially important flow details.
Computational patient-specific blood flow models are promising tools
for obtaining blood flow information with nearly arbitrary high
temporal and spatial resolution. In addition, they allow for computation of non-observable variables, such as the blood pressure and
wall-shear stresses, which are considered important factors in vascular diseases \cite{takao2012hemodynamic,samady2011coronary,dolan2013high,kulcsar2011hemodynamics,samady2011coronary}. The validity of such simulations depends on
the accuracy of the inflow and outflow velocities at
the open vessel boundaries and the segmentation of the vascular geometry.
These model parameters are typically patient-specific and partially or fully
unknown~\cite{ramalho2013,ramalho2012}.

It then seems obvious to incorporate physical measurements into the
model in order to identify the unknown parameters. The result
would be a high-resolution blood flow simulation that  best matches the available
measurements. This idea, known as variational data assimilation,
has been successfully applied to weather and ocean modelling, see for
example~\cite{gopalakrishnan2013}, and to a more general setting for the
optimal control of the Navier-Stokes equations \cite{fursikov1998, lee2011,
guerra2014b}. In the context of blood-flow simulations, this technique has been
applied to the steady-state Navier-Stokes equations, for example in
\cite{tiago2014, lee2011, guerra2014, guerra2014b, amrosi2012, veneziani,
elia2012, john2014}.  However, the steady-state assumption is inadequate for the
pulsatile blood flow in larger arteries.

In this paper we consider the variational data assimilation for \emph{transient} blood flow
models. Section \ref{sec:problem_formulation} formulates the data-assimilation
problem as a mathematical optimisation problem constrained by the Navier-Stokes
equations. Section \ref{sec:numerical_solution} discusses the numerical
details. For the optimisation method, focus will be on achieving
convergence which is independent of the choice of discretisation for the
Navier-Stokes equations. Furthermore, special considerations will be put
on the inclusion of data that are coarse with respect to the time resolution
because the number of samples per cardiac cycle is typically in the order 20-40 while
the number of time steps in a CFD simulation typically is 100-10000.
Section~\ref{sec:examples} demonstrate the feasibility of this approach through
numerical examples in two and three dimensions. This is the
first study where variational data assimilation methods have been used in a
three-dimensional transient blood flow solver.

 \section{Mathematical formulation}\label{sec:problem_formulation}
 \subsection{Blood flow model}
Most computational modelling in cerebral aneurysm studies assume Newtonian flow
with rigid walls, which appear to be adequate \cite{Steinmanreview,
evju2013study}.
Therefore, we model the blood flow through a vessel with the incompressible Navier-Stokes equations
\begin{equation}
    \begin{aligned}
    u_t + (u \cdot \nabla) u - \nu \Delta u + \nabla p & = f \qquad \textrm{in } \Omega \times (0, T], \\
    \nabla \cdot u & = 0 \qquad \textrm{in } \Omega \times (0, T].
    \end{aligned}
    \label{eq:navier_stokes}
\end{equation}
Here, $\Omega \times (0, T]$ is the space-time domain, $u$ and $p$ are the
blood velocity and (scaled) pressure fields, $\nu$ is the (kinematic) viscosity and $f$
describes external body forces.  A more complete blood flow model could
incorporate non-Newtonian effects and the fluid structure interactions between
the blood and the vessel wall~\cite{tricerri2015}.  However, for the purpose of
this paper it is sufficient to consider \eqref{eq:navier_stokes} and to note
that the proposed techniques also apply to more complex models.

\begin{figure}
    \centering
    \includegraphics[width=0.15\textheight]{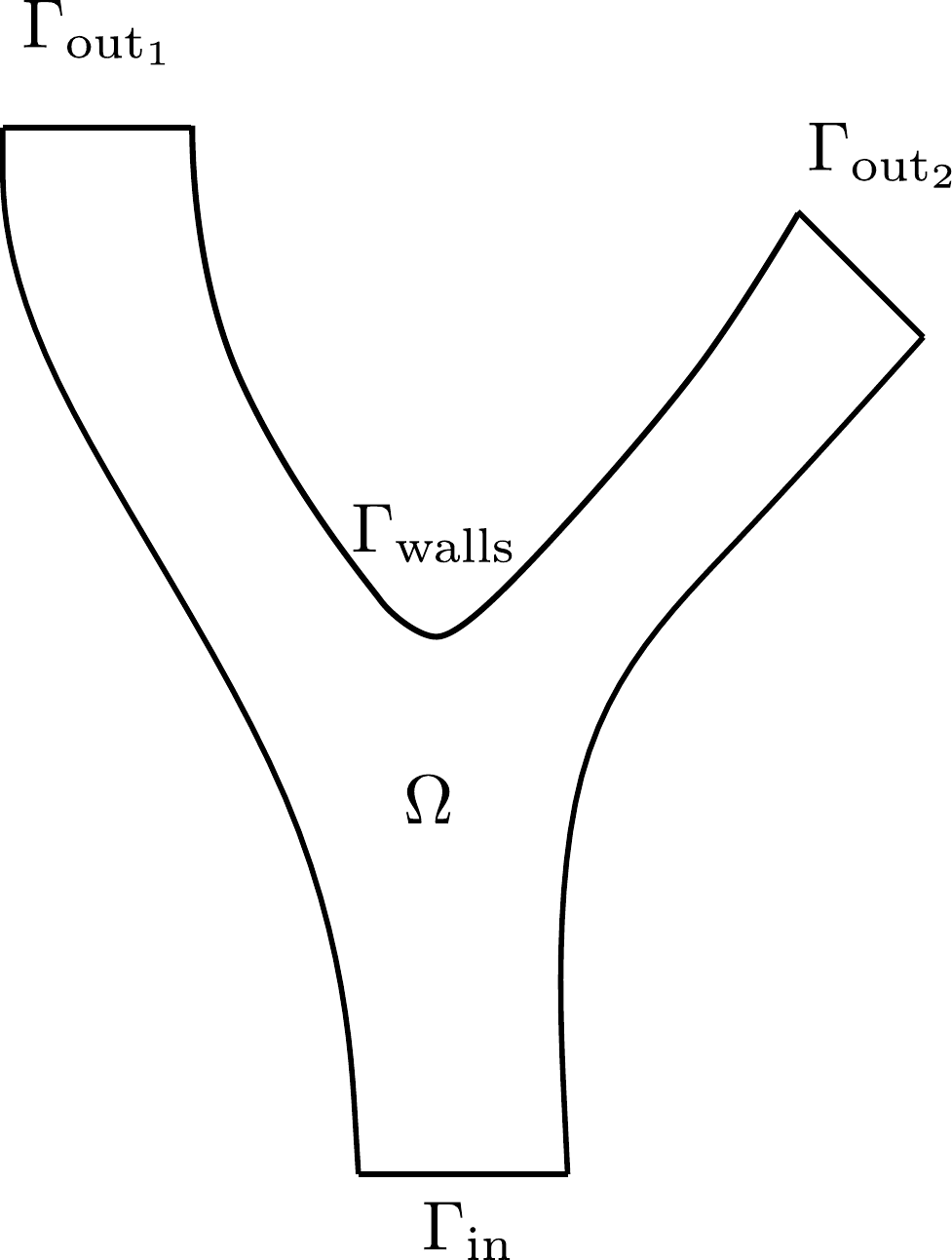}
\caption{The model scenario considered in this paper: a small subset of the artery system with one inlet and two outlet boundaries in 2D and 3D.}\label{img:domain-sketch}
\end{figure}
We only model a small subset of the artery
system and the boundaries of the computational domain consists of a
physical boundary, the vessel walls, as well nonphysical boundaries at inlets
and outlets.  Again for simplicity, we consider the common scenario of one
inlet and two outlets as sketched in figure \ref{img:domain-sketch}. To
close the system, we specify suitable initial and boundary conditions
\begin{subequations}\label{eq:ns_bc}
    \begin{alignat}{4}
        u & = u_0, \qquad           &&\text{on } \Omega\times\set{t=0}, \\
        u & = g_D,                    &&\text{on } \Gamma_D \times (0, T], \label{eq:ns_bc_dirichlet}\\
        pn - \mu \partial_n u &= 0, &&\text{on } \Gamma_{\text{out}_2} \times (0, T],\label{eq:ns_bc_nostress}\\
        u &= 0,                     &&\text{on } \Gamma_{\text{walls}} \times (0, T].
    \end{alignat}
    \label{eq:navier_stokes_bcs}%
\end{subequations}
with normal vector $n$ and a Dirichlet boundary $\Gamma_D := \Gamma_{\text{in}} \cup
\Gamma_{\text{out}_1}$. A traction free boundary is assumed on the outlet $\Gamma_{\text{out}_2}$,
which implies that the vessel is straight in the surroundings of this outlet.

\subsection{Variational data assimilation}
\begin{figure}
    \centering
    \begin{tikzpicture}[thick,node distance=3cm, every node/.style={scale=0.6}]]
      \node[minimum size=0] (x) {}; 
      \node[realitybox,below of=x] (a) {physics};
      \node[modelbox,right of=a, xshift=6cm] (b) {mathematical model};
      \node[modelbox,above of=b, yshift=-0.2cm] (i) {unknown model inputs: $u_0, g_D$};
      \node[realitybox,below of=a] (c) {physical observations: $d$};
      \node[modelbox,below of=b] (d) {modelled observations: $\mathcal{T}(u)$};

      \draw[->,thick] (i) to node[anchor=west] {} (b);
      \draw[->,thick] (a) to node[anchor=west,align=left] {measurement
      device + errors} (c);
      \draw[->,thick] (b) to node[anchor=west,align=left] {virtual
      measurement device $\mathcal{T}$} (d);
      \draw[<->,thick] (c) to node[anchor=north,align=center] {minimise $\|\mathcal{T}(u) -d \|$} (d);
    \end{tikzpicture}
    \caption{Variational data assimilation replicates the observation steps in
    a mathematical model and minimises the discrepancy between measured and
    modelled observations by varying the model inputs.}
    \label{fig:variational_data_assimilation_concept}
\end{figure}
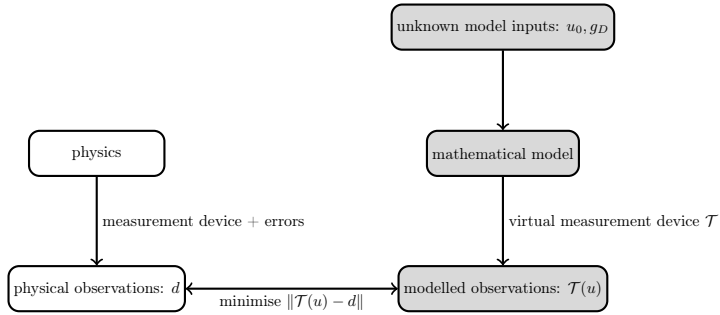

Variational data assimilation is a technique to recover unknown data from given
observation. The idea is to build a model that replicates the steps of
the measurement acquisition, and to tweak the free model parameter so that the
discrepancy between observed and modelled measurements is minimised (figure
\ref{fig:variational_data_assimilation_concept}). Variational data assimilation
for time-dependent problems, or 4DVar, has been successfully applied to problems in meteorology and oceanography and
is used in production services for weather and ocean forecasts and
retrospective analysis.

In the blood flow problem, the aim of variational data assimilation is
to recover 
flow velocity and pressure fields from observational data that is
typically noisy, as well as limited in spatial and temporal
resolution.  These recovered fields should be physically reasonable,
in the sense that they satisfy the mathematical model considered.
More specifically, the initial condition $u_0$ and Dirichlet boundary
condition $g_D$ in \eqref{eq:navier_stokes_bcs} have to be determined
in order to solve the model equations \eqref{eq:navier_stokes}.  In
the data assimilation setting these fields are unknown.  Instead, one
has some given 
observational data $d$, and aim to recover the fields $u_0$
and $g_D$ that best reproduces the data $d$.  In the present paper, it
assumed that the data is a set of $N$ measured velocity
fields on a subdomain $\Omega_{\textrm{obs}} \subset \Omega$.
To be precise, it is assumed that
$d = (d_1,\ldots, d_N) \in L^2({\Omega_\textrm{obs}})^N$. 

The objective is to recover the initial and boundary conditions by minimising the misfit between simulation and measurement data.
Hence we define a goal quantity
\begin{equation}
  J(u) = \norm{\mathcal{T} u - d}^2
  = \sum_{n=1}^{N} \int_{\Omega_\textrm{obs}} \abs{\mathcal{T}_n u - d_n}^2 \, dx,
    \label{eq:goal_quantity}
\end{equation}
where $\mathcal{T}_n$ are the observation operators which approximate the
physical measurement device by mapping the velocity solution $u$ to simulated
measurements.  For instance, if the measurement device takes instantaneously
measurements at $N$ timelevels $t_1, .., t_N$, then the observation operators
should be pointwise evaluations in time. Other common choices for the
observation operators are time and space averaged evaluations of the velocity
state.  The observations might be available only in parts of the domain, hence
the functional integrates over $\Omega_\textrm{obs} \subseteq \Omega$.

We can now formulate the data assimilation problem as an optimisation problem constrained by the Navier-Stokes equations:
\begin{equation}
    \min_{\substack{(u, p) \in Y \\(u_0, g_D)\in M}}\,   J(u) + \mathcal{R}(u_0, g_D) \quad \text{subject to  \eqref{eq:navier_stokes}-\eqref{eq:navier_stokes_bcs}},
    \label{eq:opt_problem}
\end{equation}
where $M$ and $Y$ are suitable function spaces, to be determined later in section \ref{sec:discretisation} below.
The Tikhonov regularisation $\mathcal{R}$ enforce smoothness of the controls:
\begin{equation}
    \mathcal{R}(u_0, g_D)
    = \frac{\alpha}{2} \norm{g_D}_{\Gamma_D\times (0, T]}^2
    + \frac{\gamma}{2} \norm{u_0}_{\Omega}^2,
    \label{eq:regularisation}
\end{equation}
where the coefficients $\alpha$ and
$\gamma$ determine how strongly the problem is regularised in the given norms.

\subsubsection{Choice of norms}\label{sec:norms}
The choice of norms in the regularisation term \eqref{eq:regularisation} specifies the expected regularity of the
reconstructed blood flow. For instance, \cite{john2014} has shown that a unsuitable choice can have a negative impact on
the quality of on the reconstructed data.

The norm used for the initial data is
\begin{align*}
    \norm{u_0}_{\Omega} = \norm{u_0}_{H^1(\Omega)} = \left(\int_\Omega \abs{u_0}^2 + \abs{\grad u_0}^2\, dx\right)^{\frac{1}{2}},
\end{align*}
and for the boundary
\begin{align*}
    \norm{g_D}_{\Gamma_D\times (0, T]}
    = \left(\int_0^T \int_\Omega \abs{g_D}^2 + \abs{\grad g_D}^2 + \abs{\dot g_D}^2  +\abs{\grad \dot g_D}^2 \,dxdt\right)^{\frac{1}{2}}.
\end{align*}
We remark that the norms used require more smoothness on the boundary and initial data than is usually required for the variational formulation of the Navier-Stokes equations, in particular for the time derivative $\dot g_D$.

\section{Numerical solution}\label{sec:numerical_solution}
\subsection{Formulation of the reduced problem}\label{sec:reduced_problem}
Multiple strategies exist for solving the data assimilation problems. One
option is to derive and solve the first order optimality system of
\eqref{eq:opt_problem}. This leads to a large, non-linear system that couples
all spatial and temporal degrees of freedoms of the discretised Navier-Stokes
and adjoint Navier-Stokes equations. Solving this system is numerically
challenging and requires the development of specialised solvers.

The approach taken here is based on the reduced optimisation problem of
\eqref{eq:opt_problem}. The reduced problem is formed by considering the
velocity solution as an implicit function of the initial and boundary controls
by solving the Navier-Stokes equations \eqref{eq:navier_stokes} and
\eqref{eq:navier_stokes_bcs}. We denote this velocity operator as $u(u_0, g_D)$.


To simplify notation, let $m = (u_0, g_D)\in M$
denote the controlled variable.  The functional \eqref{eq:goal_quantity} now has the reduced form
\begin{equation}\label{eq:reduced_functional}
    \hat J(m) := J(u(m)) + \mathcal{R}(m).
\end{equation}
The \emph{reduced optimisation problem} reads
\begin{equation}
    \min_{m\in M}\,   \hat J(m).
    \label{eq:reduced_problem}
\end{equation}
Note that in contrast to \eqref{eq:opt_problem}, the reduced problem is an
unconstrained optimisation problem. As a consequence, it can be solved with
established unconstrained optimisation methods. Also, evaluating the
reduced functional requires the solution of a Navier-Stokes system. Here,
standard solver techniques can be directly applied.

\subsection{Optimisation }
\label{sec:optimisation_method}
The reduced minimisation problem \eqref{eq:reduced_problem} is solved
with the Broyden-Fletcher-Goldbarb-Shanno (BFGS) algorithm. In this
section we present brief overview of the method and its
implementation.

The minimisation problem \eqref{eq:reduced_problem} is iteratively
solved by generating a sequence of points $m_0, m_1,\ldots$, approximating a miniser of
$\hat J$. In each iteration, evaluations of the derivative $D \hat J(m_k) \in M^*$
are used to determine a direction $d_k \in M$ in which the functional is
decreasing. This general descent algorithm in Hilbert spaces is
formulated in algorithm \ref{alg:descent}.

\begin{algorithm}
  \caption{A general descent algorithm in Hilbert spaces, applied to the reduced minimisation problem \eqref{eq:reduced_problem}.}
  \label{alg:descent}
    \begin{algorithmic}
        \STATE{Choose an initial point $x_0$}
  \FOR{$x = 0, 1, \ldots$}
        \STATE{Choose a search direction $d_k = - H_{k} D\hat J(x_k)$}
        \STATE{Choose a step length $\alpha_k>0$ such that $\hat J(x_k+\alpha_kd_k) < \hat J(x_k)$}
        \STATE{Set $x_{k+1} = x_k + \alpha_k d_k$}
    \IF{converged}
        \STATE{return}
    \ENDIF
  \ENDFOR
    \end{algorithmic}
\end{algorithm}

The search direction $d_k = H_k D\hat J(m_k)$ is a descent direction
if the operator $H_k: M^*\rightarrow M$ is positive definite and
self-adjoint. The choice of operators $H_k$ mapping derivatives to
search directions essentially characterises the method.  For example,
taking $H_k = H$ as the Riesz operator for $M$ (i.e. choosing $d_k$ to
be the gradient of $\hat J$ at $m_k$) results in a steepest
algorithm. Setting $H_k = D^2\hat J(m_k)^{-1}$, assuming $\hat J$ is
convex, results in a Newton algorithm.

In the present paper, the algorithm used is the
Broyden-Fletcher-Goldbarb-Shanno (BFGS) algorithm, which is a descent
method of quasi-Newton type.  Quasi-Newton methods have good
convergence properties and do not require evaluations of the Hessian.
Instead, such methods maintains an iteratively constructed
approximation to the inverse of the Hessian. The update formula
specific to the BFGS algorithm is
\begin{equation}
  H_{k+1} = \left(1 - \frac{s_{k+1}\otimes y_{k+1}}{\rho_{k+1}}\right)
               H_k
              \left(1 -  \frac{y_{k+1}\otimes s_{k+1}}{\rho_{k+1}}\right)
              + \frac{s_{k+1}\otimes s_{k+1}}{\rho_{k+1}}
  \label{eq:BFGS_H}
\end{equation}
see e.g. \cite[Chapter 6]{nocedal-wright}. Here, $\otimes:X\times Y \rightarrow \mathcal{B}(Y^*, X)$ denotes the outer product defined by
$(x\otimes y)(z) = x \brack{z, y}_{Y^*,Y}$, for $x\in X$  and $y\in Y$, where $\brack{\cdot, \cdot}_{Y^*, Y}$ denotes duality coupling, and
\begin{equation*}
  \begin{aligned}
    s_{k+1}    &= m_{k+1} - m_k,
    \\
    y_{k+1}  &= D\hat{J}(m_{k+1}) - D\hat{J}(m_k),
    \\
    \rho_{k+1} &= \brack{y_{k+1}, s_k}_{M^*, M}.
  \end{aligned}
\end{equation*}

Note that the initial $H_0:M^*\rightarrow M$ has to prescribed. A natural choice
is to take $H_0$ to be Riesz operator for the space $M$. That is, $H_0$ is the unique operators
such that
\begin{equation}
  (m_0,m_1)_M = \brack{H_0^{-1} m_0, m_1}_{M^*,M} \label{eq:riesz}
\end{equation}
for all $m_0,m_1\in M$. This definition of $H_0$ allows for mesh-independent convergence \cite{schwedes2016},
and is readily seen to coincide with the second order partial derivative of $\hat J$ with respect to $M$. If $D^2 \hat J(m) - H_0^{-1}$ is compact, the method converges superlinearly, see e.g. \cite{kelley-sachs1987}.


For practical implementations, it is common to truncate the update formula
\eqref{eq:BFGS_H} and store only the last $3-10$ pairs of vectors $y_k$ and $s_k$.
The step lengths $\alpha_k$ in algorithm
\ref{alg:descent} should chosen to satisfy the Wolfe conditions, which ensures the convergence of the method \cite[chapter 6]{nocedal-wright}.

\subsection{Discretisation}\label{sec:discretisation}
The optimisation method in section \ref{sec:optimisation_method} requires
evaluations of the reduced functional $\hat J(m)$ and
its derivative $D \hat J(m)$.
Evaluating the reduced functional requires the numerical solution of the Navier-Stokes
equations. This is described in section~\ref{sec:discretisation_ns}.
The derivatives are computed by solving the adjoint equations, described in section~\ref{sec:adjoint_equations}.

\subsubsection{Discretisation of the Navier-Stokes equations}\label{sec:discretisation_ns}
The Navier-Stokes equations are discretised with a $\theta$
time-stepping scheme and the finite element method. The controlled Dirichlet
boundary conditions are weakly enforced with a variant of the Nitsche method \cite{nitsche1971}. An advantage of the
Nitsche approach is that the boundary values are explicitly included in the variational
formulation, which simplifies the (automated) derivation of the adjoint equations.
This is exploited in the implementation.

For the spatial discretisation, we consider conforming finite element spaces
\begin{equation}\label{eq:space_discretisation}
  \begin{aligned}
    V_h &\subset H^1_{0,\walls}(\Omega) = \set{u \in H^1(\Omega)\, :  u|_\walls = 0 }
    \\
    Q_h &\subset L^2(\Omega).
  \end{aligned}
\end{equation}
For the time discretisation we assume a partition of the interval $[0,T]$ with a constant timestep $\delta t$.  Applying a standard $\theta$
time-stepping scheme to the Navier-Stokes equations
\eqref{eq:navier_stokes}, we obtain a sequence of nonlinear problems:
For $k = 0, \ldots, N-1$, let $u^{k+\theta}= \theta u^{k+1} +
(1-\theta)u^{k}$ and find $(u^{k+1},p^{k+1})\in V_h \times Q_h$
such that
\begin{equation}\label{eq:ns_semi-discrete}
\begin{aligned}
  \frac{u^{k+1}-u^k}{\delta t} -\nu \laplace u^{k+\theta} +(u^{k+\theta}\cdot\del)u^{k+\theta} -\grad p^{k+1} &= 0,
  \\
  \div u^{k+1} &= 0,
\end{aligned}
\end{equation}
subject to the boundary conditions \eqref{eq:ns_bc}.  The equations
\eqref{eq:ns_semi-discrete} are integrated against test functions $v\in
V_h$ and $q\in Q_h$ in order to obtain a nonlinear variational
problem at each time $t_k$,
\begin{equation}\label{eq:nonlinear_var_problem_k}
  \begin{aligned}
    0
    & =  \int_\Omega \left(\frac{u^{k+1}-u^k}{\delta t}\right)\cdot v
      +  \nu \grad u^{k+\theta} \colon \grad v\, dx
    \\ & \qquad
      + \int_\Omega (u^{k+\theta}\cdot\del)u^{k+\theta}\cdot v\,dx
    \\ & \qquad
      + \int_\Omega  q \div u^{k+1} +    p^{k+1} \div v \, dx
    \\ & \qquad
      - \int_{\Gamma_D} \left(\nu \frac{\partial u^{k+\theta}}{\partial n}
      -  p^{k+1} n\right) \cdot v\, ds
    \\ & \qquad
      - \int_{\Gamma_D} \left(\theta \nu \frac{\partial v}{\partial n}
      -  q n\right) \cdot (u^{k+1}-g^{k+1})\, ds
    \\ & \qquad
      + \int_{\Gamma_D} \frac{\nu\sigma}{h} (u^{k+1}-g^{k+1})\cdot v \, ds.
  \end{aligned}
\end{equation}
The nonlinear variational problem
\eqref{eq:nonlinear_var_problem_k} consists of a volume integral and a
boundary integral over $\Gamma_D$.  The volume integral coincides
with the ``standard'' variational form of \eqref{eq:ns_semi-discrete}
obtained when the boundary condition \eqref{eq:ns_bc_dirichlet} is
strongly imposed. The second, boundary integral part of the
variational problem arises from the weakly imposing the Dirichlet
boundary condition \eqref{eq:ns_boundary_form} with Nitsche's method,
and is discussed in detail below.

The discrete spaces for the state and control variables are
\begin{align*}
  Y &= V_h^N\times Q_h^N
  \\
  M &= V_h \times (T_{\Gamma}V_h)^N,
\end{align*}
and we introduce the notation
\begin{align*}
  y &= (u, p) \in Y,
  \\
  u &= (u_1, \ldots, v_N) \in V_h^N
  \\
  p &= (p_1, \ldots, q_N) \in Q_h^N
  \\
  m &= (u_0, g_1, \ldots, g_N) \in M
\end{align*}
The sequence of variational problems \eqref{eq:nonlinear_var_problem_k}
is reformulated as an operator equation combining all the time steps,
\begin{equation}\label{eq:ns_operator_form}
  \mathcal{F}(m, y) = \sum_{k=0}^{N-1}\big\{ \mathcal{F}_{k,\Omega}(m, y) + \mathcal{F}_{k,\Gamma_D}(m, y)\big\} = 0,
\end{equation}
where $\mathcal{F}_{k,\Omega}:Y\rightarrow Y^*$ 
is the operator combining all the
volume integrals in \eqref{eq:nonlinear_var_problem_k}, i.e.
\begin{equation}\label{eq:ns_volume_form}
\begin{aligned}
  \mathcal{F}_{k,\Omega}(m, y; v, q)
  &=  \int_\Omega \left(\frac{u^{k+1}-u^k}{\delta t}\right)\cdot v^{k+1}
  +\nu \grad u^{k+\theta} \colon \grad v^{k+1}\, dx
  \\ & \qquad + \int_\Omega (u^{k+\theta}\cdot\del)u^{k+\theta}\cdot v^{k+1}\,dx
  \\ & \qquad + \int_\Omega  q^{k+1} \div u^{k+1} +    p^{k+1} \div v^{k+1} \, dx ,
\end{aligned}
\end{equation}
for all $(v,q) = \set{(v_k, q_k)}_{k=1}^{N} \in Y^*$, for $k=0, \ldots,  N-1$. Note that that this part only involves the initial data $u_0$ from $m$.
The operator $\mathcal{F}_{k,\Gamma_D}:Y\times M \rightarrow Y^*$
combines all the boundary integrals in \eqref{eq:nonlinear_var_problem_k} and reads
\begin{equation}\label{eq:ns_boundary_form}
  \begin{aligned}
      \mathcal{F}_{k, \Gamma_D}(m, y; v, q)
      & = - \int_{\Gamma_D} \left(\nu \frac{\partial u^{k+\theta}}{\partial n} -  p^{k+1} n\right) \cdot v^{k+1}\, ds
      \\
      & \qquad - \int_{\Gamma_D} \left(\theta \nu \frac{\partial v^{k+1}}{\partial n} -  q^{k+1} n\right) \cdot (u^{k+1}-g^{k+1})\, ds
      \\
      & \qquad + \int_{\Gamma_D} \frac{\nu\sigma}{h} (u^{k+1}-g^{k+1})\cdot v^{k+1} \, ds.
  \end{aligned}
\end{equation}
The first integral in \eqref{eq:ns_boundary_form} arises when the
integration by parts formula is applied to \eqref{eq:ns_semi-discrete}, and the integral would vanish if the Dirichlet boundary condition
\eqref{eq:ns_bc_dirichlet} were strongly imposed on the space $V_h$.
The remaining terms are added to obtain a variational problem that is
consistent and stable, see e.g. \cite{becker2009nitsche, burman2007}.
The form \eqref{eq:ns_boundary_form} is linear
and symmetric, and positive definite provided that the parameter
$\sigma$ is sufficiently large. We must also require $\theta > 0$ to
apply the Nitsche method.


The numerical examples in section \ref{sec:examples} use two common finite element
pairs: P2-P1 (Taylor-Hood) and P1-P1.  The lowest order discretisation
does not satisfy the LBB conditions, and hence requires stabilisation. We
used the stabilisation $-\beta h^2 \left(\nabla p, \nabla
q\right)_\Omega$ where $h$ is the local mesh element size and $\beta = 10^{-3}$
is the stabilisation coefficient.


\subsubsection{Adjoint equations}\label{sec:adjoint_equations}
The adjoint equations are used to efficiently compute the functional derivative
$\textrm{d}J/\textrm{d}m: Y \times M \to M^*$, at a cost
of roughly one linearised Navier-Stokes solve.

To derive the adjoint equations consider the Navier-Stokes equations
in the operator form $\mathcal F(m; y) =0 \in Y^*$ and a functional $J(y, m) \in \mathbb R$.
The total derivative of the functional in direction $\tilde m$ is
\begin{equation}
    \left<\frac{\textrm{d}J}{\textrm{d}m}, \tilde m\right>_{M^*, M} = \left< \frac{\partial J}{\partial y}, \frac{\textrm{d}y}{\textrm{d}m} \tilde m\right>_{Y^*, Y} + \left<\frac{\partial J}{\partial m}, \tilde m\right>_{M^*, M}.
\label{eq:functional_derivative}
\end{equation}
Evaluating \eqref{eq:functional_derivative} directly is challenging because computing
${\textrm{d}y}/{\textrm{d}m}(m) \in \mathcal L(M, Y)$ is computationally expensive. The adjoint
approach eliminates this term by taking the derivative of the PDE equation
\begin{equation}
    \frac{\partial \mathcal F}{\partial y} \frac{\textrm{d}y}{\textrm{d}m} + \frac{\partial \mathcal F}{\partial m} = 0.
\end{equation}
and substituting it into \eqref{eq:functional_derivative}:
\begin{equation}
    \left<\frac{\textrm{d}J}{\textrm{d}m}, \tilde m\right>_{M^*, M} = - \left<\frac{\partial \mathcal F}{\partial m}  \tilde m, \left(\frac{\partial \mathcal F}{\partial y}\right)^{-*}\frac{\partial J}{\partial y}\right>_{Y^*, Y} + \left<\frac{\partial J}{\partial m}, \tilde m\right>_{M^*, M}.
\end{equation}
The functional derivative is then computed in two steps:
\begin{enumerate}
    \item  Compute the adjoint solution $\lambda \in Y$ by solving the adjoint PDE
        \begin{equation}\label{eq:adjoint_eqn}
    \left(\frac{\partial \mathcal F}{\partial y}\right)^* \lambda =  -\frac{\partial J}{\partial y}
\end{equation}
\item Evaluate the derivative with
    \begin{equation}\label{eq:gradient_evaluation}
    \frac{\textrm{d}J}{\textrm{d}m} = \left(\frac{\partial \mathcal F}{\partial m} \right)^* \lambda + \frac{\partial J}{\partial m}
\end{equation}
\end{enumerate}
The computational expensive part is the solution of \eqref{eq:adjoint_eqn}, which involves the solution of a linear PDE.

The adjoint equations \eqref{eq:adjoint_eqn} can be derived before or after the discretisation of the
Navier-Stokes equations.
Here, we chose the discretise-then-adjoint approach,
which has the advantage that the discretised derivative is the exact derivative
of the discretised system. The alternative approach does not guarantee this, and simple descent methods like algorithm \ref{alg:descent} may fail, as demonstrated in~\cite{gunzburger2002}. As a
consequence a more robust optimisation algorithm would need to be implemented.

The adjoint system \eqref{eq:adjoint_eqn} for the discretised Navier-Stokes operator
\eqref{eq:ns_operator_form} is
\begin{equation}\label{eq:adjoint_equation}
  \begin{split}
      \left \langle \left(
        \frac{\partial \mathcal{F}}{\partial y}
      \right)^* \lambda, w \right\rangle
      = \left \langle \left(\frac{\partial \mathcal{F}}{\partial y } \right)w, \lambda \right\rangle
      = - \left\langle \frac{\partial J}{\partial u}, 
        v\right\rangle,
  \end{split}
\end{equation}
for all $w = (v,q)\in Y$. 
Note that the derivative of the regularisation term in the functional vanishes because it does not depend on the state.
Since the adjoint operator is linear, 
it can be written in matrix form:
\begin{equation}\label{eq:adjoint_matrix_form}
    \left(\frac{\partial \mathcal F}{\partial y} \right)^* =
\begin{pmatrix}
    \frac{\partial \mathcal F_{0}}{\partial y^{1}} & 0 & 0& \cdots \\
    \frac{\partial \mathcal F_{1}}{\partial y^{1}} & \frac{\partial \mathcal F_{1}}{\partial y^{2}} & 0 & \ddots \\
    0& \frac{\partial \mathcal F_{2}}{\partial y^{2}} & \ddots & \ddots
\end{pmatrix}^*
=
\begin{pmatrix}
    \frac{\partial \mathcal F^*_{0}}{\partial y^{1}} &  \frac{\partial \mathcal F^*_{1}}{\partial y^{1}} & 0& \cdots \\
   0 & \frac{\partial \mathcal F^*_{1}}{\partial y^{2}} &  \frac{\partial \mathcal F^*_{2}}{\partial y^{2}} & \ddots \\
    0&0 & \ddots & \ddots
\end{pmatrix}
\end{equation}
or more compactly, 
\begin{equation*}
  \left(\frac{\partial \mathcal F}{\partial y^k} \right)^* \lambda
  = \begin{cases}
    \left(\frac{\partial \mathcal F_{k-1}}{\partial y^k} \right)^*\lambda_{k}
  + \left(\frac{\partial \mathcal F_{k\vphantom{-1}}}{\partial y^k} \right)^*\lambda_{k+1}
    \vphantom{\vdots}
  & \mbox{if } k < N \\
    \left(\frac{\partial \mathcal F_{k-1}}{\partial y^k} \right)^*\lambda_k
    \vphantom{\vdots}
  & \mbox{if } k = N. \\
  \end{cases}
\end{equation*}
The system \eqref{eq:adjoint_matrix_form} is upper-triangular, hence the adjoint \eqref{eq:adjoint_equation} is solved by backwards substitution.
Written explicitly, the volume integrals in the equation for $\lambda_k$, $k=0,\ldots, N$, are
\begin{equation}
\begin{aligned}
    & \int_\Omega \left(\frac{\lambda^{k} - \lambda^{k+1}}{\delta t}\right)\cdot v^k\, dx
    + \int_\Omega \nu \grad \lambda^{k+\tilde \theta} \colon \grad v^k\, dx
    \\
    + &\theta \int_\Omega (u^{k-1+\theta}\cdot\del)v^k \cdot \lambda_u^{k}\,dx
    + \theta \int_\Omega (v^k\cdot\del)u^{k-1+\theta}\cdot \lambda_u^{k}\,dx
    \\
    + &\tilde \theta \int_\Omega (u^{k+\theta}\cdot\del)v^k \cdot \lambda_u^{k+1}\,dx
    + \tilde \theta \int_\Omega (v^k\cdot\del)u^{k+\theta}\cdot \lambda_u^{k+1}\,dx
  \\  + &\int_\Omega  \lambda_p^{k} \div v^k +  q^k \div \lambda_u^{k} \, dx,
\end{aligned}
\end{equation}
with $(\lambda^k_u, \lambda^k_p) = \lambda^k$,  $\tilde \theta = 1 - \theta$ and setting $\lambda_u^{N+1}=0$.
Similarly, the boundary integrals are
\begin{equation}
  \begin{aligned}
      - & \int_{\Gamma_D} \left(\nu \frac{\partial v^k}{\partial n} \right) \cdot \lambda_u^{k+\tilde \theta}\, ds
       + \int_{\Gamma_D} \left(q^k n\right) \cdot \lambda_u^{k}\, ds
      \\
      - & \int_{\Gamma_D} \left(\theta \nu \frac{\partial \lambda_u^{k}}{\partial n} -  \lambda_p^{k} n\right) \cdot v^k\, ds
       + \int_{\Gamma_D} \frac{\nu\sigma}{h} v^k\cdot \lambda_u^{k} \, ds.
  \end{aligned}
\end{equation}
The adjoint equations are solved backwards in time, starting from a zero final condition.
The timestepping scheme is the same $\theta$-scheme as for the forward discretisation, but with  a modified advective velocity.
The homogeneous Dirichlet boundary conditions on the controlled surfaces are enforced with a Nitsche like approach.


\subsection{Implementation and verification}
The Navier-Stokes solver was implemented in the FEniCS finite element
framework \cite{logg2012}. The adjoint solver was automatically derived via
the algorithmic differentiation tool dolfin-adjoint \cite{farrell2013}. The
correctness of the adjoint equations, and the resulting derivatives of the goal
functional, were verified using the Taylor remainder convergence test. This test checks that
for a sufficiently smooth functional $\hat J$, a correct
implementation should satisfy
\begin{equation}
    \left|\hat{J}(m + h\delta m) - \hat{J}(m) - h\left<\frac{\text{d} \hat{J}(m)}{\text{d} m}, \delta m\right>\right| = O(h^2),
\end{equation}
where $\delta m$ is a control perturbation and $h>0$ the perturbation size.
The Taylor remainder convergence test was performed for different $m$ and
random directions $\delta m$. Second order convergence was consistently
observed, giving confidence that the adjoint implementation is correct.

\section{Experiments}\label{sec:examples}
In this section, the data assimilation is applied to two different experiments. The
first experiment uses an aneurysm-like domain in 2D with known exact solution (section~\ref{sec:example_2d}).
The second experiment aims to reconstruct the flow conditions in a real geometry in 3D with observations from an 4D MRI scan (section~\ref{sec:example_3d}).

The implementation and files needed to reproduce the results of this section are available on bitbucket:
\url{https://bitbucket.org/biocomp/navier_stokes_data_assimilation}.
This website contains a Readme file with instructions for the
installation and how to reproduce the paper results.

\subsection{2D Aneurysm}\label{sec:example_2d}
This experiment tests the variational data assimilation under idealised
conditions where the blood flow to be reconstructed is known a priori. This is used
to study the robustness of the reconstruction against incomplete data (both in space and time), noise and the
the choice of the regularisation. We also compare two different types of observation operators.

\begin{figure}
    \centering
    \includegraphics[width=0.2\textwidth]{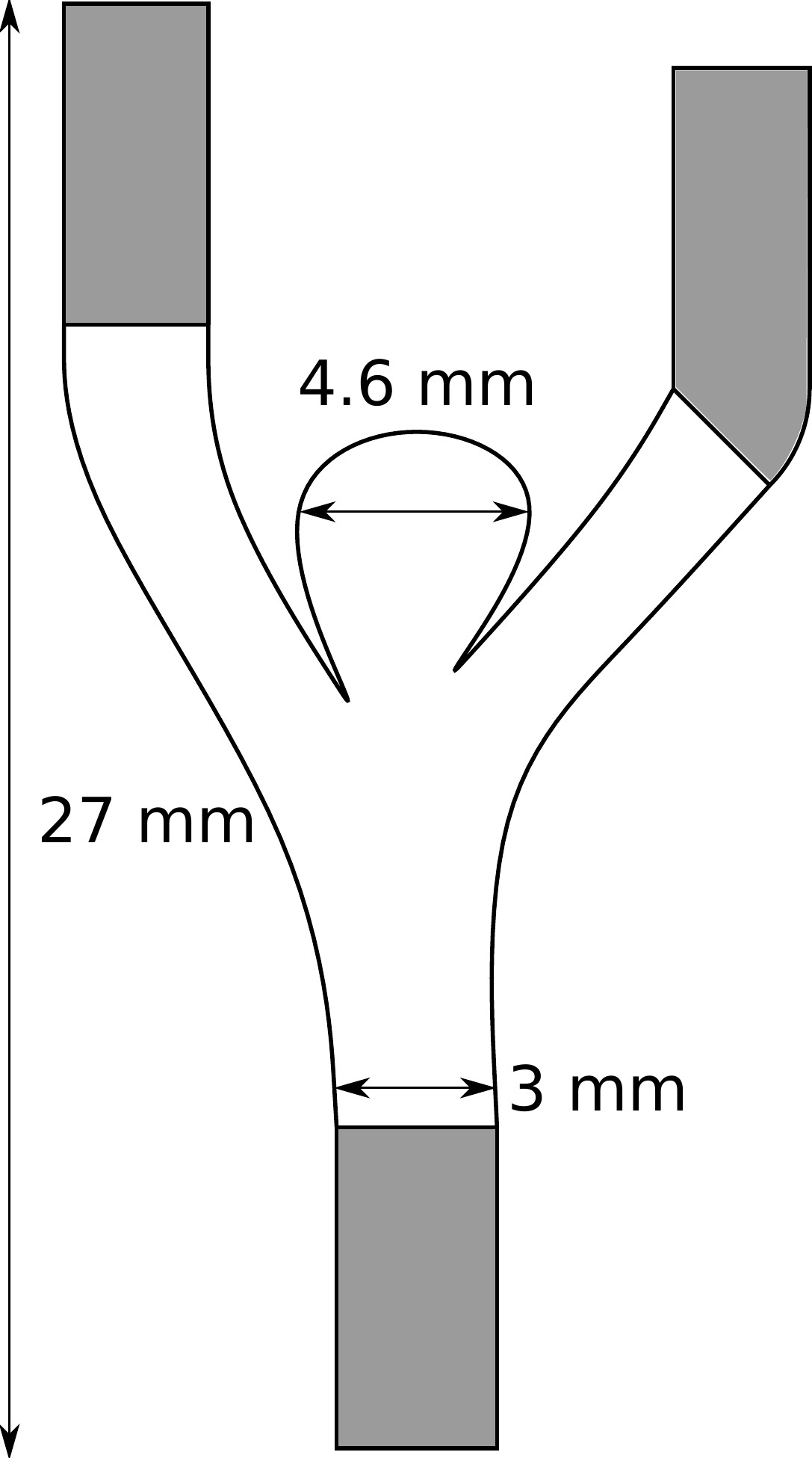}
    \caption{The computational domain for the 2D example. The grey area indicates
    the extended area used to generate the measurement data.}\label{img:domain-2D}
\end{figure}

\begin{table}
    \centering
    \begin{tabular}{ l c r }
        Parameter & Symbol & Value\\
        \hline
        Viscosity & $\nu$ & $3.5$ \\
        Model timestep  & $\Delta t$ & $0.004625$ s \\
        End time  & $T$ & $0.555$ s \\
        Time discretisation & $\theta$ & $0.5$ \\
        Spatial discretisation &  & $P1$-$P1$ \\
        Mesh triangles &  & $20,989$ \\
        Nitsche coefficient & $\sigma$ & $100.0$ \\
        \hline
        Number of observations & N & $16$ \\
        Regularisation & $\alpha = \gamma$ & $10^{-5}$\\
        \hline
    \end{tabular}
    \caption{The settings for the 2D aneurysm flow reconstruction. The first parameters specify the model setup, while the final two parameters configure the data assimilation.}
    \label{tab:2d_setup}
\end{table}

The computational domain, shown in figure \ref{img:domain-2D}, resembles a blood
vessel bifurcation with an aneurysm in 2D. The observations were generated with the same numerical model that was used in
the data assimilation procedure. That is, the Navier-Stokes
equations was first solved on an extended domain (including the gray area in
figure~\ref{img:domain-2D}) and the velocity solution used to generate the observations. For this setup, the initial velocity was set to zero. On the
inlet and right outlet boundaries a parabolic velocity profiles was enforced
with peak values of $1000$ mm/s (inlet) and $870$ mm/s (right outlet),
multiplied by $\sin(\pi (1-t)^3)$ to obtain a pulse like flow
pattern\footnote{Note that real flow in cerebral arteries will never go to
zero, but rather pulsate between around 0.5 m/s and one third of 0.5 m/s.}. The simulation started from a zero velocity at $t=0$ and
was terminated after the peak velocity at $t = 0.629$.
The remaining model settings are listed in
table~\ref{tab:2d_setup}.

The observation operator $\mathcal{T}$ was applied to the resulting velocity
to obtain $N=16$ observations. We compared two observation
operators: the \emph{instantaneous observation operator}, which takes
instantaneous measurements at evenly distributed times $t_n$:
\begin{equation}
    \mathcal{T}^{\text{inst}} u := \mathcal{R}_{\Omega_\text{obs}} u(t_n),
\end{equation}
where $\mathcal{R}_{\Omega_\text{obs}}$ restricts the velocity to
the observation domain (the white area in figure~\ref{img:domain-2D}). The
restriction avoids the ``inverse crime'' and simulates the incompleteness of real measurement data.
The \emph{time-averaging observation operator} takes pointwise time-averaged observations over each observation time interval:
\begin{equation}
    \mathcal{T}^{\text{avg}} u := \frac{1}{t_n-t_{n-1}}\int_{t_{n-1}}^{t_n}\mathcal{R}_{\Omega_\text{obs}} u(t)\text{d}t.
\end{equation}

The data assimilation was then applied to recover the
original flow from the observations. The reconstructions were performed on the restricted domain
$\Omega_\text{obs}$, that is without any knowledge about the geometry of the
extended domain. Furthermore, the outflow Dirichlet
boundaries were swapped between the data generation and reconstruction to further avoid the ``inverse crime''.
The optimisation was terminated when the relative change of the functional in one iteration
dropped below $|J(x_k) - J(x_{k-1})|/|J(x_0) \le 10^{-4}$ or if the number of iteration exceeded $100$.

The results of the data assimilation with $\mathcal{T}^{\text{inst}}$ and
$\mathcal{T}^{\text{avg}}$ are shown in figures
\ref{fig:results_aneurysm_inst_noise} and \ref{fig:results_aneurysm_avg_noise}
(left column), respectively.  The first three plots show the observed and
reconstructed velocities at $t=0.296$ s.  Visually, the observed and
assimilated velocities agree well.

Since the true velocity is known from the initial simulation, the reconstruction error can
also be quantified more rigorously.
We define following two error measures:
the first measures the relative error of the reconstructed velocity in the aneurysm
\begin{equation}
    \mathcal E_{\Omega_{\text{ane}}} = \frac{\|u_{true} - u \|_{\Omega_{\text{ane}}\times (0, T]}}{
                                      \|u_{true}     \|_{\Omega_{\text{ane}}\times (0, T]}},
                         \label{error_measure1}
\end{equation}
where $u_{true}$ is the true velocity. 
The second measures the relative error of the reconstructed wall shear stress on the aneurysm wall,
motivated by the fact that this an important diagnostic value in blood flow simulations:
\begin{equation}
    \mathcal E_{\text{WSS}} = \frac{\|\text{WSS}(u_{true})- \text{WSS}(u)\|_{\Gamma{\text{ane}}\times (0, T]}}{
                             \|\text{WSS}(u_{true})                       \|_{\Gamma{\text{ane}}\times (0, T]}},
                         \label{error_measure2}
\end{equation}
with $\text{WSS}(u) = |\sigma n - (\sigma n \cdot n)n|$ and $\sigma = \rho \left(-pI + \nu (\nabla u + (\nabla u)^T)\right)$ with $\rho = 1060$ $\text{kg}/\text{m}^3$.
The timeplots in figures \ref{fig:results_aneurysm_inst_noise} and \ref{fig:results_aneurysm_avg_noise}
(left column) visualise these error measures over the simulation period. The results
show a good agreement throughout the simulation period.

\subsubsection{Sensitivity of reconstruction with respect to parameter changes}

In this section, we investigate how
the quality of the reconstruction depends on noise in the observations and the choice of the reconstruction parameters,
such as the amount of regularisation.
Since the exact solution is known for this example, we can visualise the reconstructed and the ``true'' velocities
and compute the error measures \eqref{error_measure1} and \eqref{error_measure2}.
The following tests are based on the configuration listed in table
\ref{tab:2d_setup}, and in each test one parameter is varied and the
quality of the reconstruction investigated.

\begin{figure}
    \centering
        \footnotesize
    \hspace{1cm}
    \parbox[b][7.5em][t]{0.29\textwidth}{
        \emph{Base setup with no noise \\ {}}
        \vspace{0.5em}\\
        \input{results_aneurysm/nsassimilation/source/results_aneurysm/instant/assimilated_H1H1_0_noise/metrics}
    }
    \parbox[b][7.5em][t]{0.29\textwidth}{
        \emph{Base setup with\\signal-to-noise ratio of 2}
        \vspace{0.5em}\\
        \input{results_aneurysm/nsassimilation/source/results_aneurysm/instant/assimilated_H1H1_255_noise/metrics}
    }
    \parbox[b][7.5em][t]{0.29\textwidth}{
        \emph{Base setup with\\signal-to-noise ratio of 1}
        \vspace{0.5em}\\
        \input{results_aneurysm/nsassimilation/source/results_aneurysm/instant/assimilated_H1H1_360_noise/metrics}
    }
    \\
    \rotatebox{90}{\qquad \quad Observation}\quad
    \includegraphics[width=0.29\textwidth]{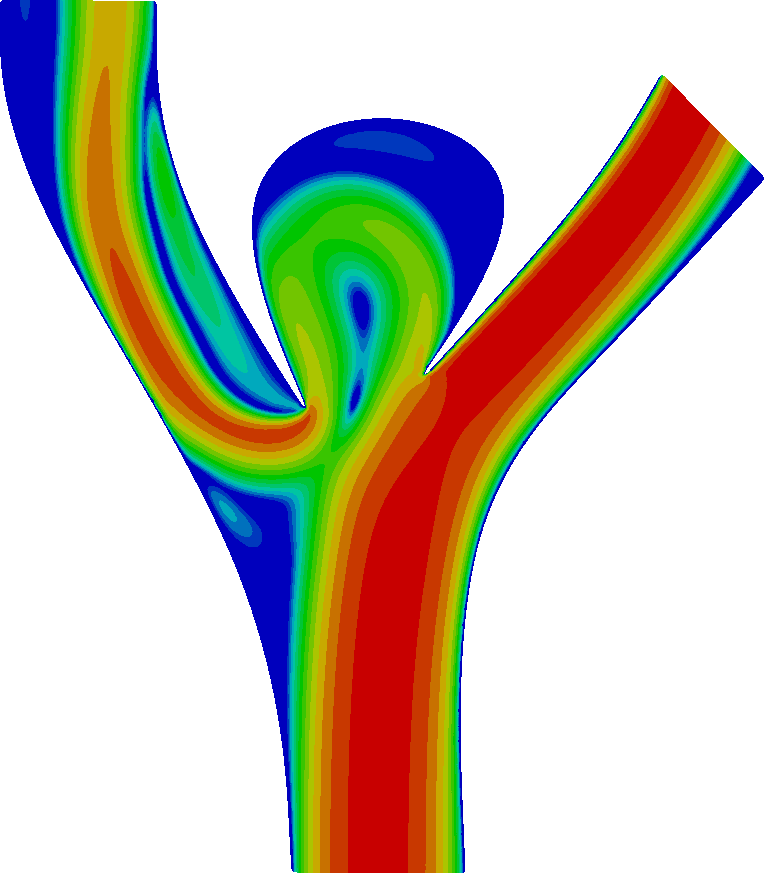}
    \includegraphics[width=0.29\textwidth]{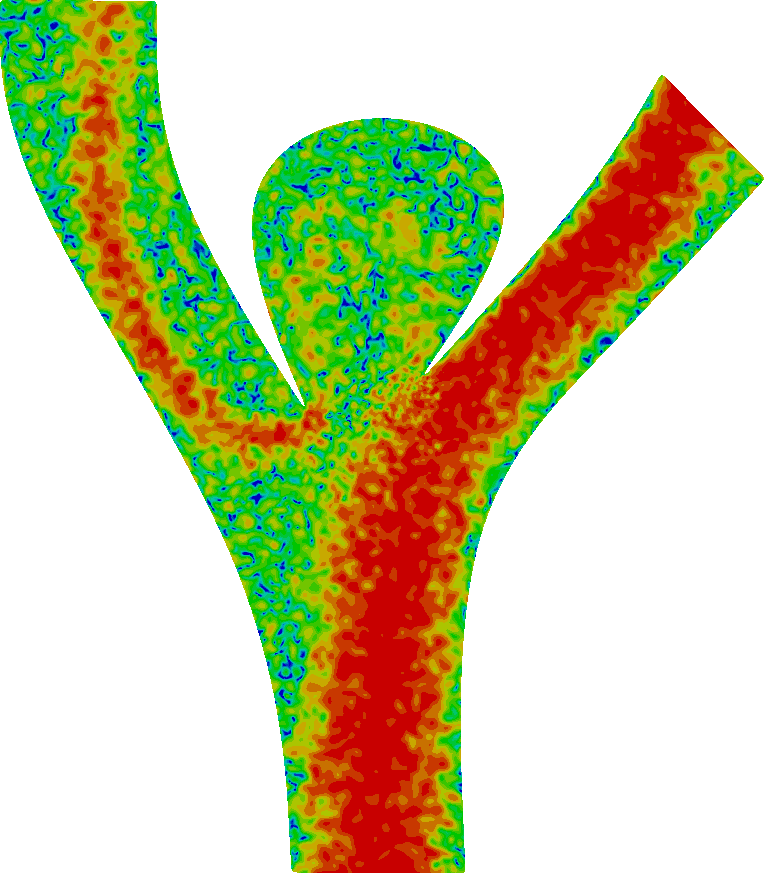}
    \includegraphics[width=0.29\textwidth]{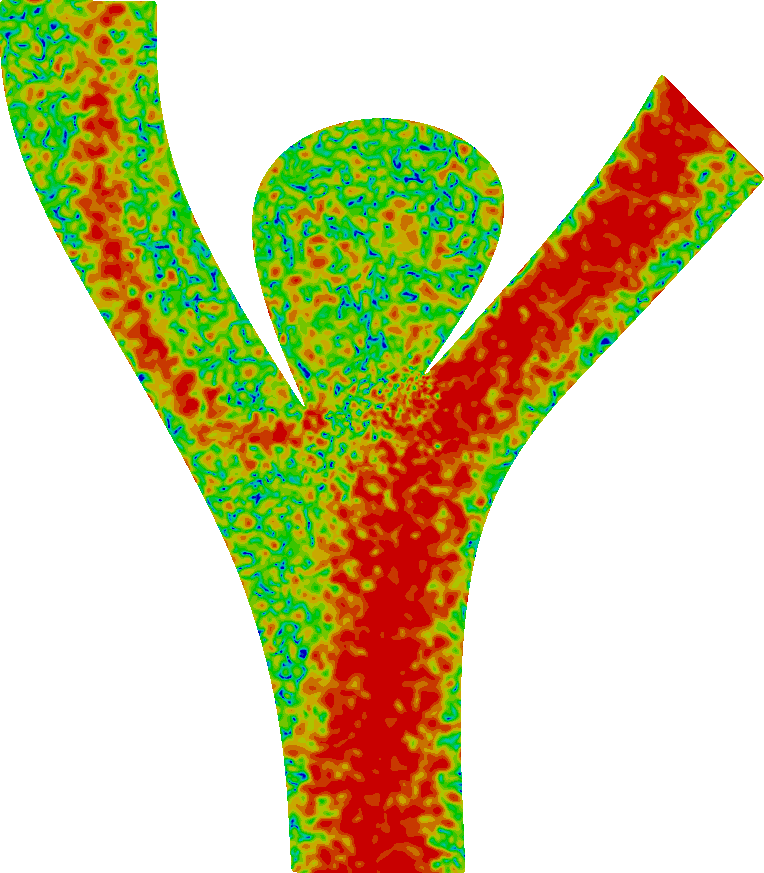} \\
    \rotatebox{90}{\qquad \quad Assimilation}\quad
    \includegraphics[width=0.29\textwidth]{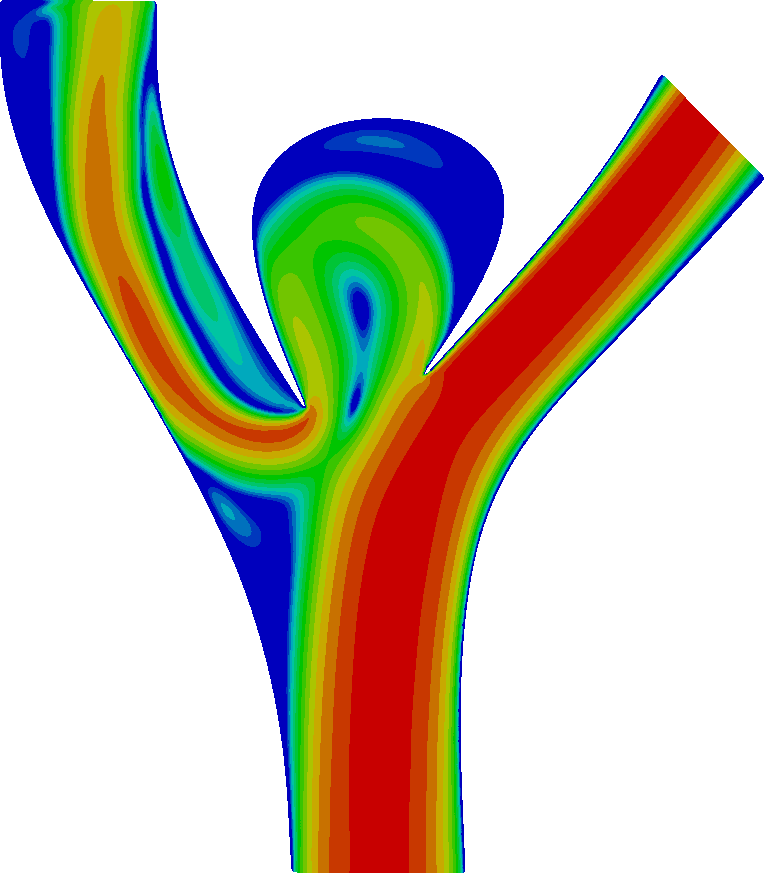}
    \includegraphics[width=0.29\textwidth]{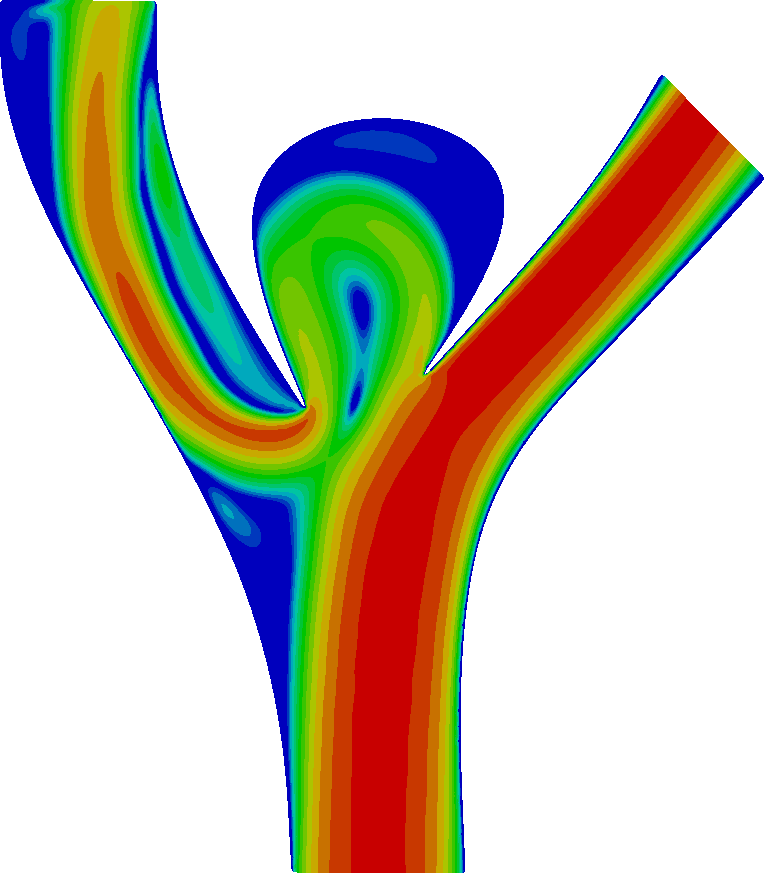}
    \includegraphics[width=0.29\textwidth]{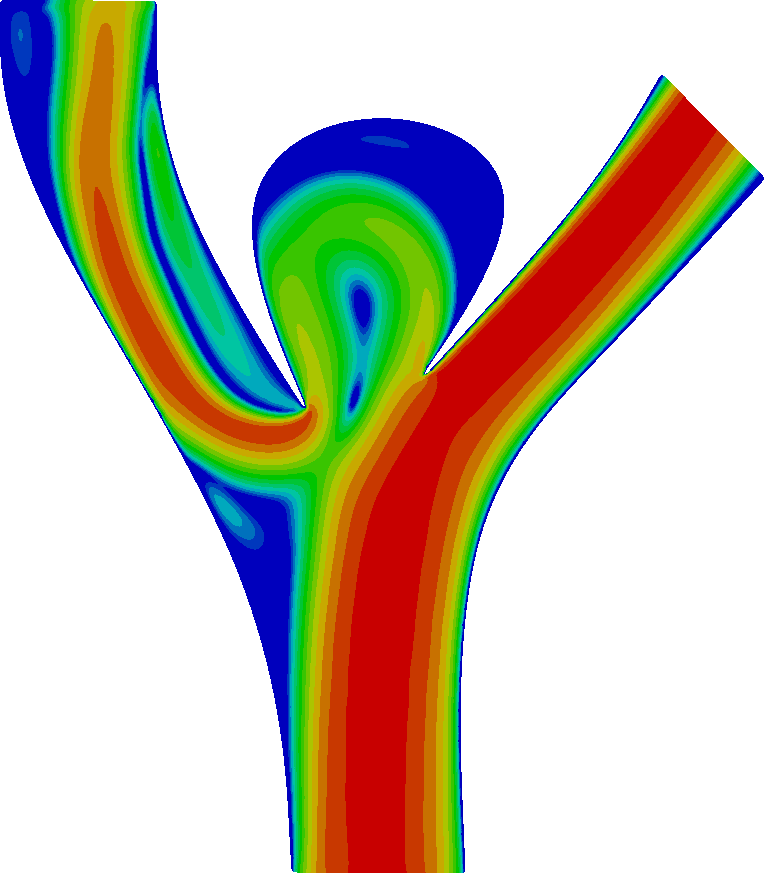} \\
    \begin{center}
        \includegraphics[width=0.3\textwidth]{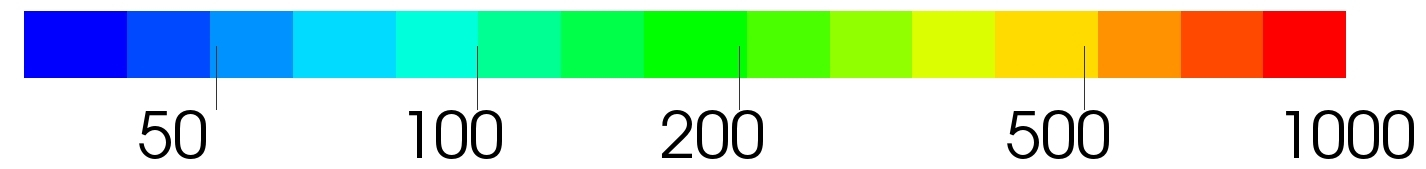}
        \vspace{-0.5cm}
    \end{center}
    \rotatebox{90}{Velocity norm in ${\Omega_{\text{ane}}}$}\quad
    \includegraphics[width=0.29\textwidth]{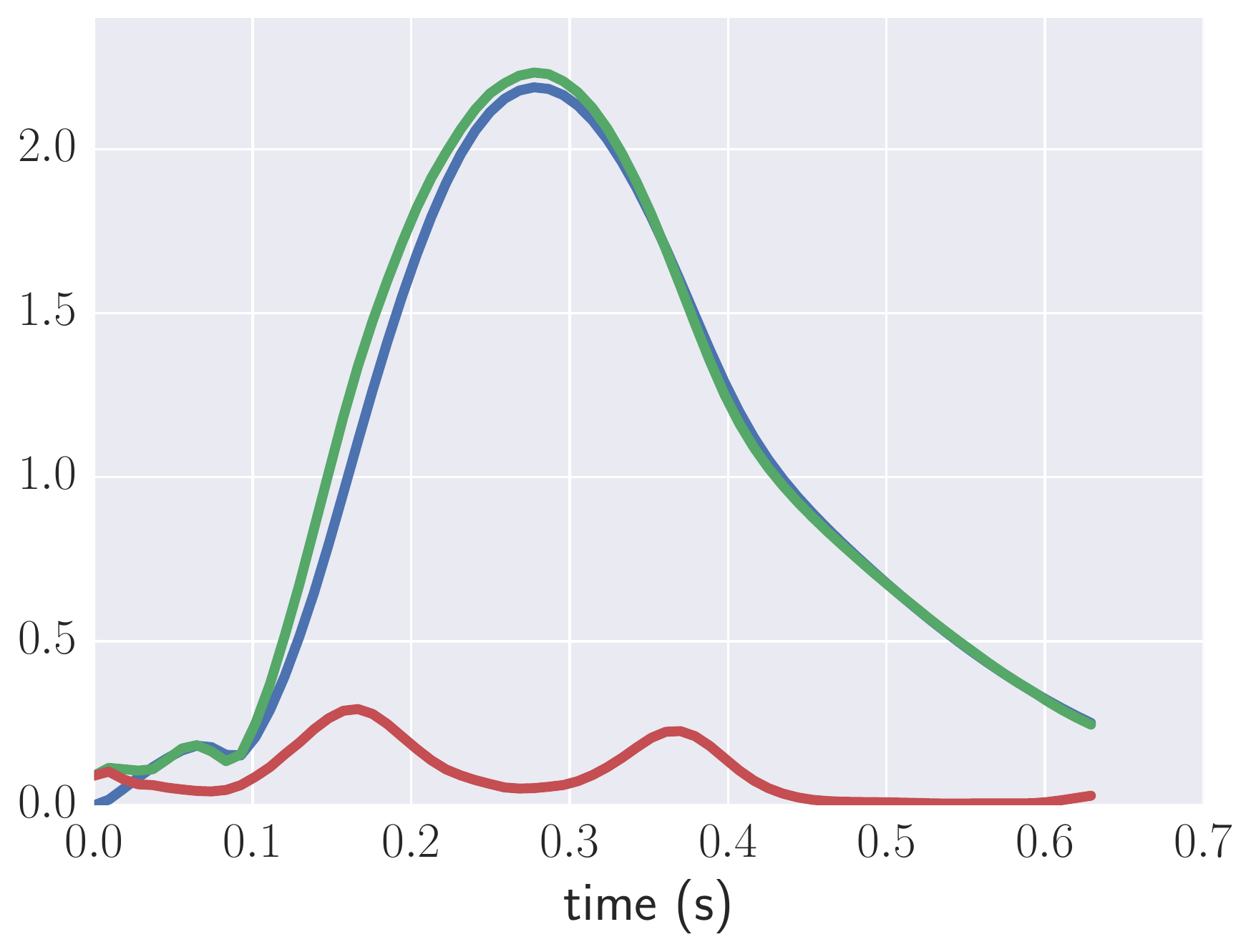}
    \includegraphics[width=0.29\textwidth]{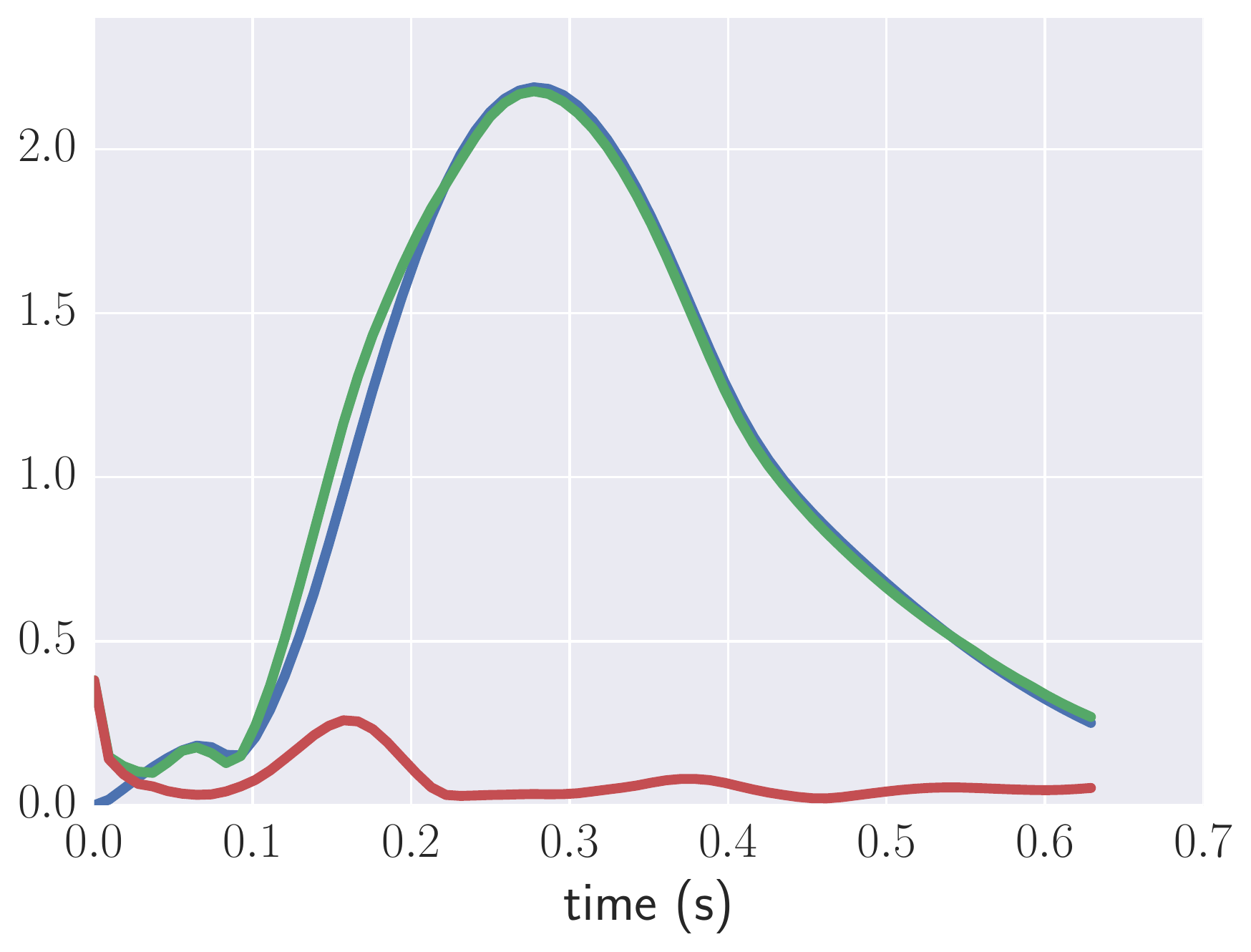}
    \includegraphics[width=0.29\textwidth]{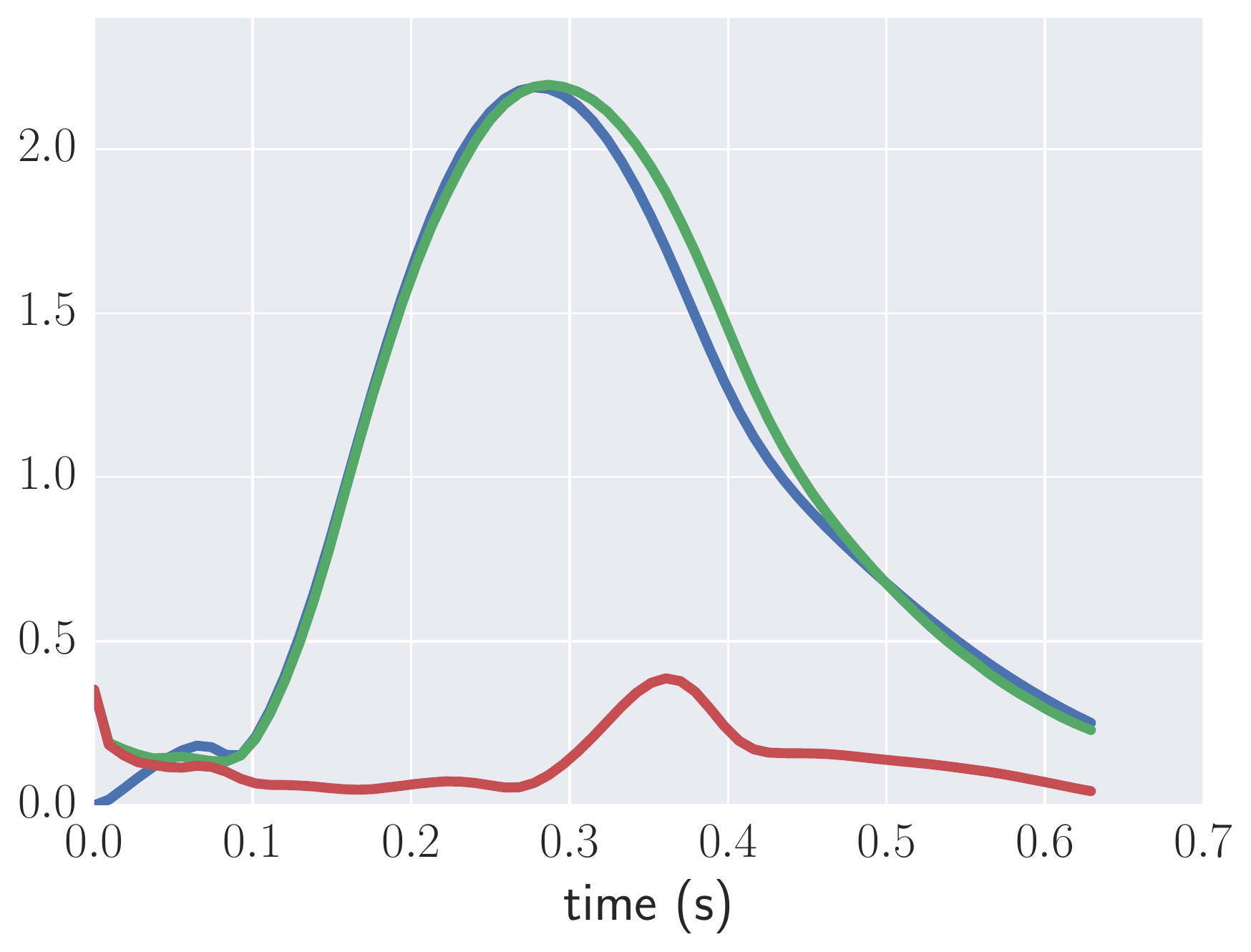} \\
    \rotatebox{90}{WSS norm on $\Gamma_{\text{ane}}$}\quad
    \includegraphics[width=0.29\textwidth]{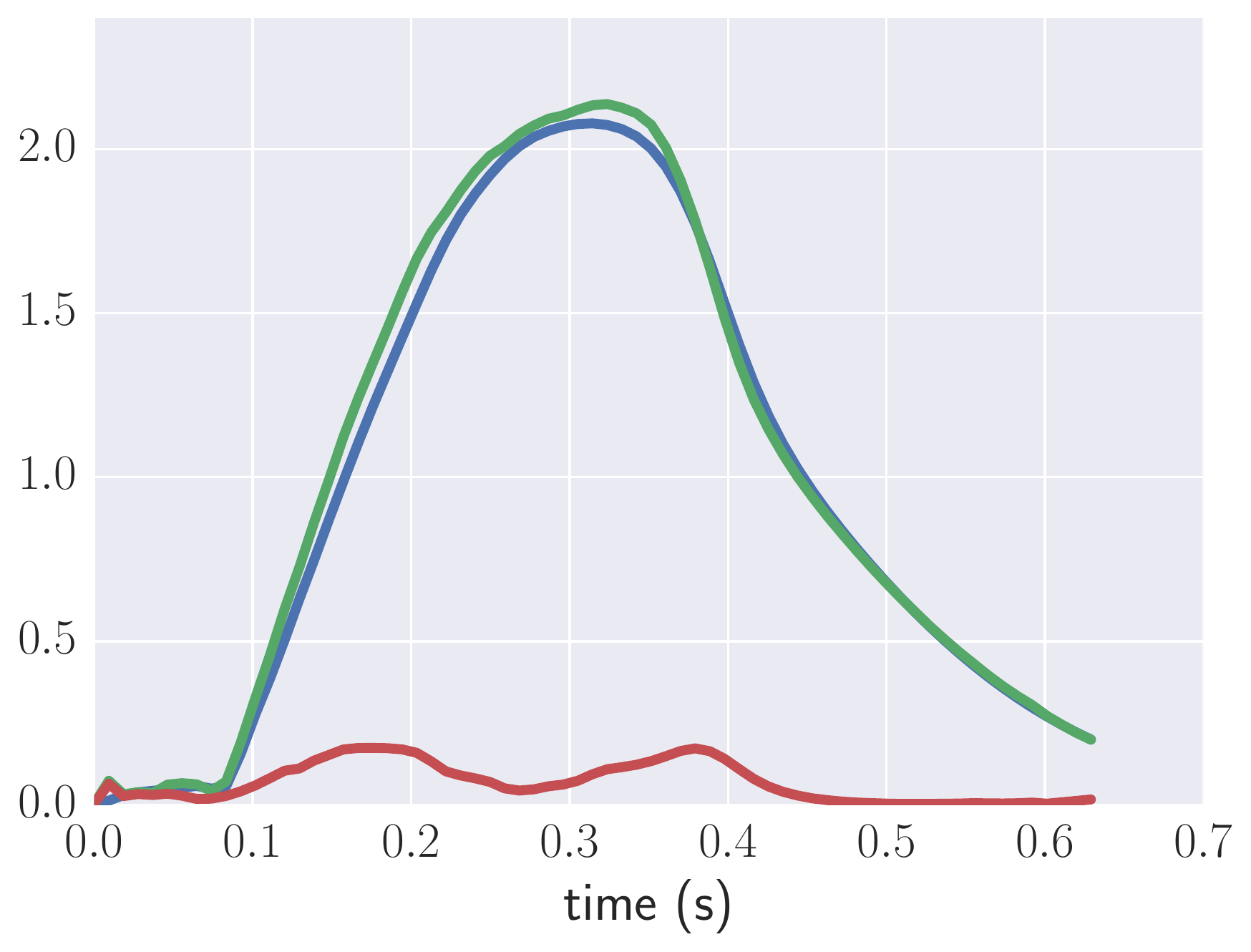}
    \includegraphics[width=0.29\textwidth]{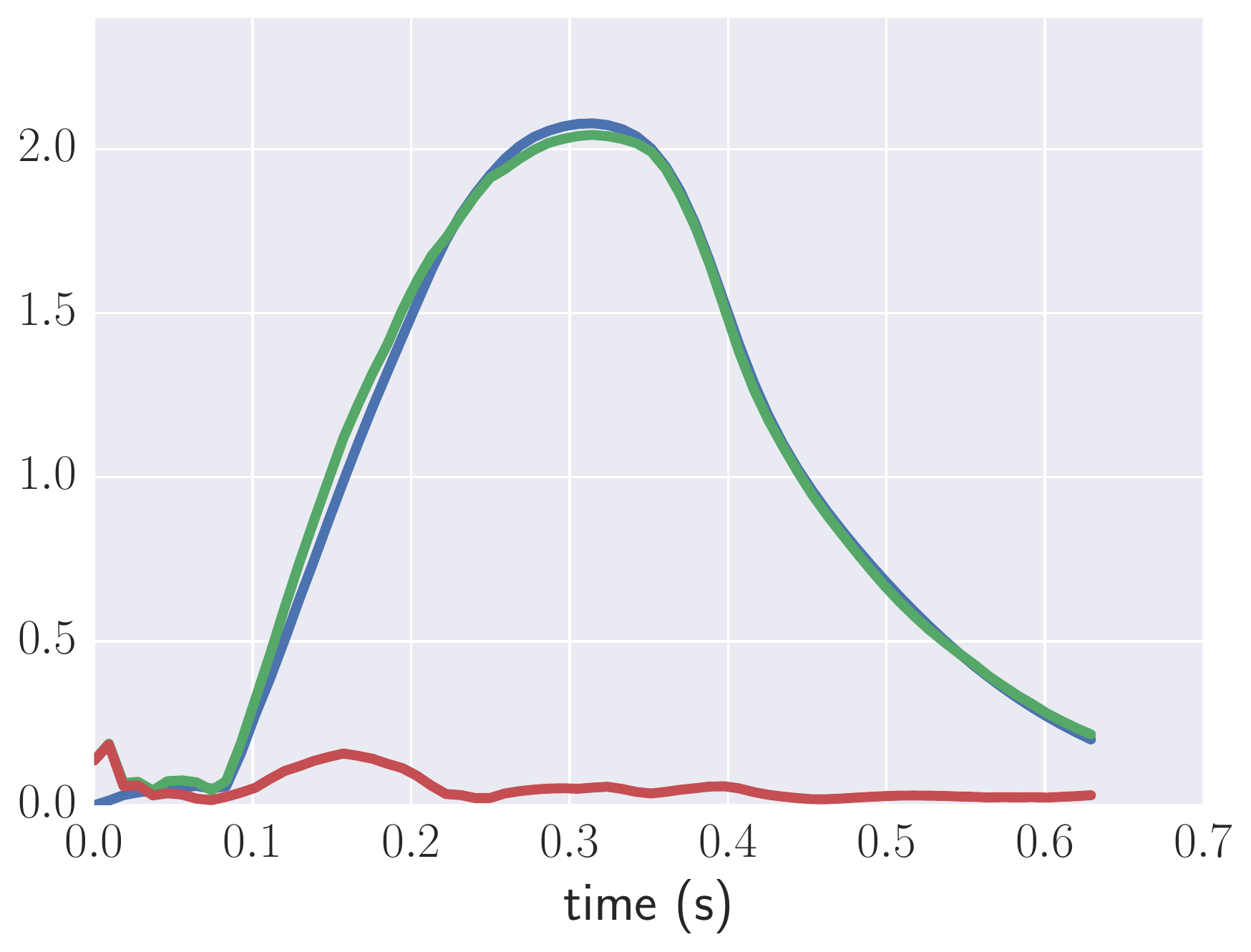}
    \includegraphics[width=0.29\textwidth]{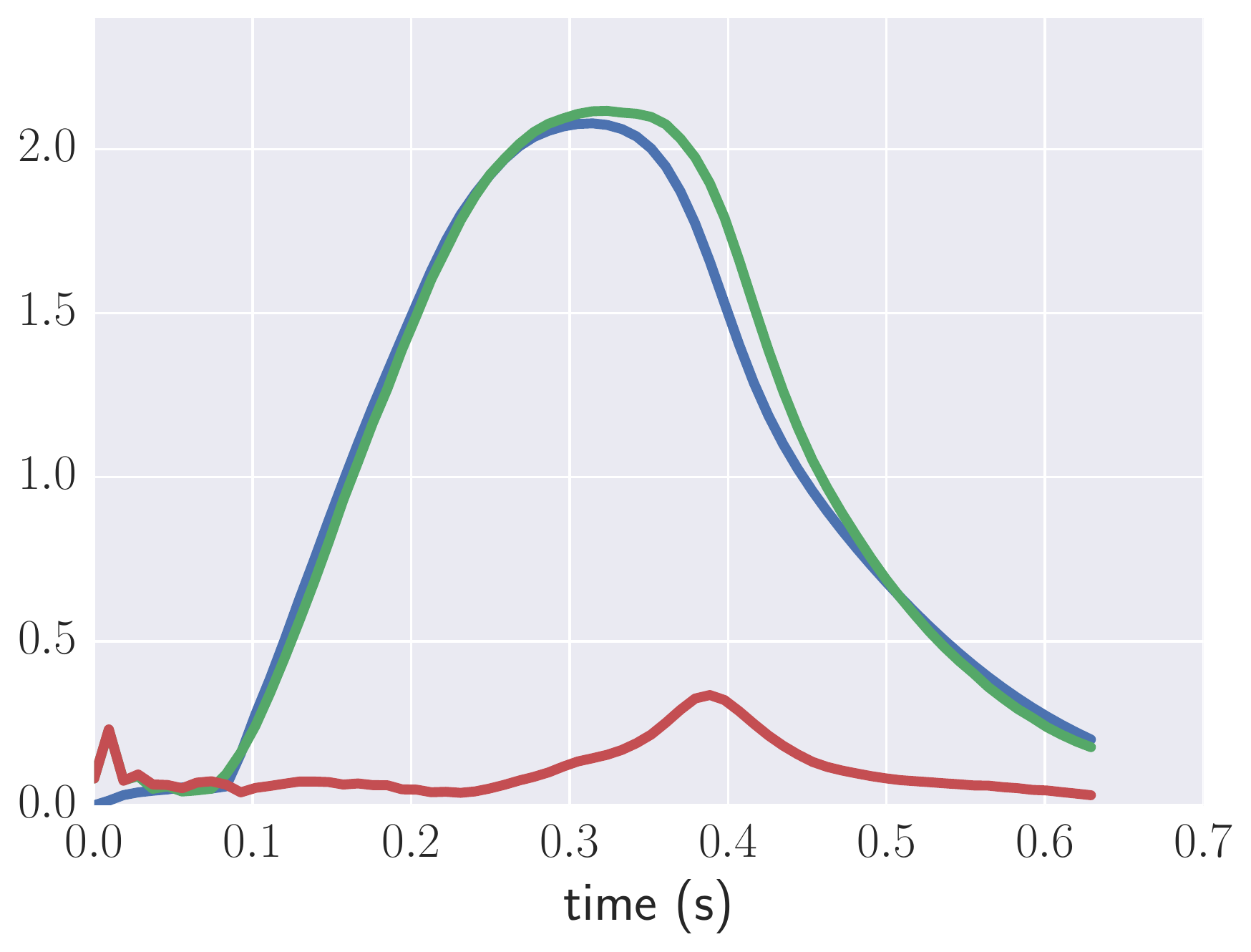} \\
    \includegraphics[width=0.59\textwidth]{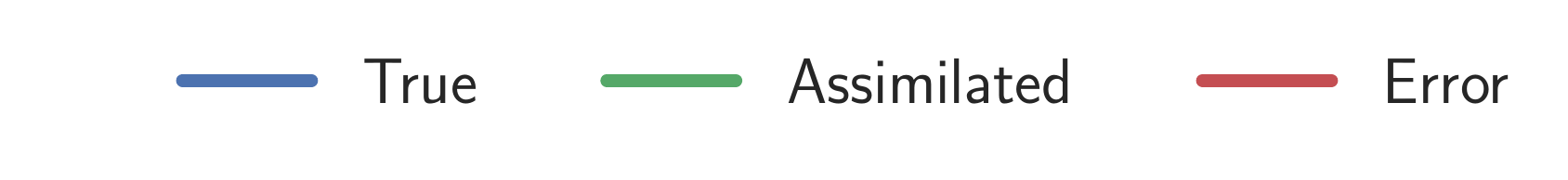}
    \caption{Results using the \textbf{instantaneous observation operator} with pointwise additive Gaussian white noise.
The signal-to-noise-ratio was computed as $\norm{\mathcal{T}^{\text{inst}}u_{\text{true}}}^2 / \norm{\mathcal{T}^{\text{inst}}u_{\text{true}}-d}^2$, where $d$ is the noisy data.
The snapshots on the top three rows are taken at $t=0.296$s.}
    \label{fig:results_aneurysm_inst_noise}
\end{figure}

\begin{figure}
    \centering
        \footnotesize
    \hspace{1cm}
    \parbox[b][7.5em][t]{0.29\textwidth}{
        \emph{Base setup with no noise\\{}}
        \vspace{0.5em}\\
        \input{results_aneurysm/nsassimilation/source/results_aneurysm/averaged/assimilated_H1H1_0_noise/metrics}
    }
    \parbox[b][7.5em][t]{0.29\textwidth}{
        \emph{Base setup with\\signal-to-noise ratio of 2}
        \vspace{0.5em}\\
        \input{results_aneurysm/nsassimilation/source/results_aneurysm/averaged/assimilated_H1H1_255_noise/metrics}
    }
    \parbox[b][7.5em][t]{0.29\textwidth}{
        \emph{Base setup with\\signal-to-noise ratio of 1}
        \vspace{0.5em}\\
        \input{results_aneurysm/nsassimilation/source/results_aneurysm/averaged/assimilated_H1H1_360_noise/metrics}
    }
    \\
    \rotatebox{90}{\qquad \quad Observation}\quad
    \includegraphics[width=0.29\textwidth]{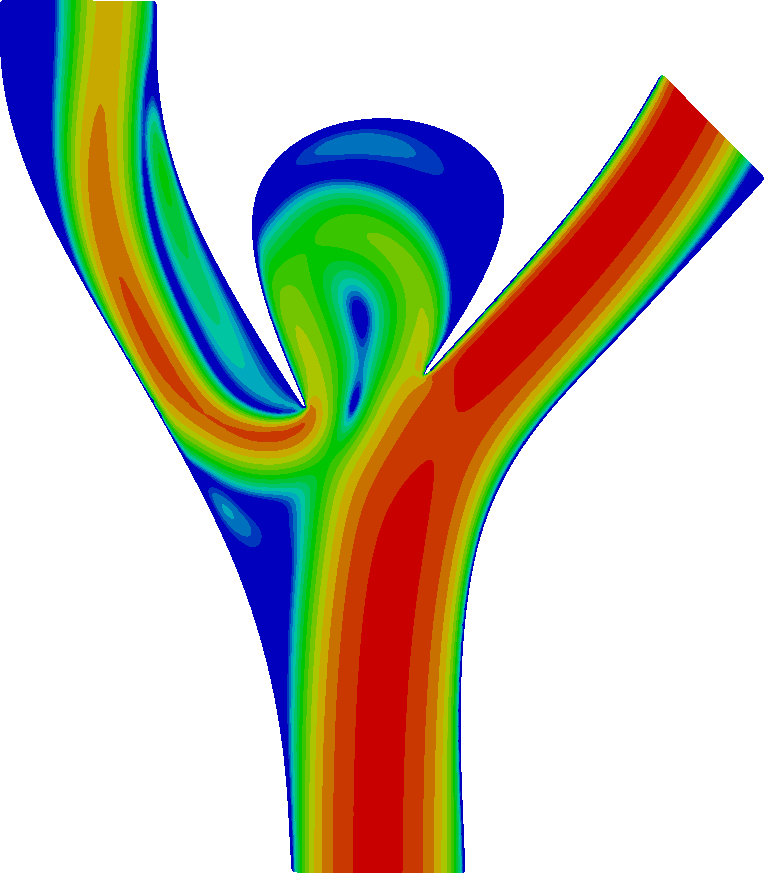}
    \includegraphics[width=0.29\textwidth]{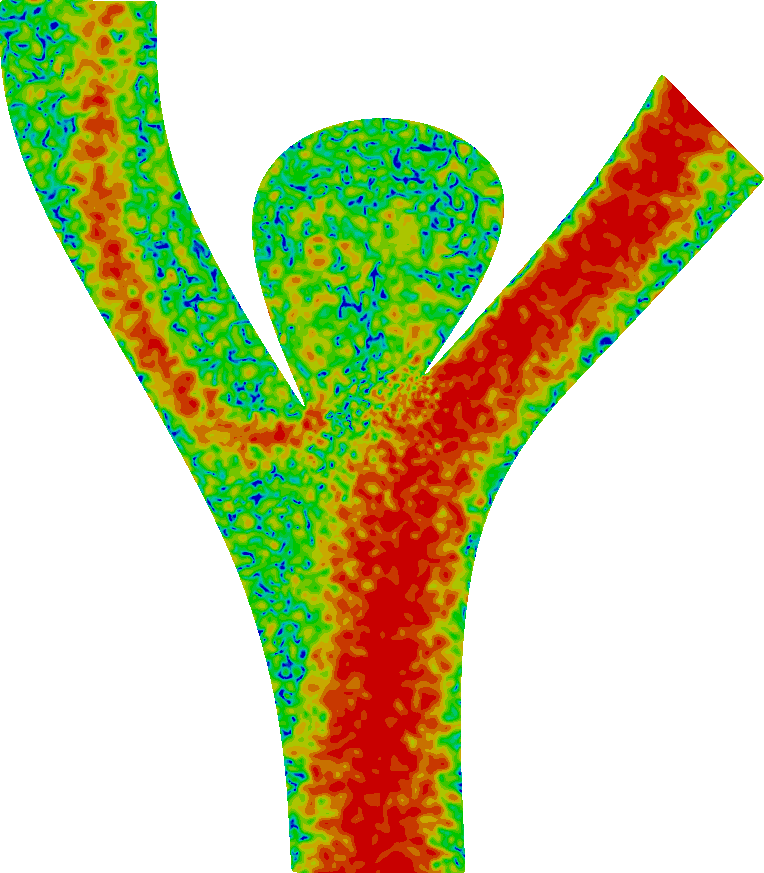}
    \includegraphics[width=0.29\textwidth]{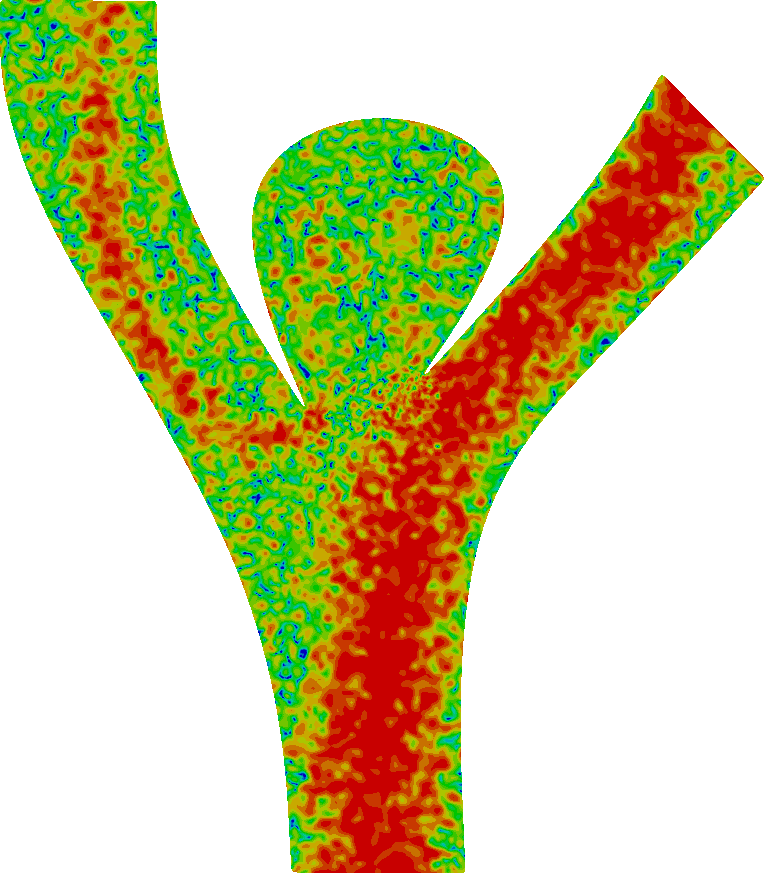} \\
    \rotatebox{90}{\qquad \quad Assimilation}\quad
    \includegraphics[width=0.29\textwidth]{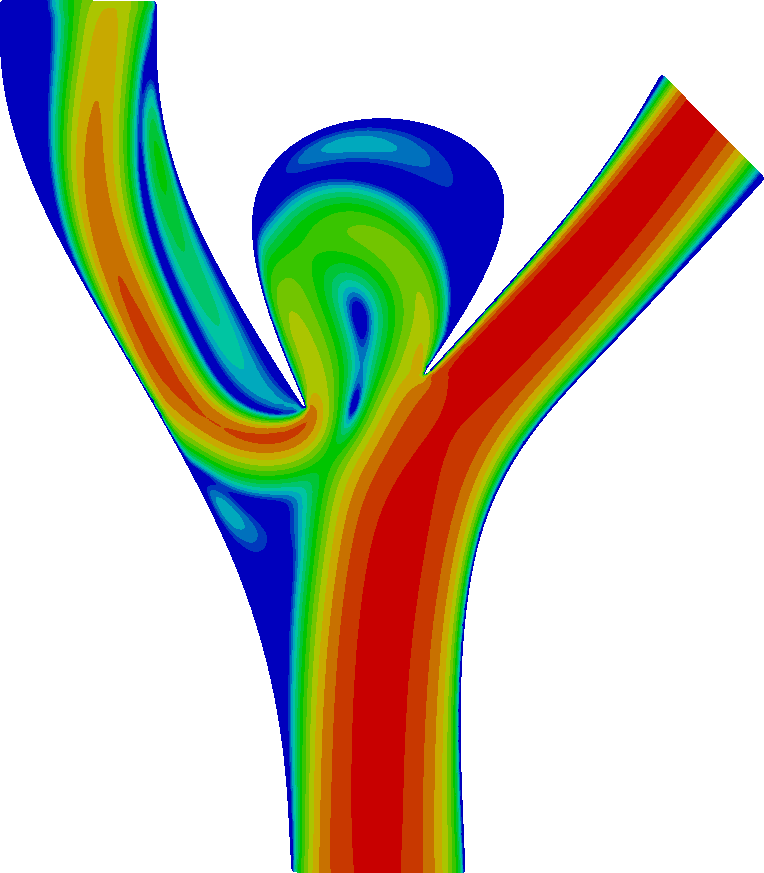}
    \includegraphics[width=0.29\textwidth]{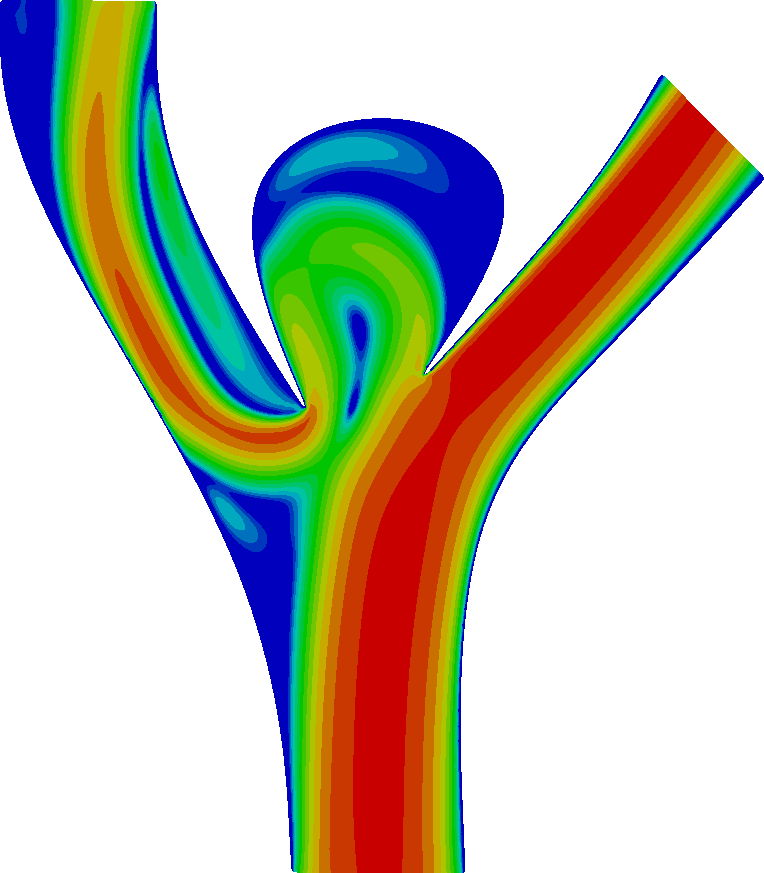}
    \includegraphics[width=0.29\textwidth]{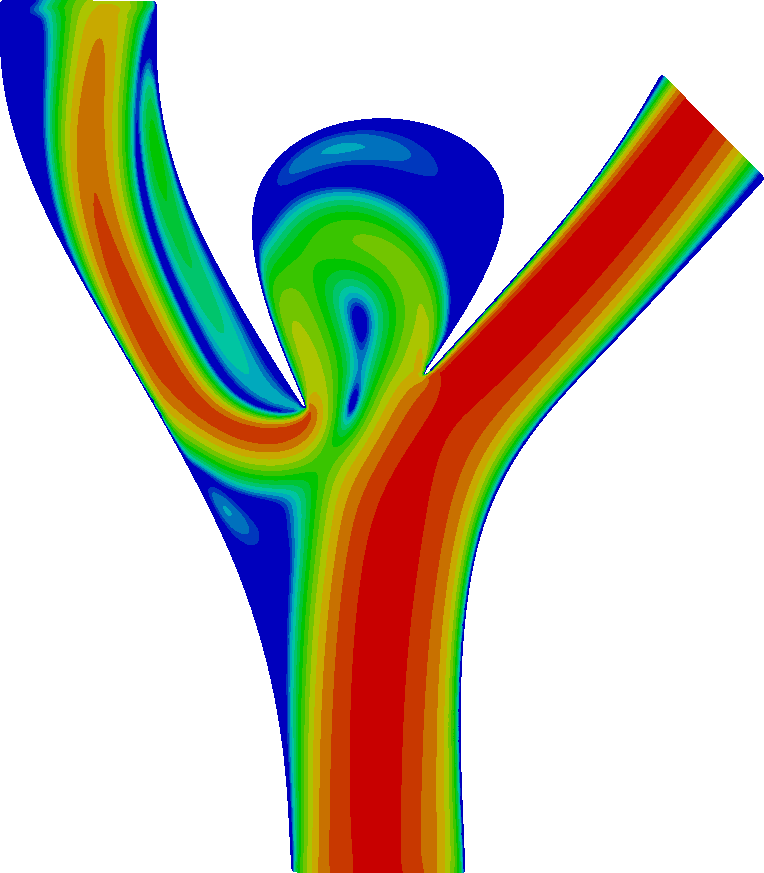} \\
    \begin{center}
        \includegraphics[width=0.3\textwidth]{results_aneurysm/scale}
    \end{center}
    \rotatebox{90}{Velocity norm in ${\Omega_{\text{ane}}}$}\quad
    \includegraphics[width=0.29\textwidth]{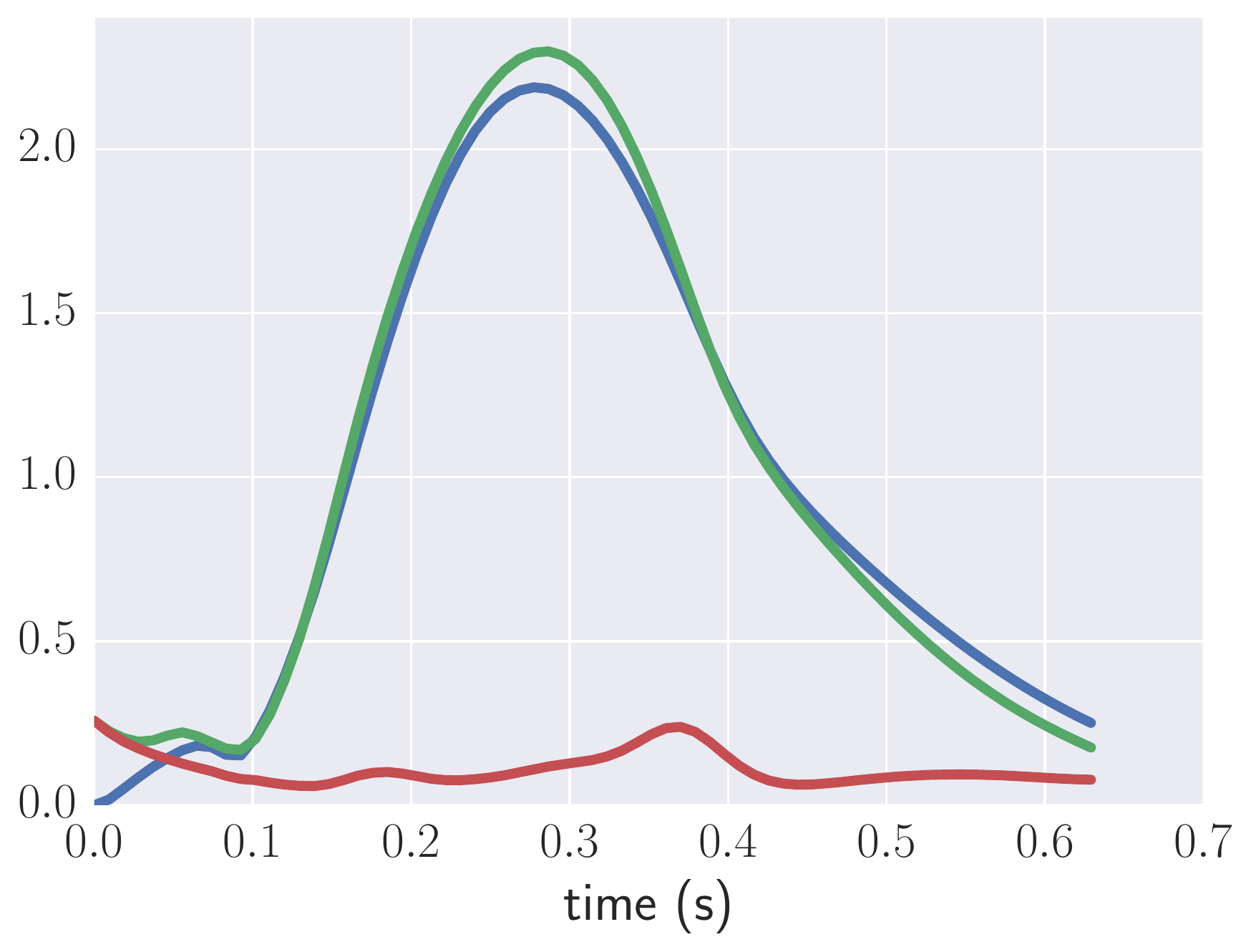}
    \includegraphics[width=0.29\textwidth]{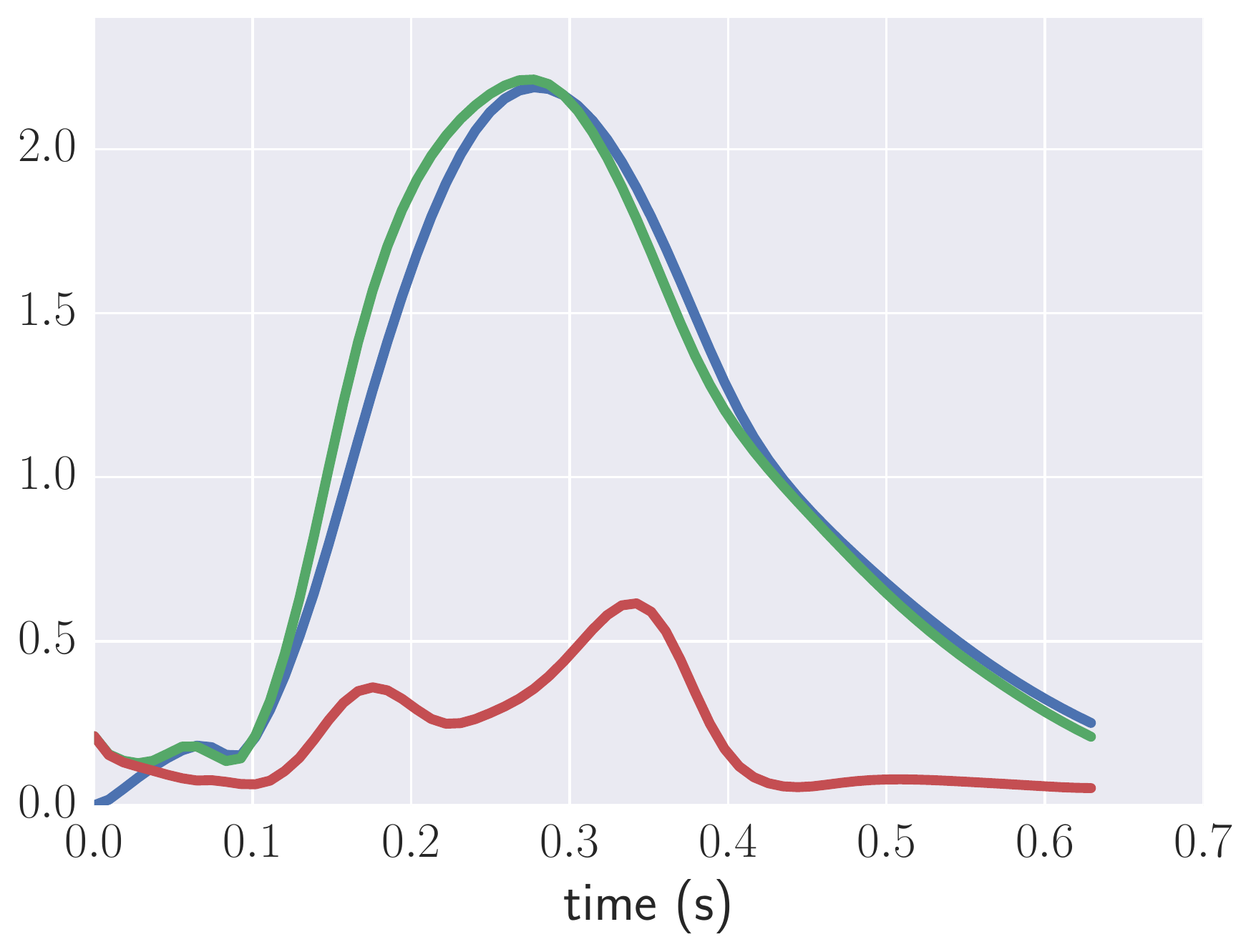}
    \includegraphics[width=0.29\textwidth]{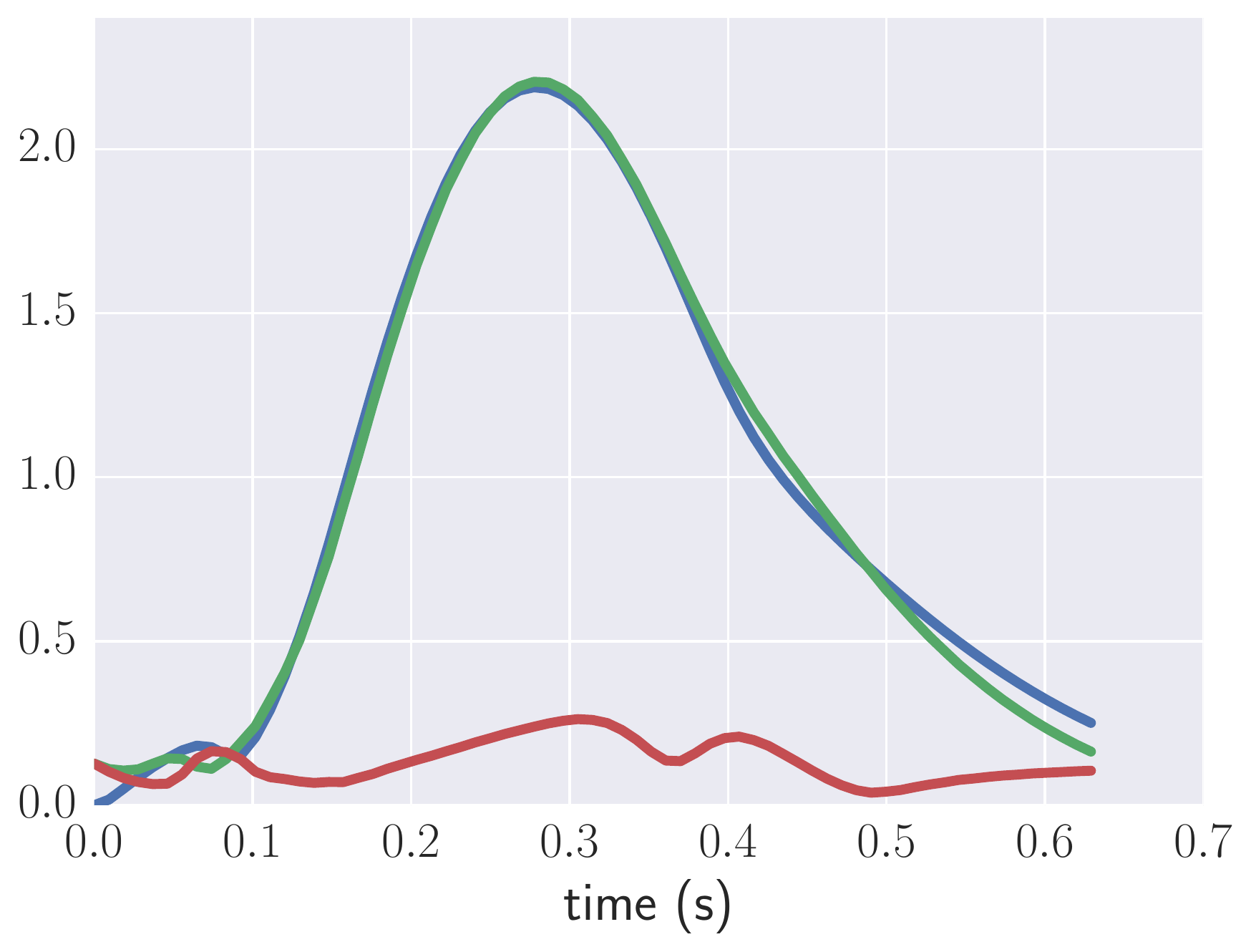} \\
    \rotatebox{90}{WSS norm on $\Gamma_{\text{ane}}$}\quad
    \includegraphics[width=0.29\textwidth]{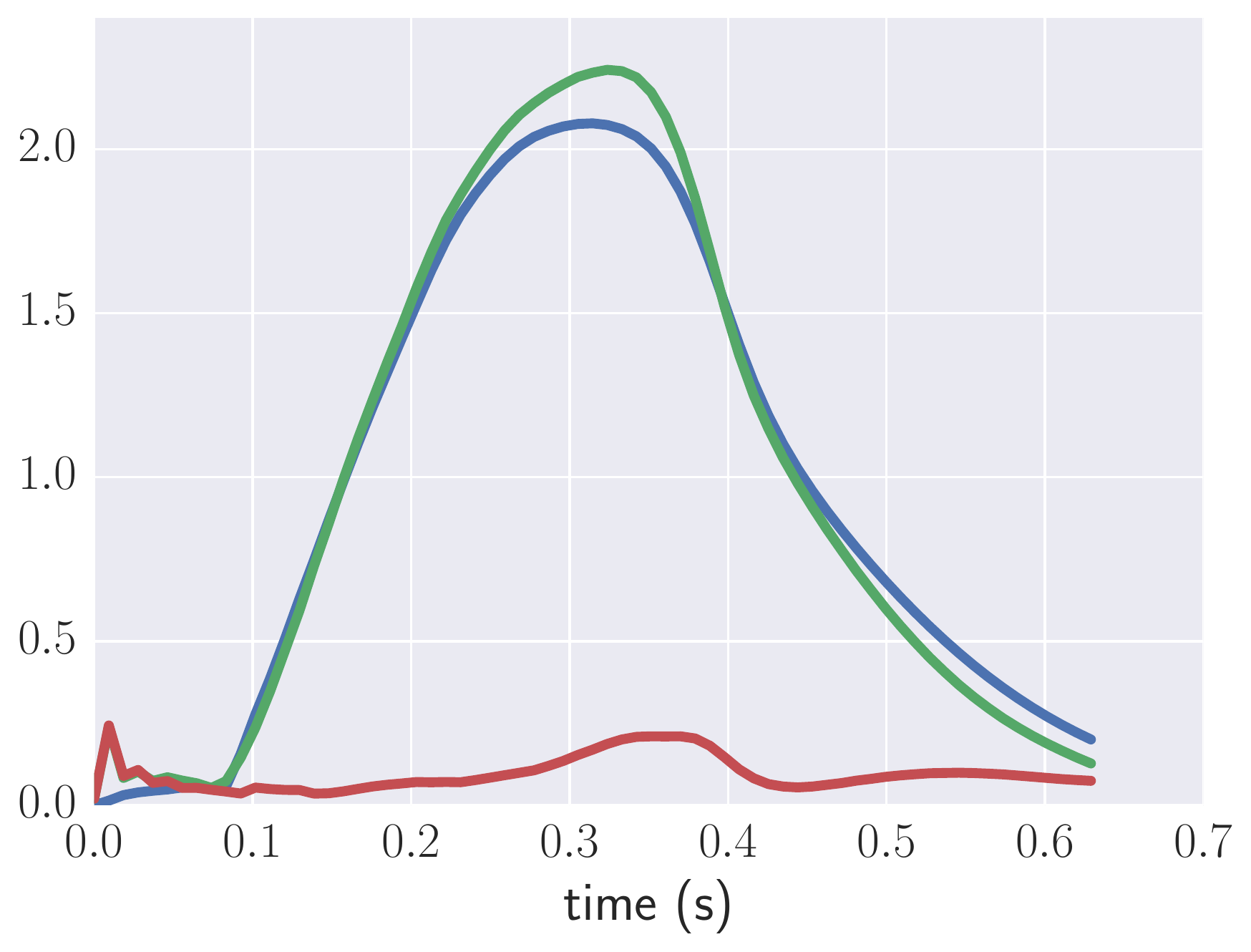}
    \includegraphics[width=0.29\textwidth]{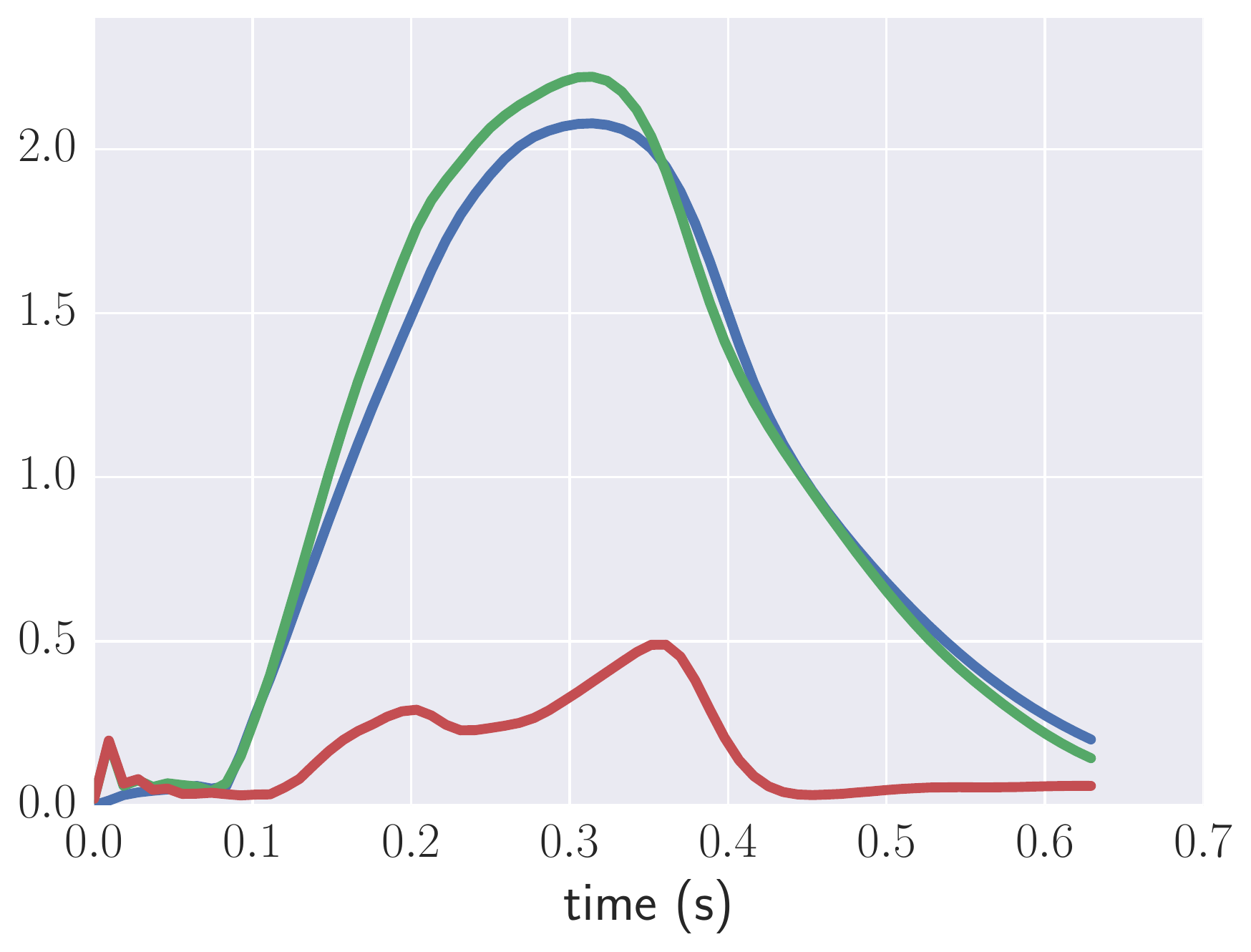}
    \includegraphics[width=0.29\textwidth]{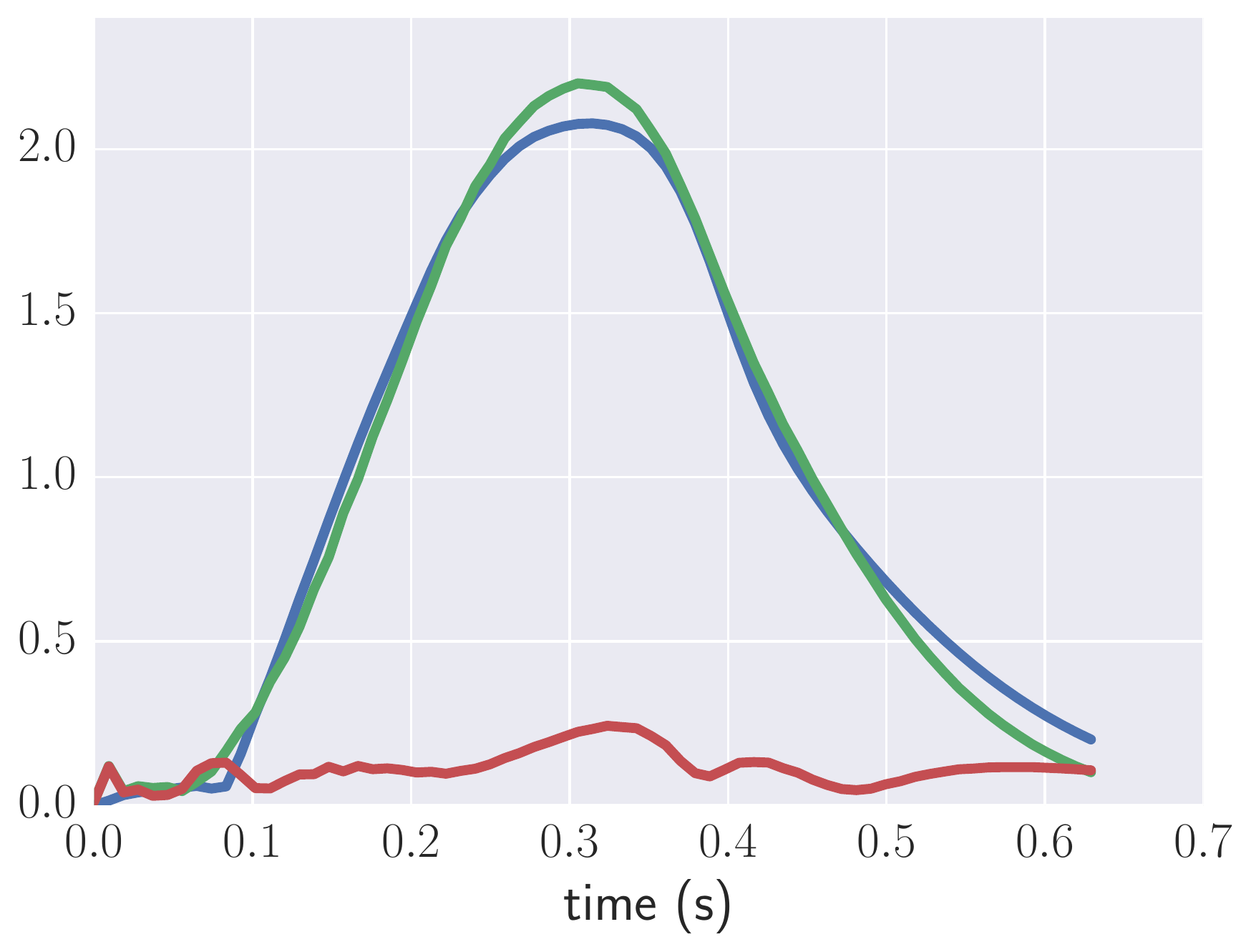} \\
    \includegraphics[width=0.59\textwidth]{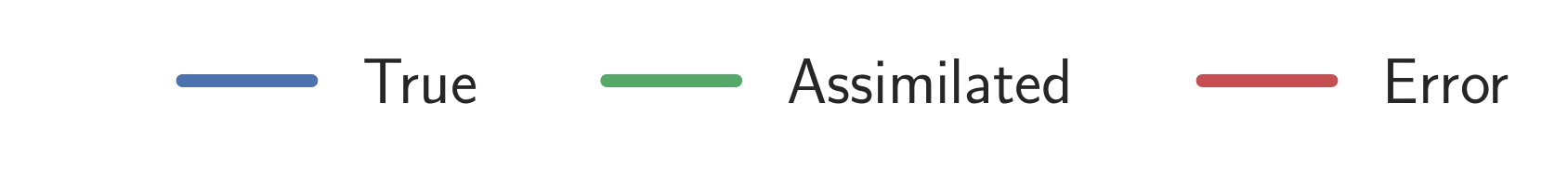}
    \caption{Results using the \textbf{time-averaging observation operator} with pointwise additive Gaussian white noise.
The signal-to-noise ratio was computed as
 $\norm{\mathcal{T}^{\text{avg}}u_{\text{true}}}^2 / \norm{\mathcal{T}^{\text{avg}}u_{\text{true}}-d}^2$, where $d$ is the noisy data.
The snapshots on the top three rows are taken at $t=0.296$s.}
    \label{fig:results_aneurysm_avg_noise}
\end{figure}

\textbf{Noisy observations (figures  \ref{fig:results_aneurysm_inst_noise} and \ref{fig:results_aneurysm_avg_noise}).}
Pointwise Gaussian white noise was added to the observations with zero mean and
varying magnitude. This type of noise is not expected in real observations, in
particular because it depends on the numerical mesh. Nevertheless, we consider
it as a suitable benchmark setup. The results of the reconstruction for different signal to noise ratios
are shown in
figure \ref{fig:results_aneurysm_inst_noise}, and
for the instantaneous observation operator $\mathcal{T}^{\text{inst}}$
and figure \ref{fig:results_aneurysm_avg_noise} for the time-averaging observation operator
$\mathcal{T}^{\text{avg}}$.
With increasing level of noise, the optimised functional value $\mathcal J + \mathcal R$ increases,
because of the increased the difference between reconstructed and observed velocity.
Nevertheless, the error measures remain small, showing that the reconstruction works reliable even
for high noise to signal ratios.
Overall, the assimilated flows and metrics agree well for all noise
levels, and one can conclude that the reconstruction is little affected by this type of
noise.

\begin{figure}
    \centering
        \footnotesize
    \hspace{1cm}
    \parbox[b][8.5em][t]{0.29\textwidth}{
        \emph{Base setup with\\
        $\alpha = \gamma = 10^{-4}$}
        \vspace{0.5em}\\
        \input{results_aneurysm/nsassimilation_alpha_beta_gamma1e4/source/results_aneurysm/instant/assimilated_H1H1_0_noise/metrics}
    }
    \parbox[b][8.5em][t]{0.29\textwidth}{
        \emph{Base setup with \\$\alpha = \gamma = 10^{-2}$}
        \vspace{0.5em}\\
        \input{results_aneurysm/nsassimilation_alpha_beta_gamma1e2/source/results_aneurysm/instant/assimilated_H1H1_0_noise/metrics}
    }
    \parbox[b][8.5em][t]{0.29\textwidth}{
        \emph{Base setup with \\$\alpha = \gamma = 1$}
        \vspace{0.5em}\\
        \input{results_aneurysm/nsassimilation_alpha_beta_gamma1e0/source/results_aneurysm/instant/assimilated_H1H1_0_noise/metrics}
    }
    \\
    \rotatebox{90}{\quad \quad Observation}\quad
    \includegraphics[width=0.29\textwidth]{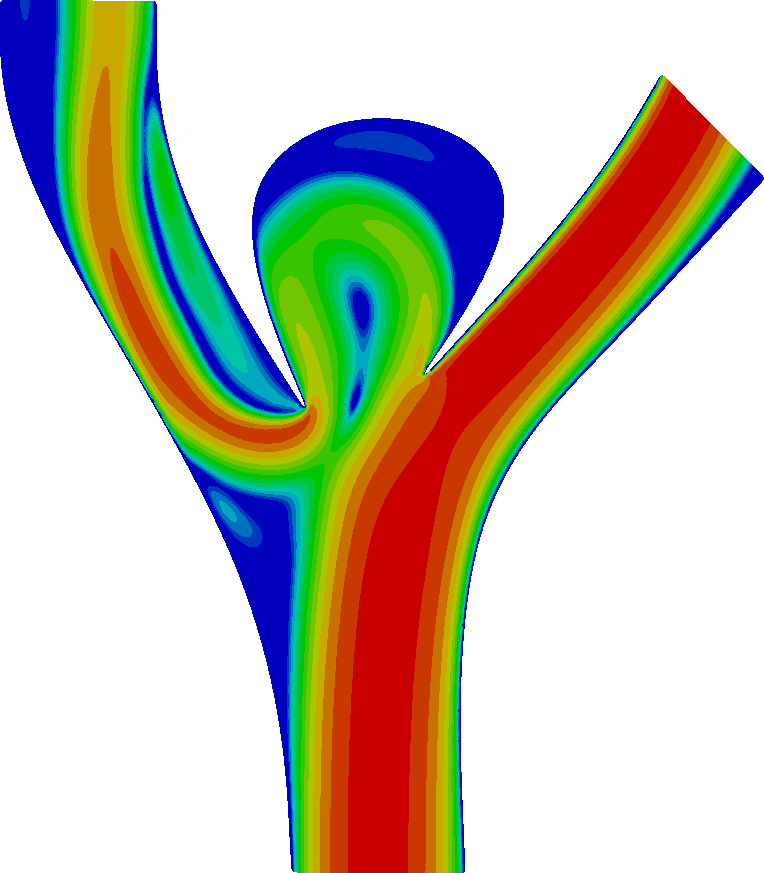}
    \includegraphics[width=0.29\textwidth]{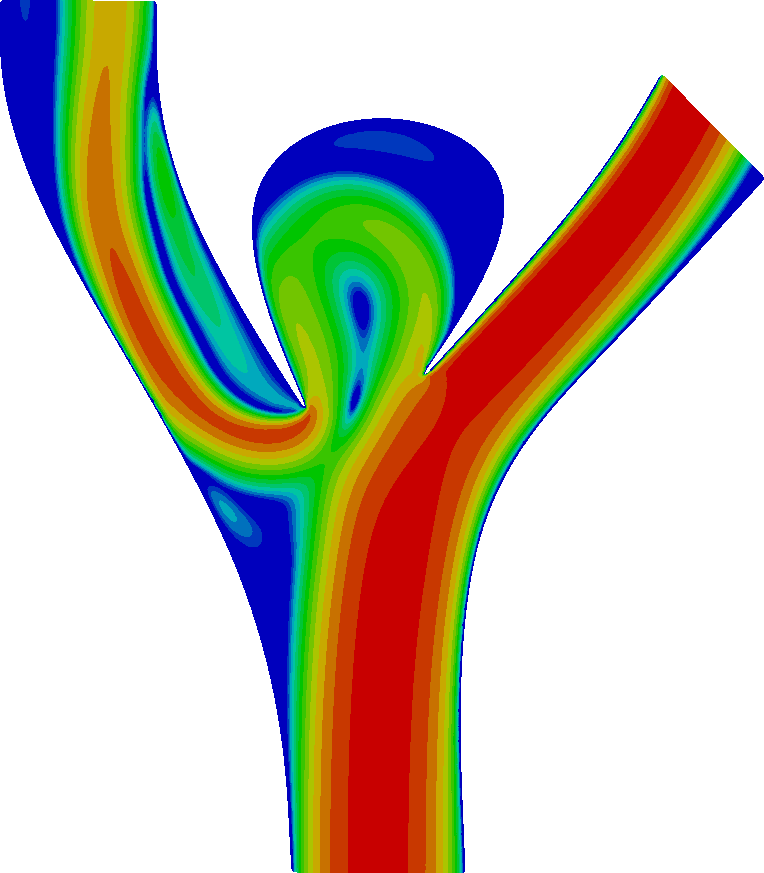}
    \includegraphics[width=0.29\textwidth]{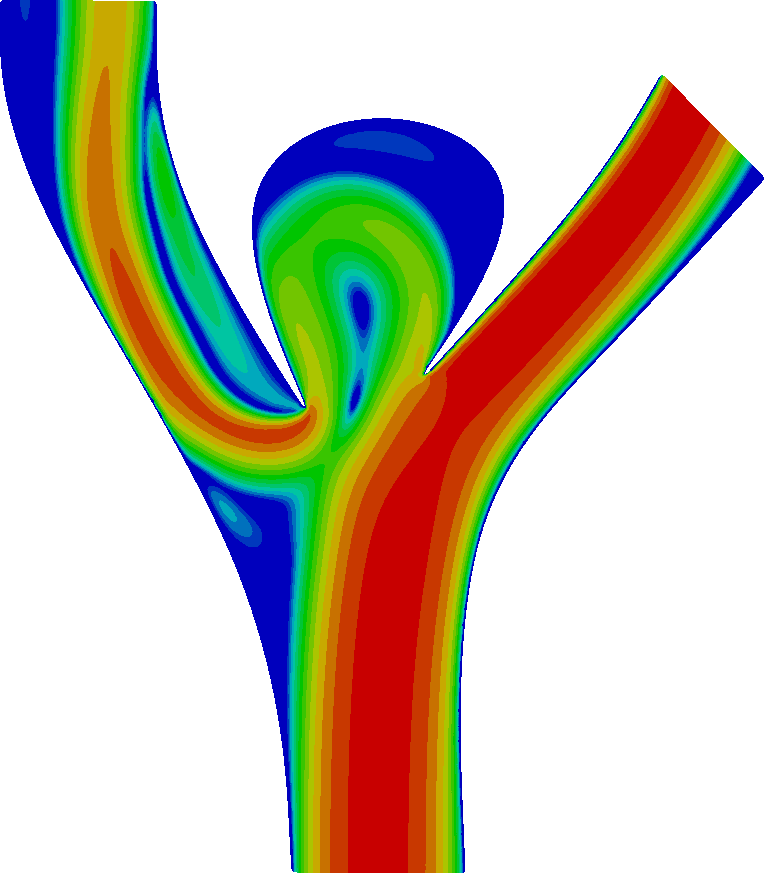}\\
    \rotatebox{90}{\qquad \quad Assimilation}\quad
    \includegraphics[width=0.29\textwidth]{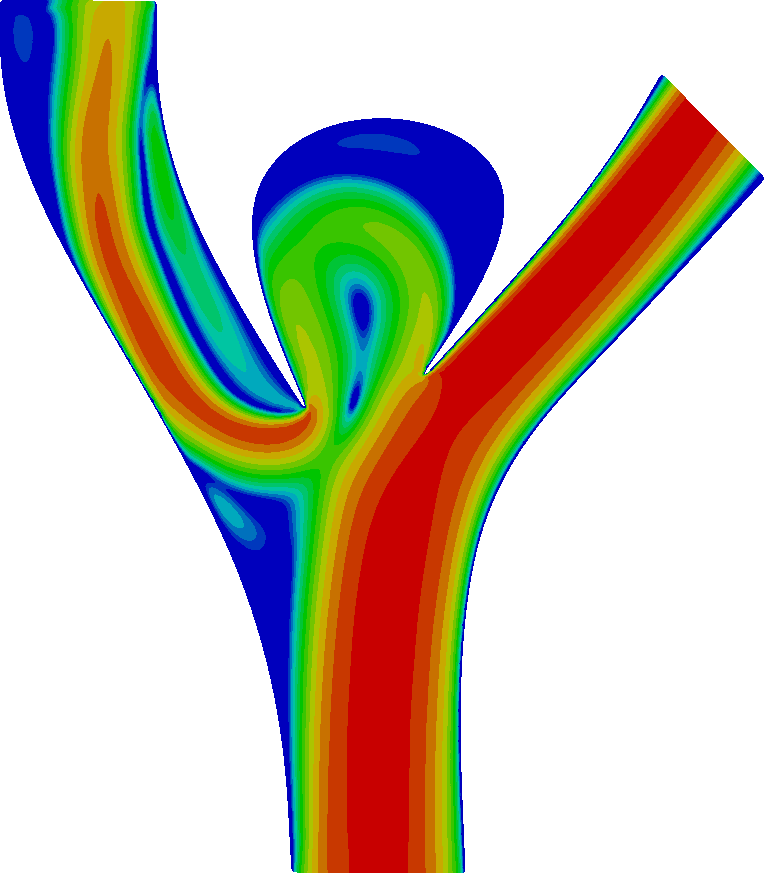}
    \includegraphics[width=0.29\textwidth]{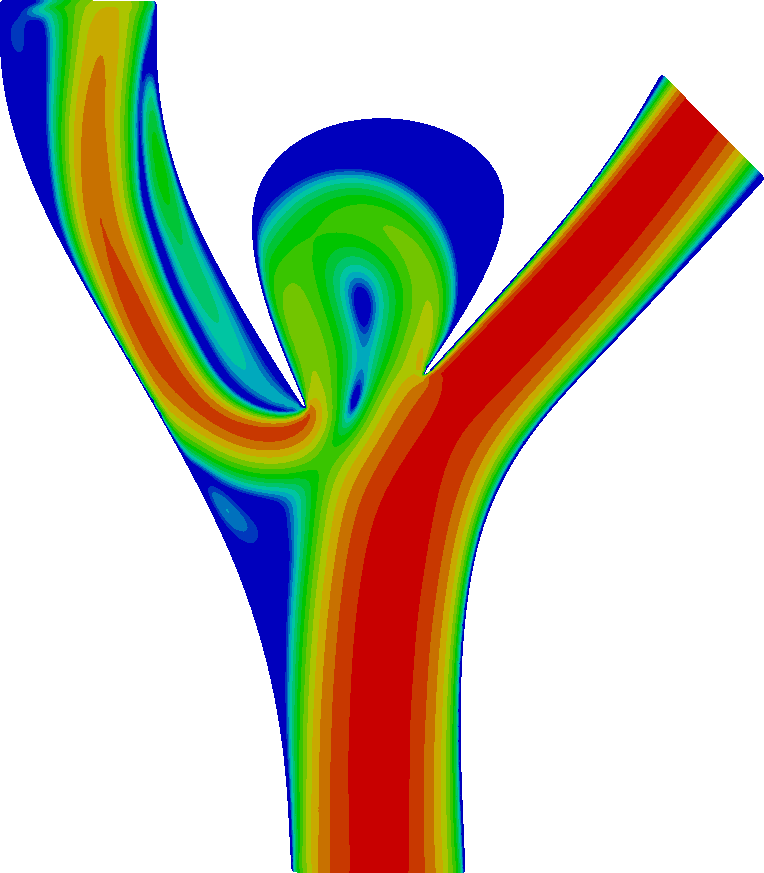}
    \includegraphics[width=0.29\textwidth]{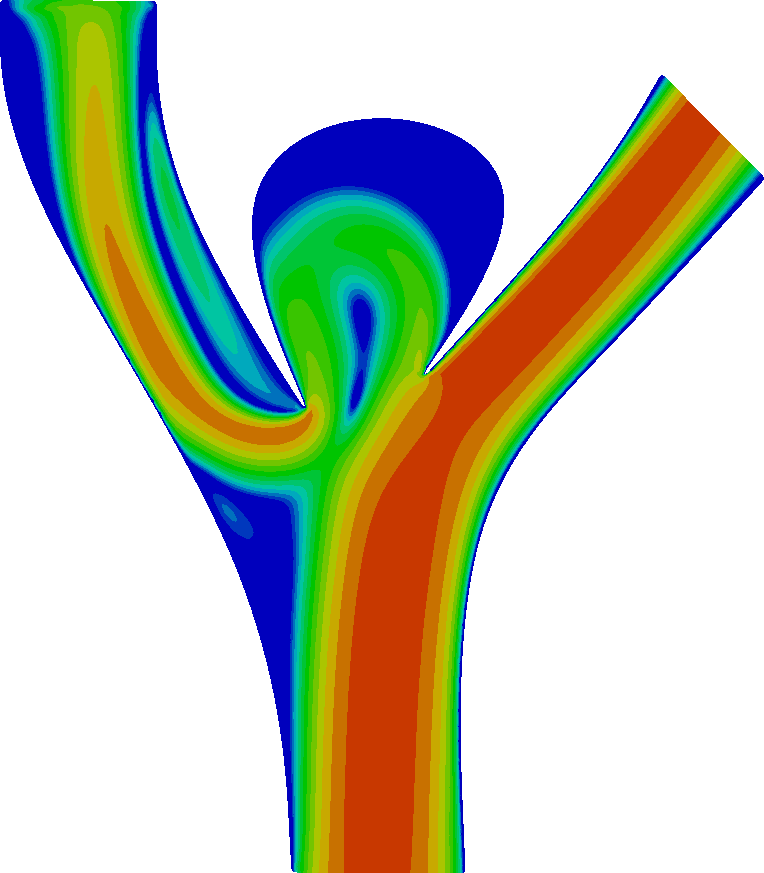}\\
    \begin{center}
        \includegraphics[width=0.3\textwidth]{results_aneurysm/scale}
    \end{center}
    \rotatebox{90}{Velocity norm in ${\Omega_{\text{ane}}}$}\quad
    \includegraphics[width=0.29\textwidth]{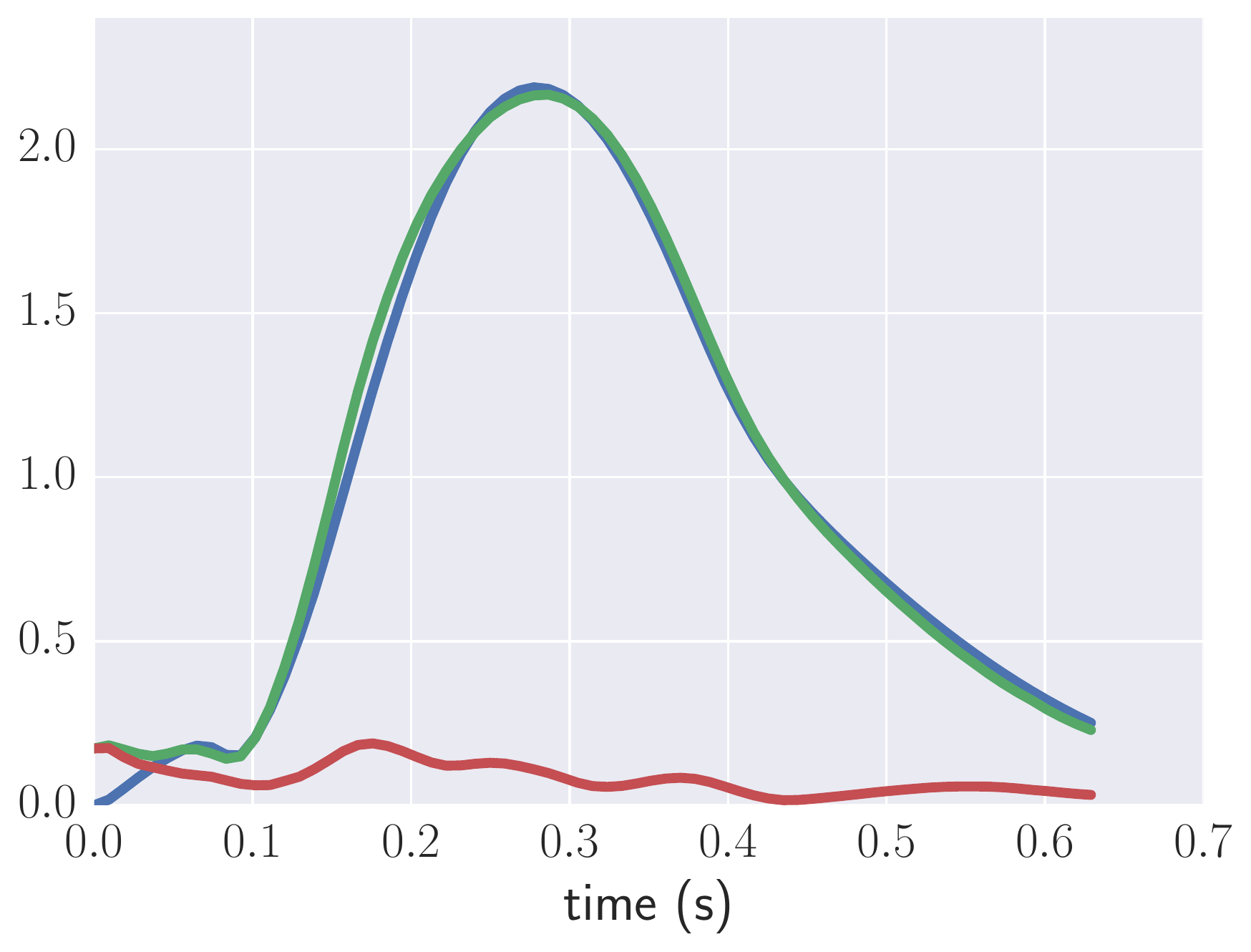}
    \includegraphics[width=0.29\textwidth]{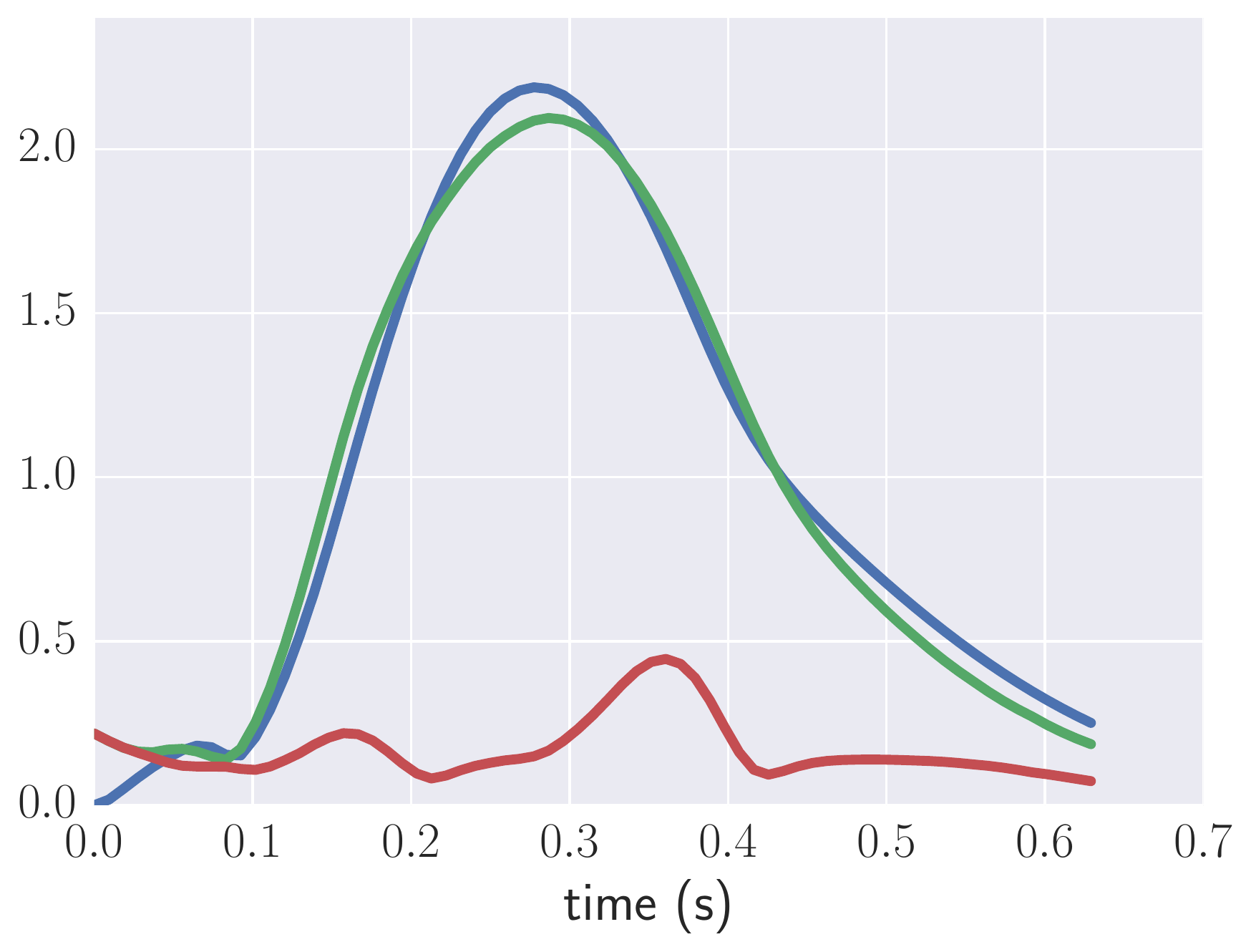}
    \includegraphics[width=0.29\textwidth]{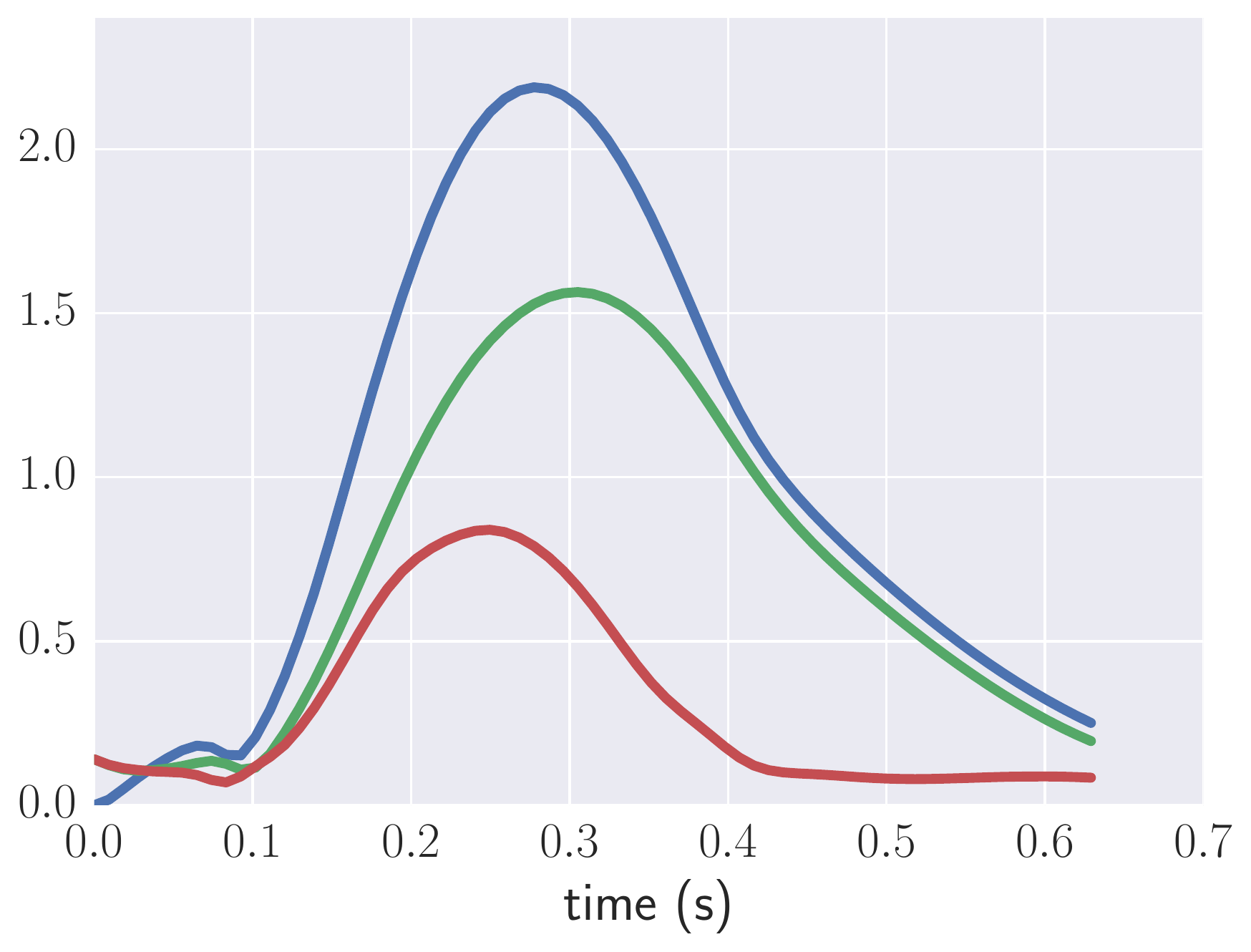}\\
    \rotatebox{90}{WSS norm on $\Gamma_{\text{ane}}$}\quad
    \includegraphics[width=0.29\textwidth]{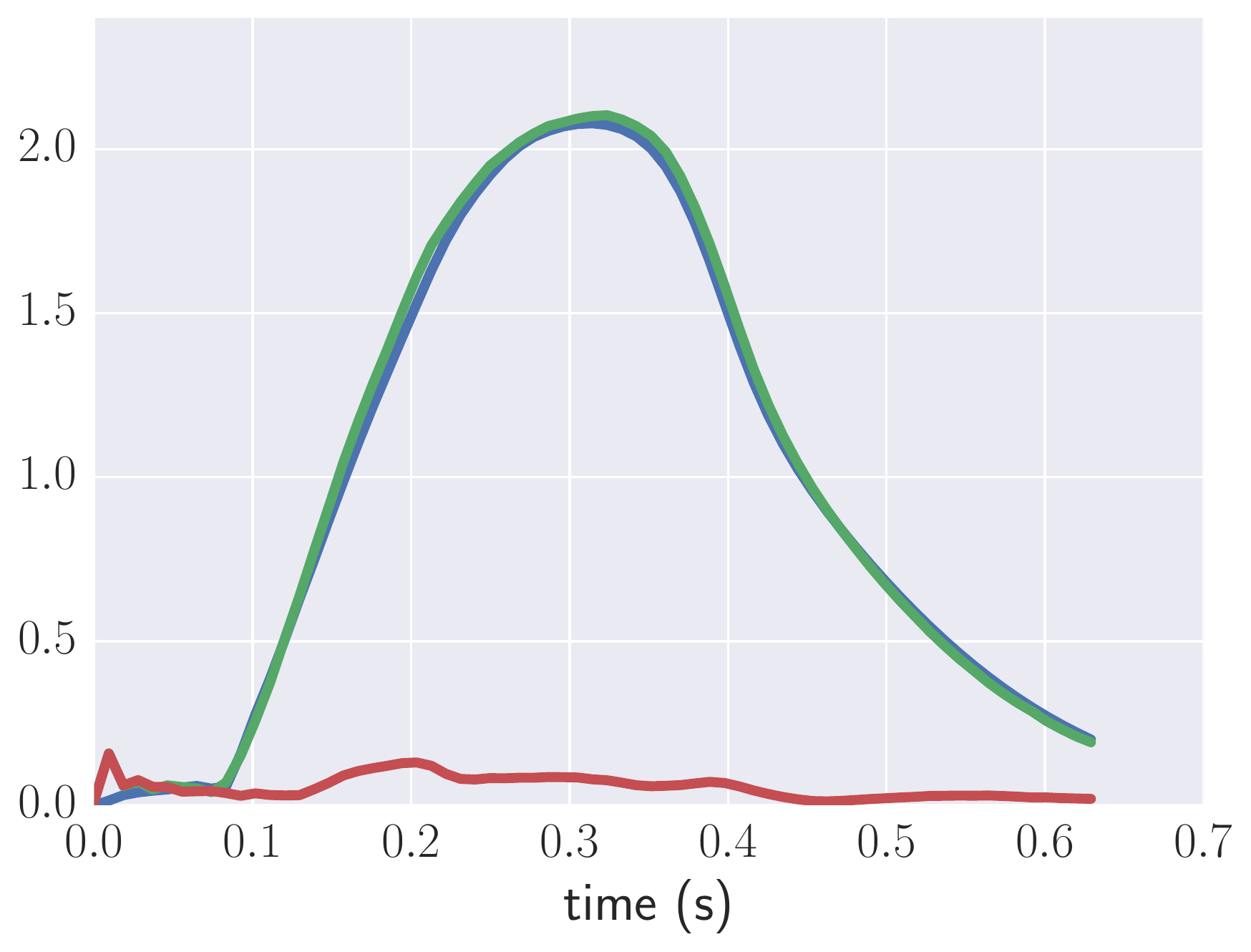}
    \includegraphics[width=0.29\textwidth]{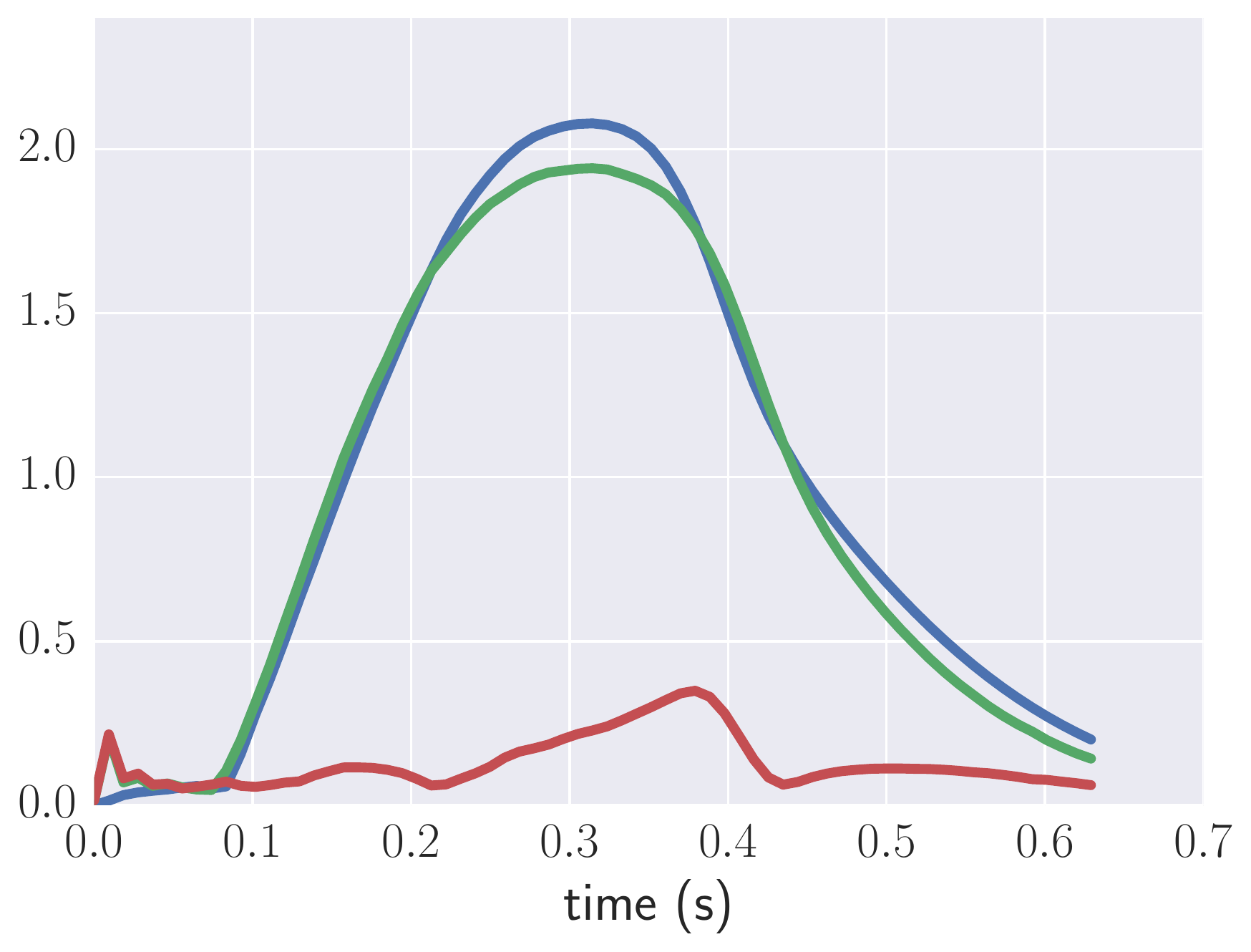}
    \includegraphics[width=0.29\textwidth]{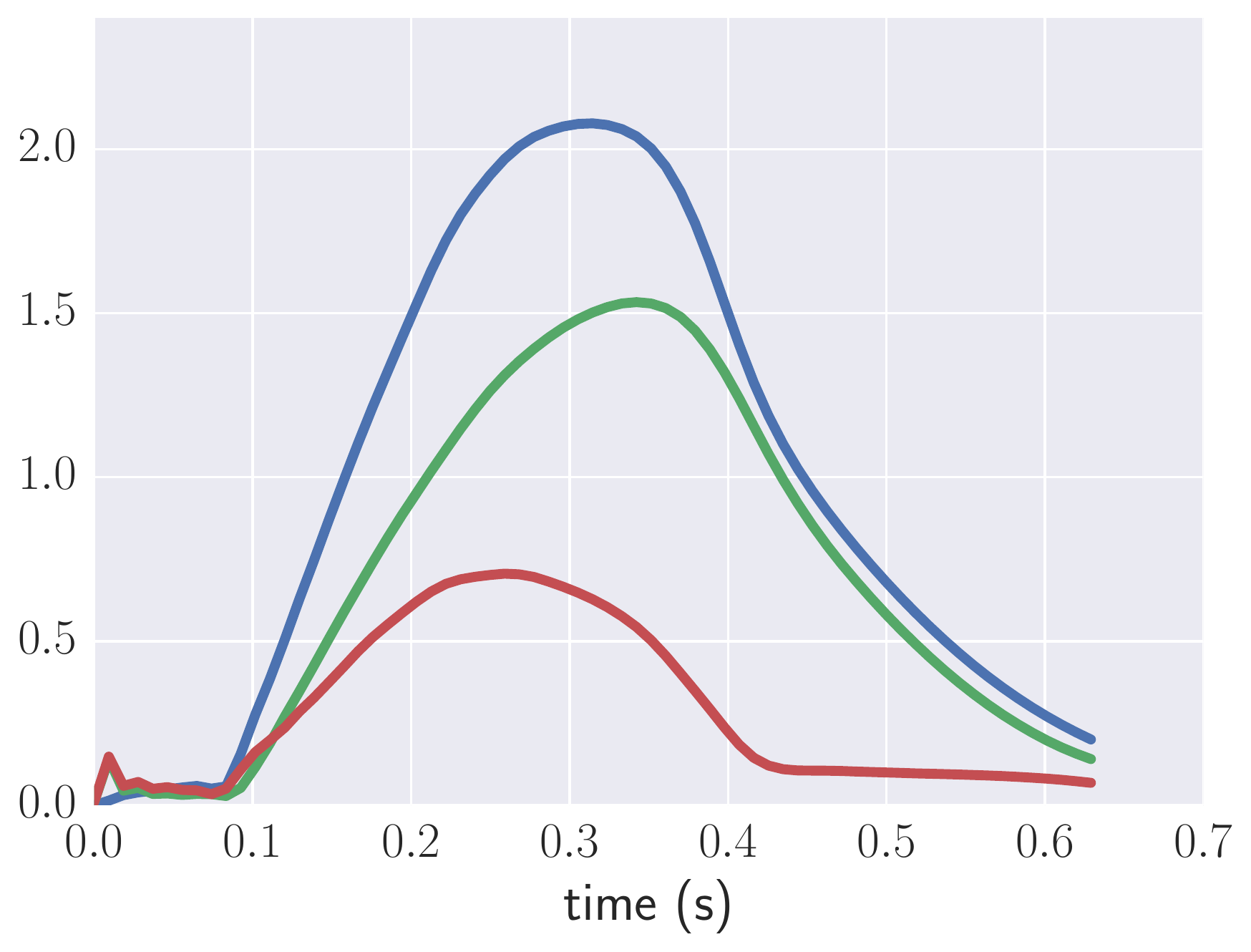}\\
    \includegraphics[width=0.59\textwidth]{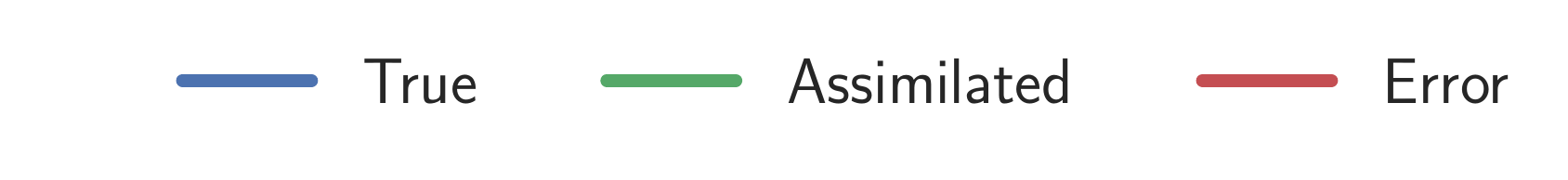}
    \caption{Results using the \textbf{instantaneous observation operator} with varying $\alpha$ and $\gamma$ regularisation coefficients.
        The base setup (figure~\ref{fig:results_aneurysm_inst_noise}, left column) uses $\alpha=\gamma=10^{-5}$.
        The snapshots on the top two rows are taken at $t=0.296$s.}
    \label{fig:results_aneurysm_instant_alpha_beta_gamma}
\end{figure}

\begin{figure}
    \centering
        \footnotesize
    \hspace{1cm}
    \parbox[b][8.5em][t]{0.29\textwidth}{
        \emph{Base setup with\\
        $\alpha = \gamma = 10^{-4}$}
        \vspace{0.5em}\\
        \input{results_aneurysm/nsassimilation_alpha_beta_gamma1e4/source/results_aneurysm/averaged/assimilated_H1H1_0_noise/metrics}
    }
    \parbox[b][8.5em][t]{0.29\textwidth}{
        \emph{Base setup with \\$\alpha = \gamma = 10^{-2}$}
        \vspace{0.5em}\\
        \input{results_aneurysm/nsassimilation_alpha_beta_gamma1e2/source/results_aneurysm/averaged/assimilated_H1H1_0_noise/metrics}
    }
    \parbox[b][8.5em][t]{0.29\textwidth}{
        \emph{Base setup with \\$\alpha = \gamma = 1$}
        \vspace{0.5em}\\
        \input{results_aneurysm/nsassimilation_alpha_beta_gamma1e0/source/results_aneurysm/averaged/assimilated_H1H1_0_noise/metrics}
    }
    \\
    \rotatebox{90}{\qquad \quad Observation}\quad
    \includegraphics[width=0.29\textwidth]{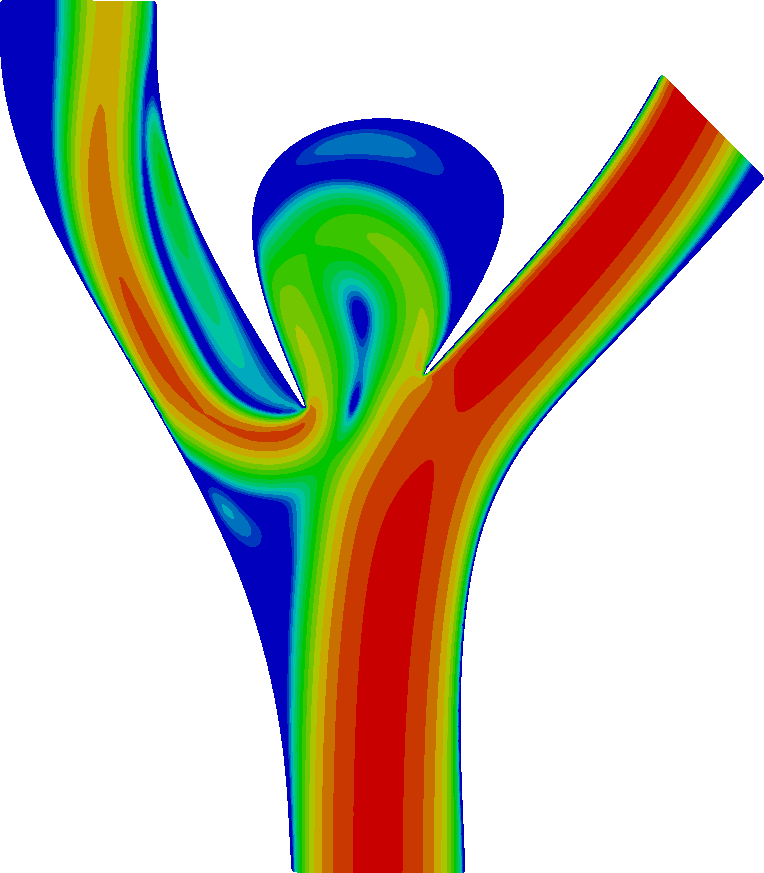}
    \includegraphics[width=0.29\textwidth]{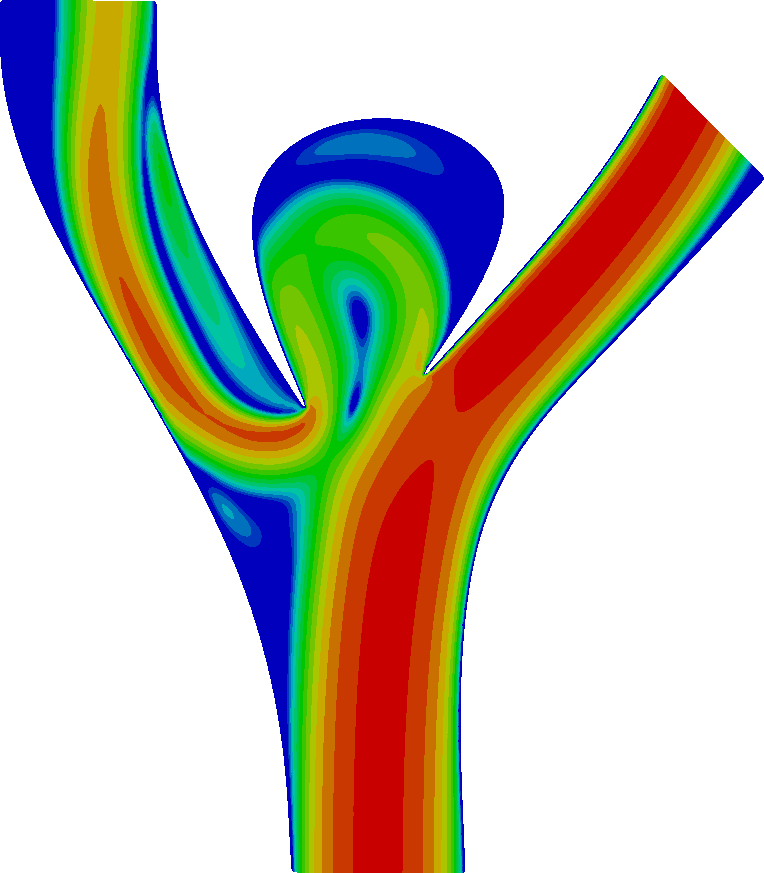}
    \includegraphics[width=0.29\textwidth]{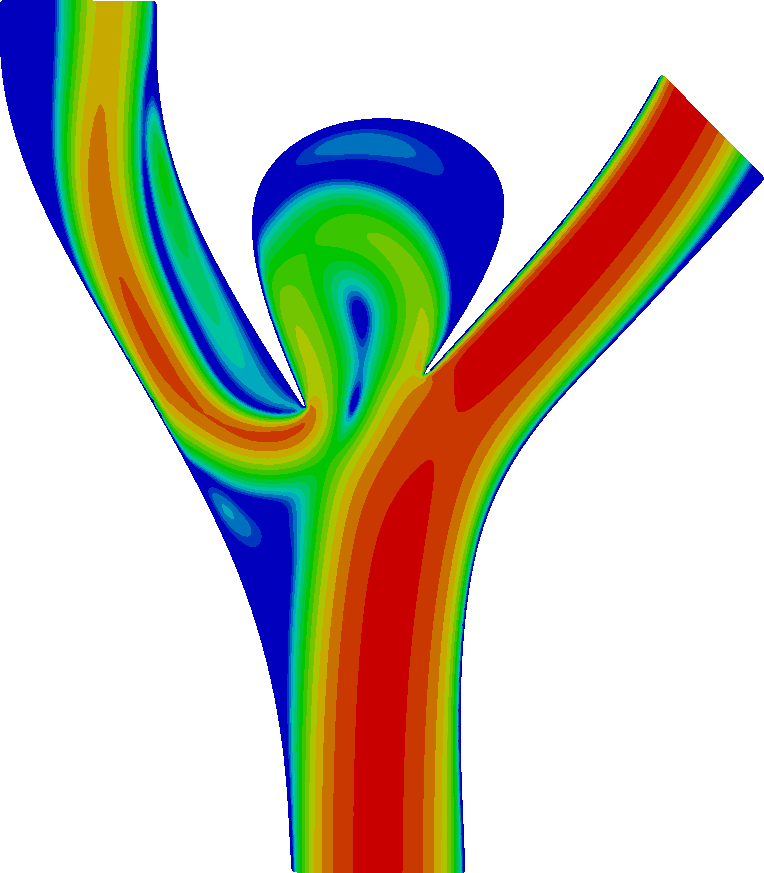}\\
    \rotatebox{90}{\qquad \quad Assimilation}\quad
    \includegraphics[width=0.29\textwidth]{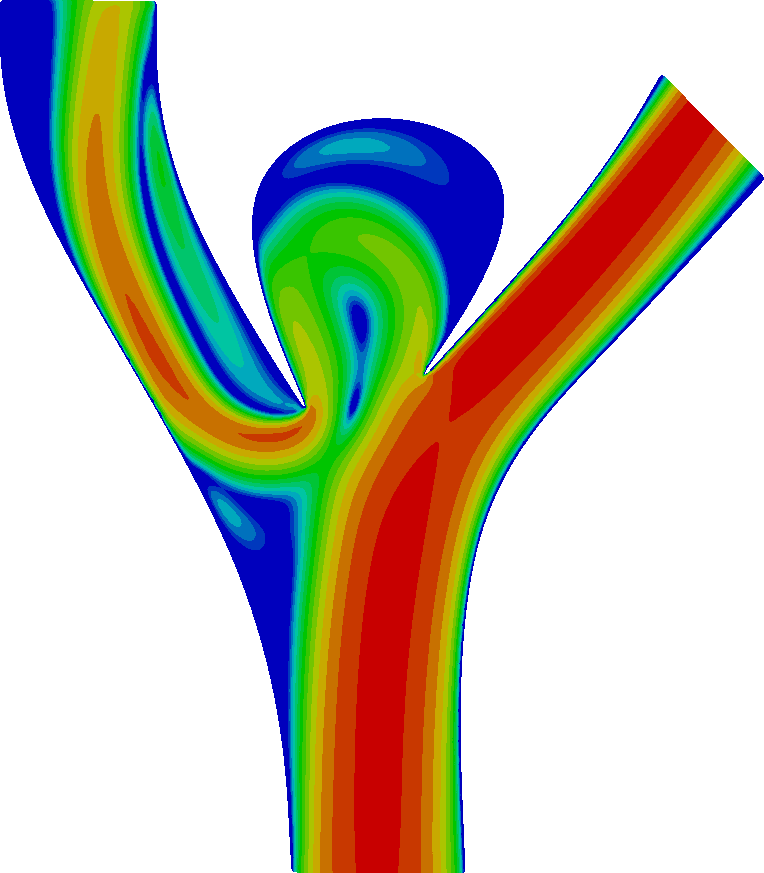}
    \includegraphics[width=0.29\textwidth]{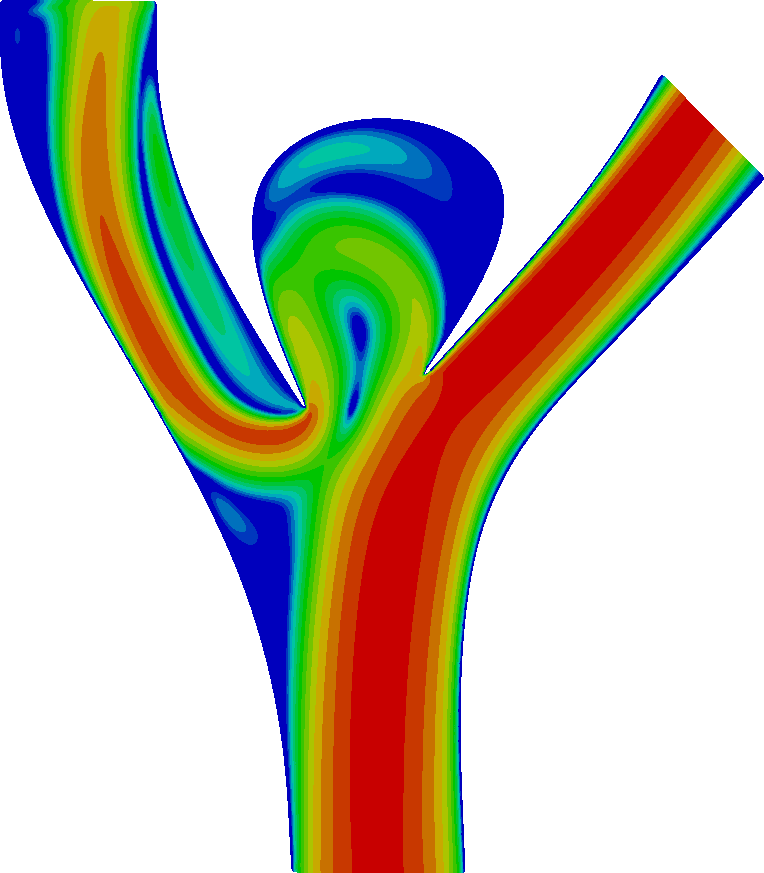}
    \includegraphics[width=0.29\textwidth]{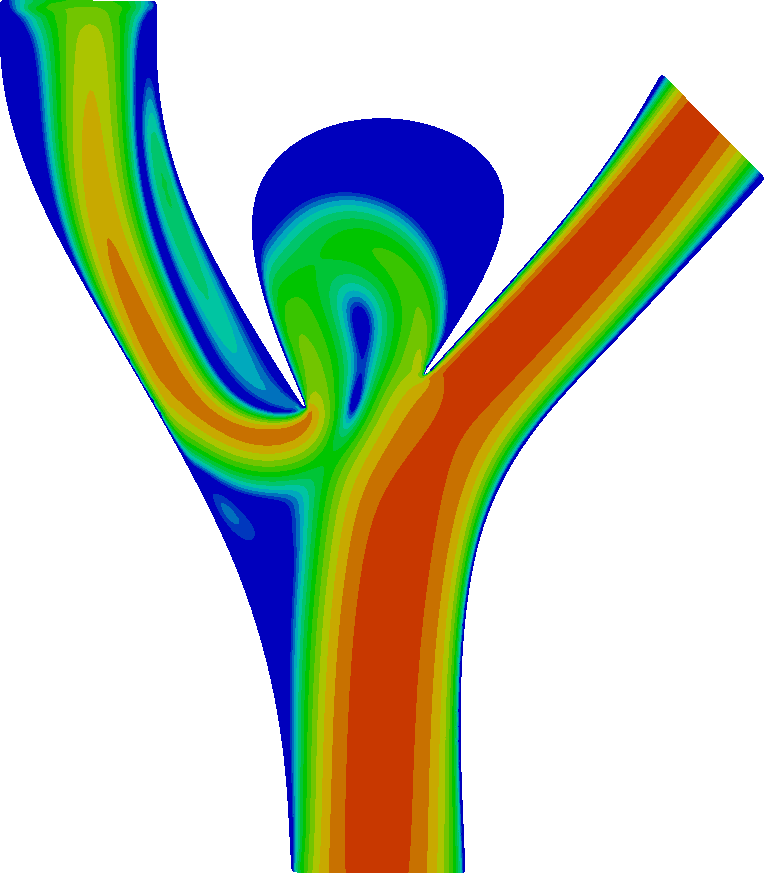}\\
    \begin{center}
        \includegraphics[width=0.3\textwidth]{results_aneurysm/scale}
    \end{center}
    \rotatebox{90}{Velocity norm in ${\Omega_{\text{ane}}}$}\quad
    \includegraphics[width=0.29\textwidth]{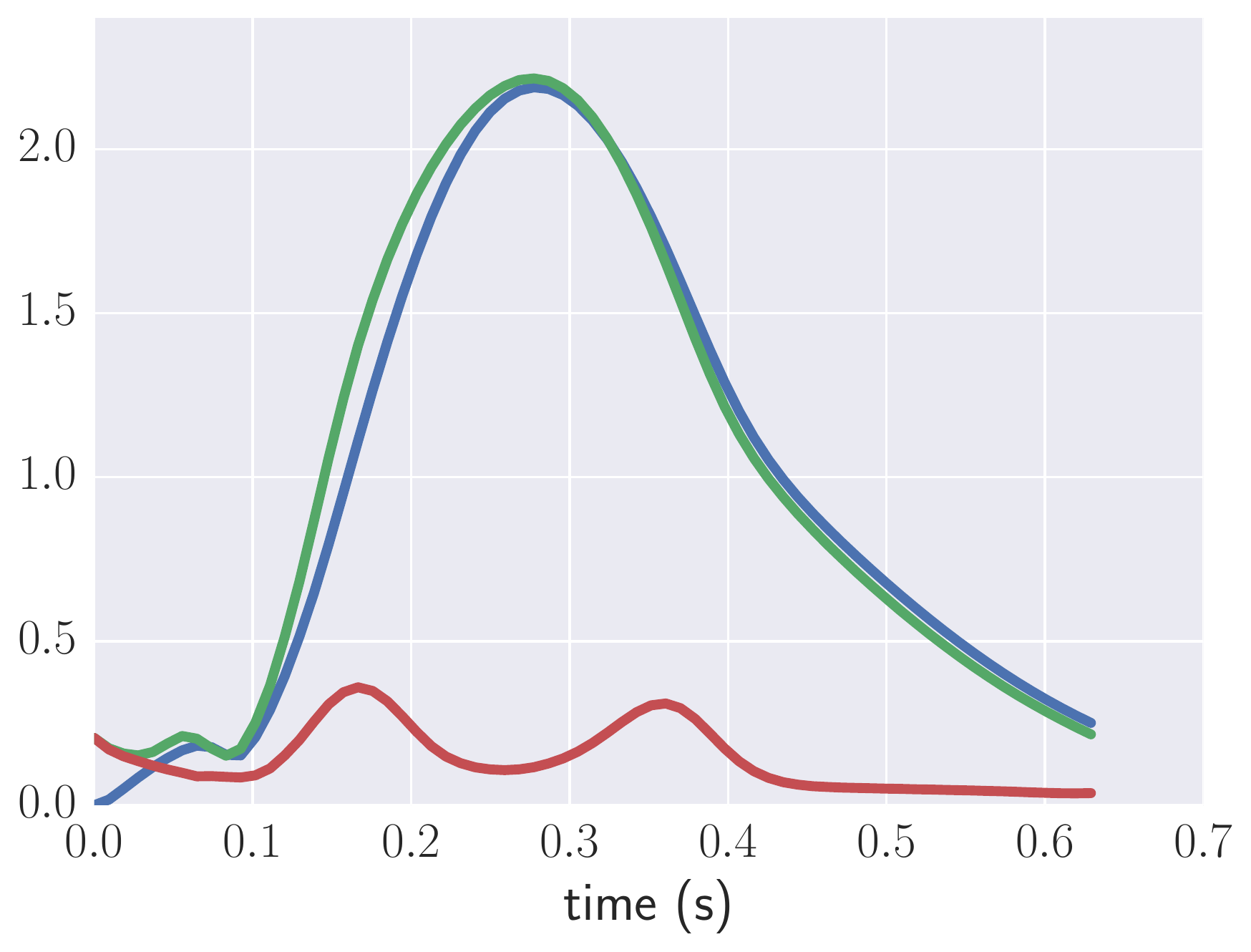}
    \includegraphics[width=0.29\textwidth]{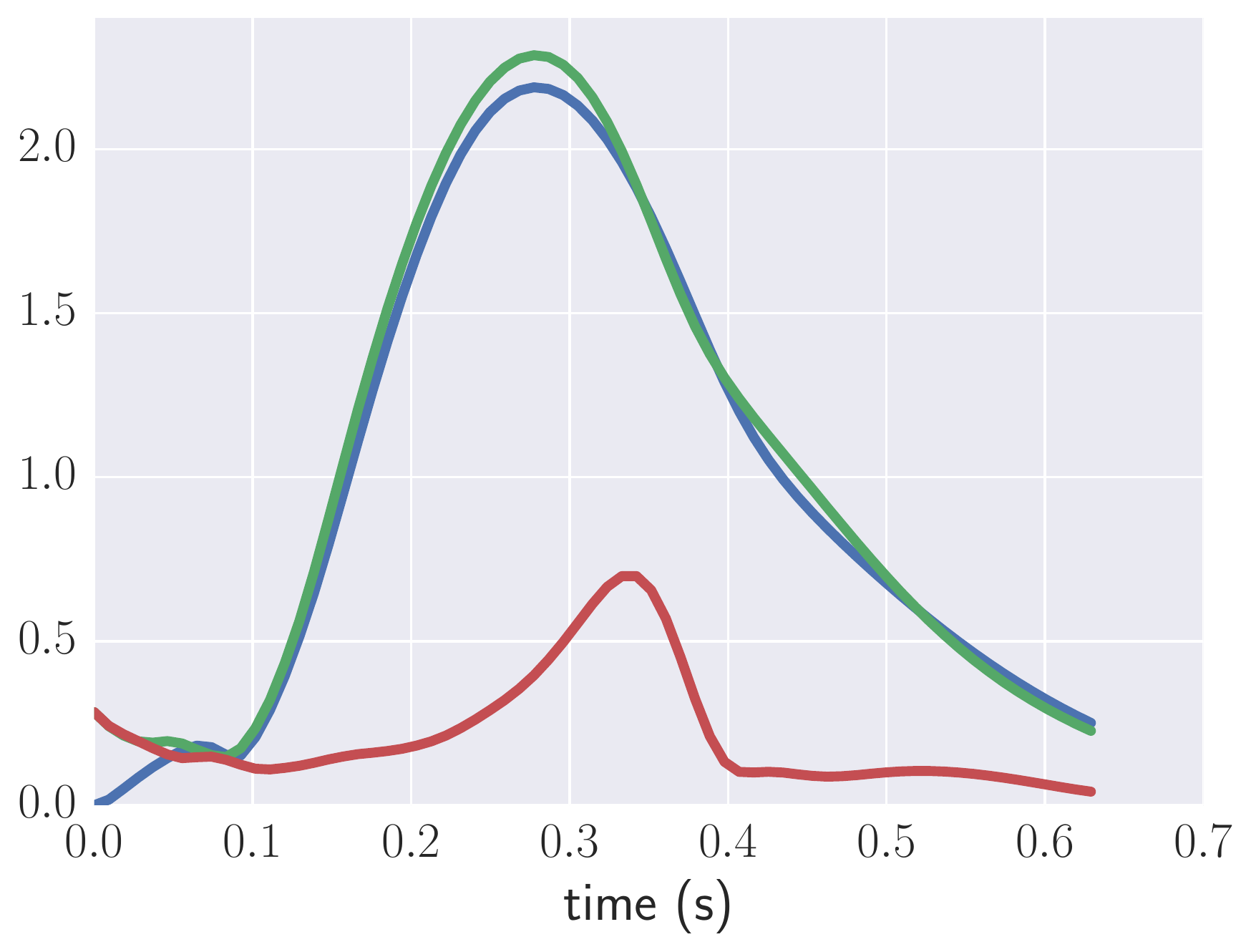}
    \includegraphics[width=0.29\textwidth]{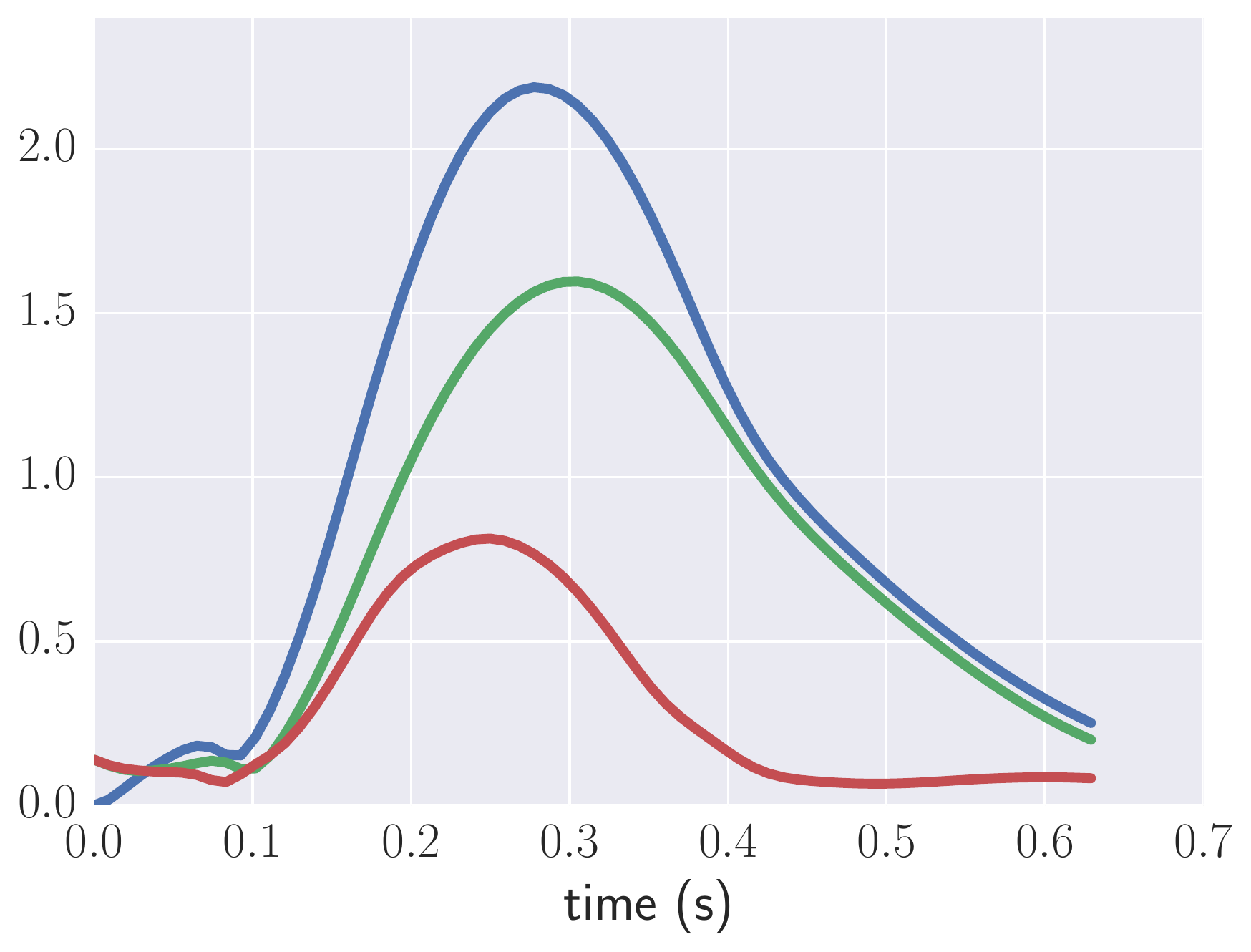}\\
    \rotatebox{90}{WSS norm on $\Gamma_{\text{ane}}$}\quad
    \includegraphics[width=0.29\textwidth]{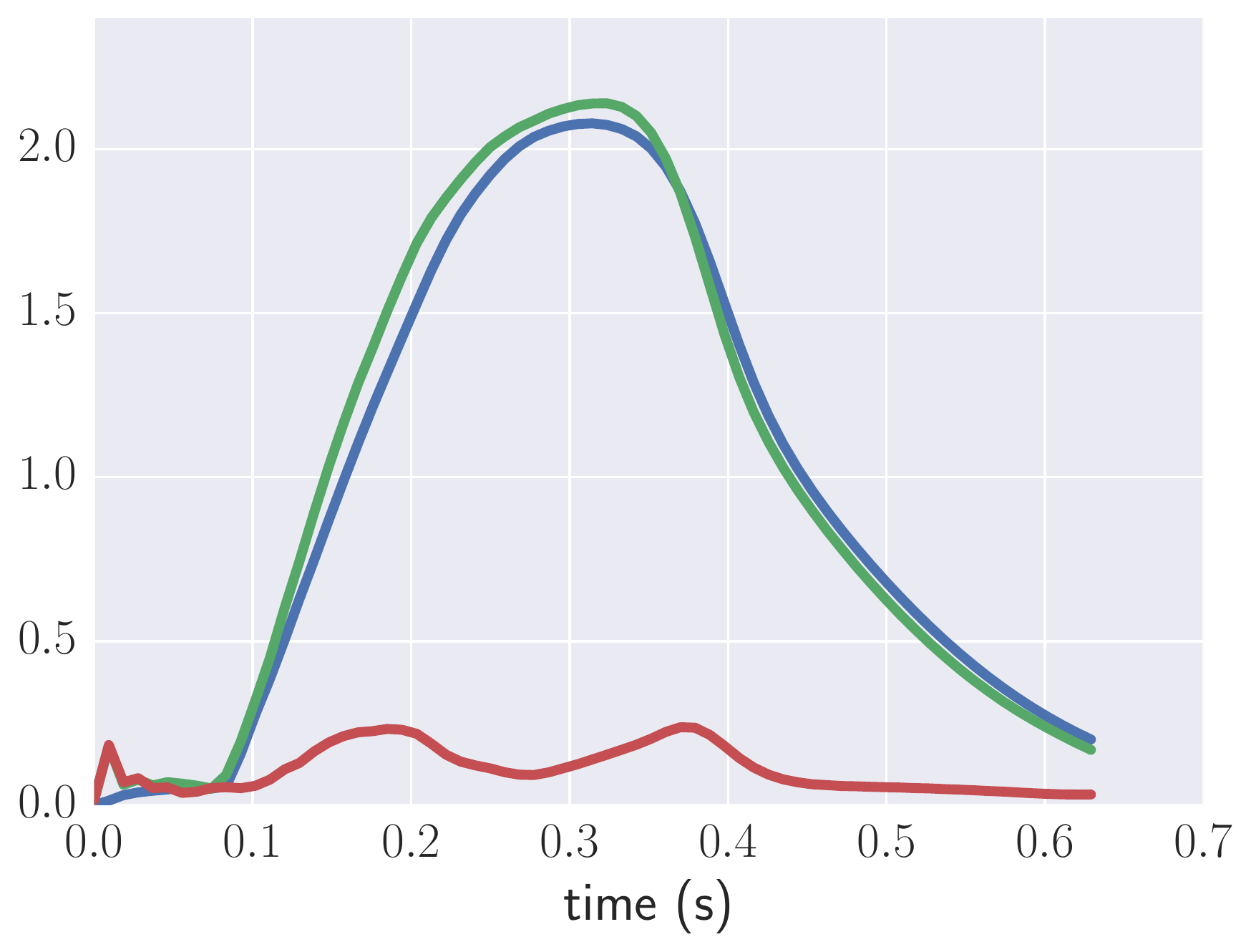}
    \includegraphics[width=0.29\textwidth]{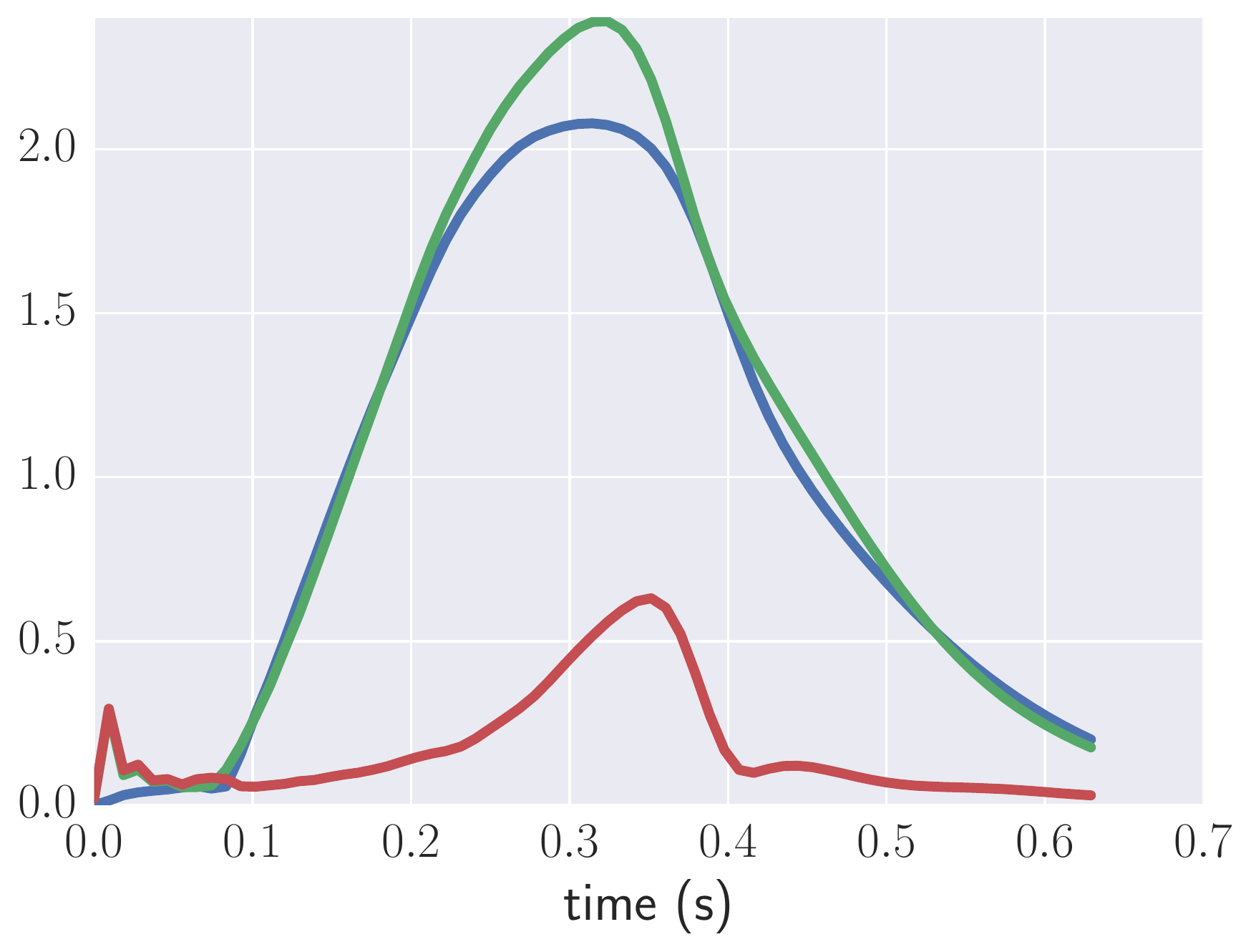}
    \includegraphics[width=0.29\textwidth]{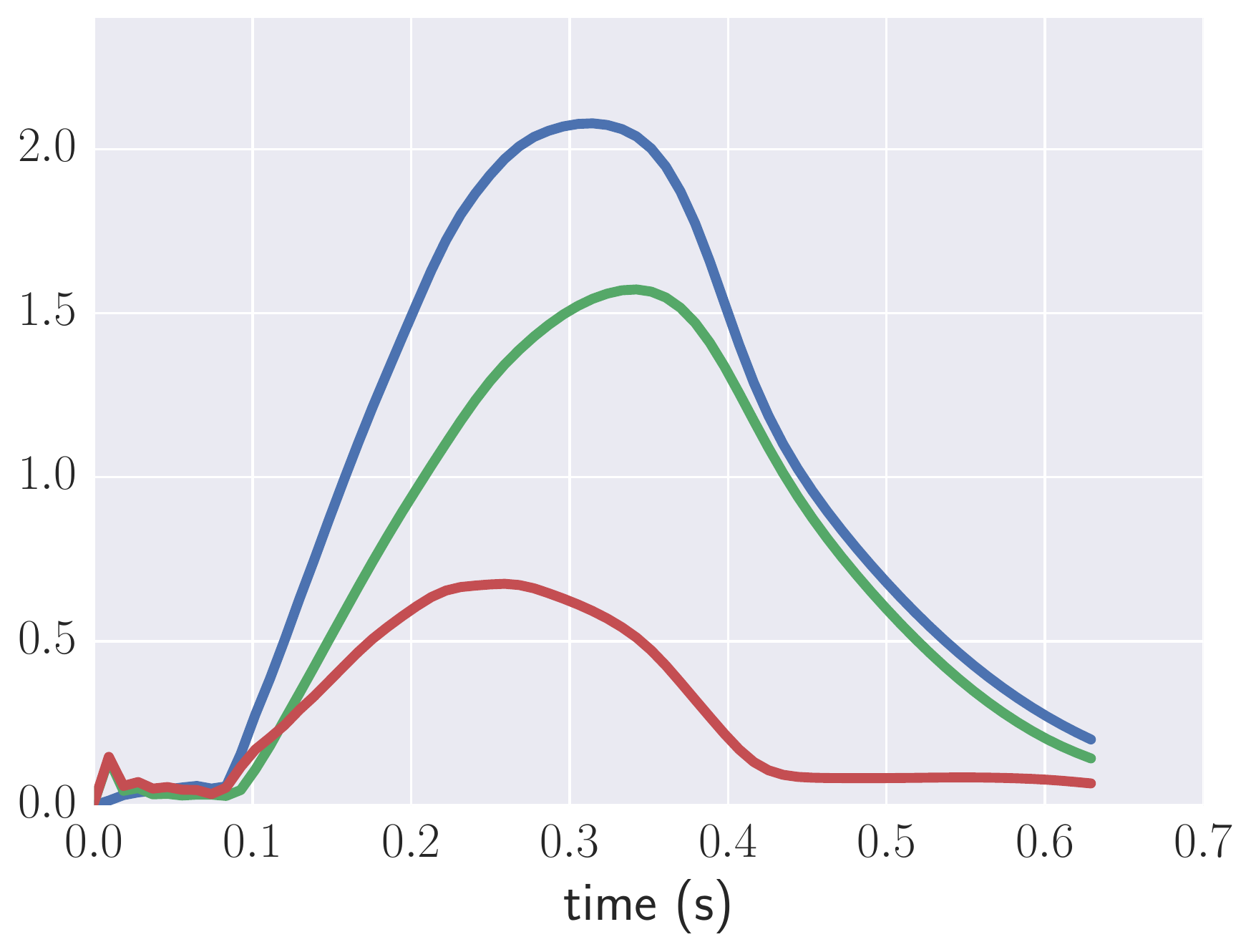}\\
    \includegraphics[width=0.59\textwidth]{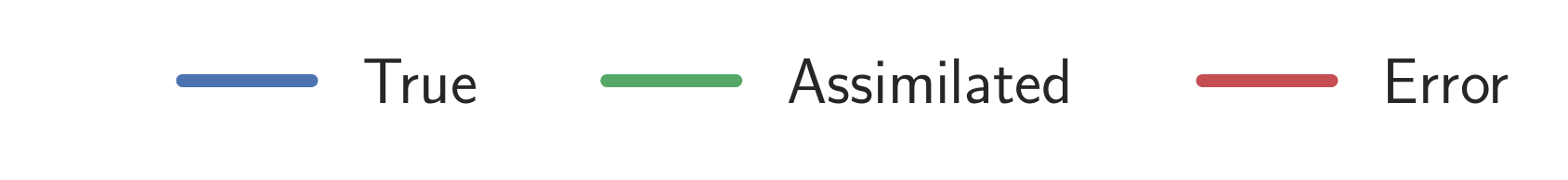}
    \caption{Results using the \textbf{time-averaging observation operator} with varying $\alpha$ regularisation parameters.
        The base setup (figure~\ref{fig:results_aneurysm_avg_noise}, left column) uses $\alpha=\gamma=10^{-5}$.
        The snapshots on the top two rows are taken at $t=0.296$s.}
    \label{fig:results_aneurysm_averaged_alpha_beta_gamma}
\end{figure}

\textbf{Regularisation (figures \ref{fig:results_aneurysm_instant_alpha_beta_gamma} and \ref{fig:results_aneurysm_averaged_alpha_beta_gamma}).}
The regularisation terms \eqref{eq:regularisation} enforce
``smoothness'' on the control functions. Hence the choice of the
regularisation coefficients $\alpha$ and $\gamma$ could have a strong influence on the
assimilation results. For the experiments, we varied $\alpha$ and $\gamma$
coefficients simultaneously to retain the balance between the two regularisation
terms. The results for different regularisation
values are shown in figures
\ref{fig:results_aneurysm_instant_alpha_beta_gamma} and
\ref{fig:results_aneurysm_averaged_alpha_beta_gamma} for $\mathcal{T}^{\text{inst}}$ and
$\mathcal{T}^{\text{avg}}$, respectively. The
reconstruction works well for values between $10^{-3}$ and
$10^{-5}$, but the quality starts to reduce visibly when $\alpha = \gamma > 10^{-2}$.
For the case $\alpha=\gamma=1$, the assimilated velocity is significantly
lower than the true velocity, because the strong regularisation enforces
spatially and temporally nearly constant controls.

\begin{figure}
    \centering
        \footnotesize
        \hspace{1cm}
    \parbox[b][8.5em][t]{0.29\textwidth}{
        \emph{Base setup with\\
        $N=32$}
        \vspace{0.5em}\\
        \input{results_aneurysm/nsassimilation_more_obs/source/results_aneurysm/instant/assimilated_H1H1_0_noise/metrics}
    }
    \parbox[b][8.5em][t]{0.29\textwidth}{
        \emph{Base setup with\\
        $N=8$}
        \vspace{0.5em}\\
        \input{results_aneurysm/nsassimilation_fewer_obs/source/results_aneurysm/instant/assimilated_H1H1_0_noise/metrics}
    }
    \parbox[b][8.5em][t]{0.29\textwidth}{
        \emph{Base setup with\\
        $N=4$}
        \vspace{0.5em}\\
        \input{results_aneurysm/nsassimilation_fewerfewerobs/source/results_aneurysm/instant/assimilated_H1H1_0_noise/metrics}
    }
    \\
    \rotatebox{90}{\qquad \quad Observation}\quad
    \includegraphics[width=0.29\textwidth]{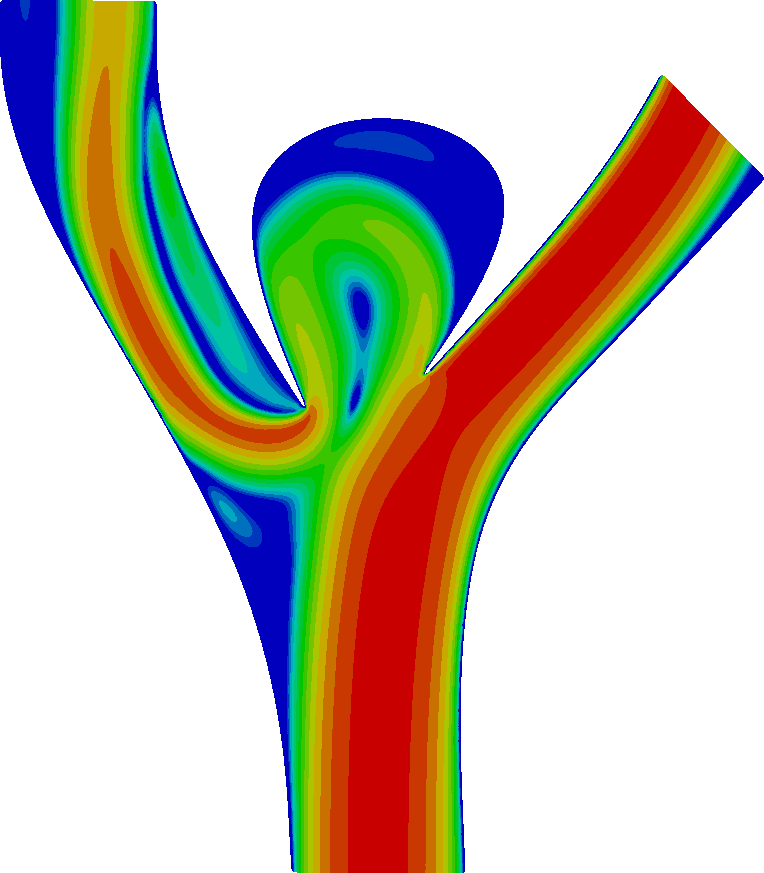}
    \includegraphics[width=0.29\textwidth]{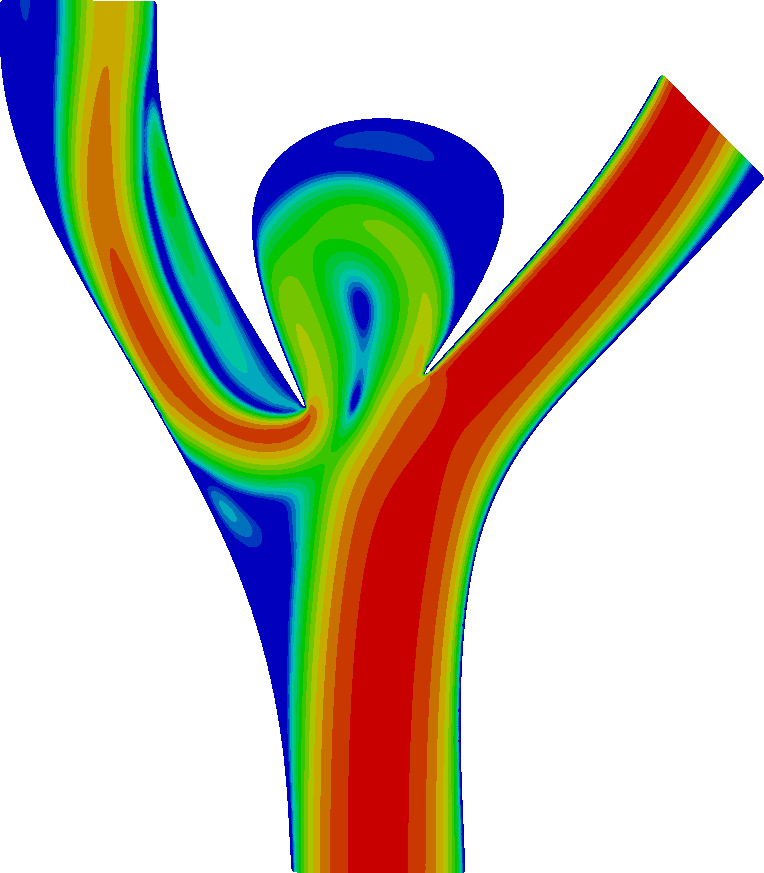}
    \includegraphics[width=0.29\textwidth]{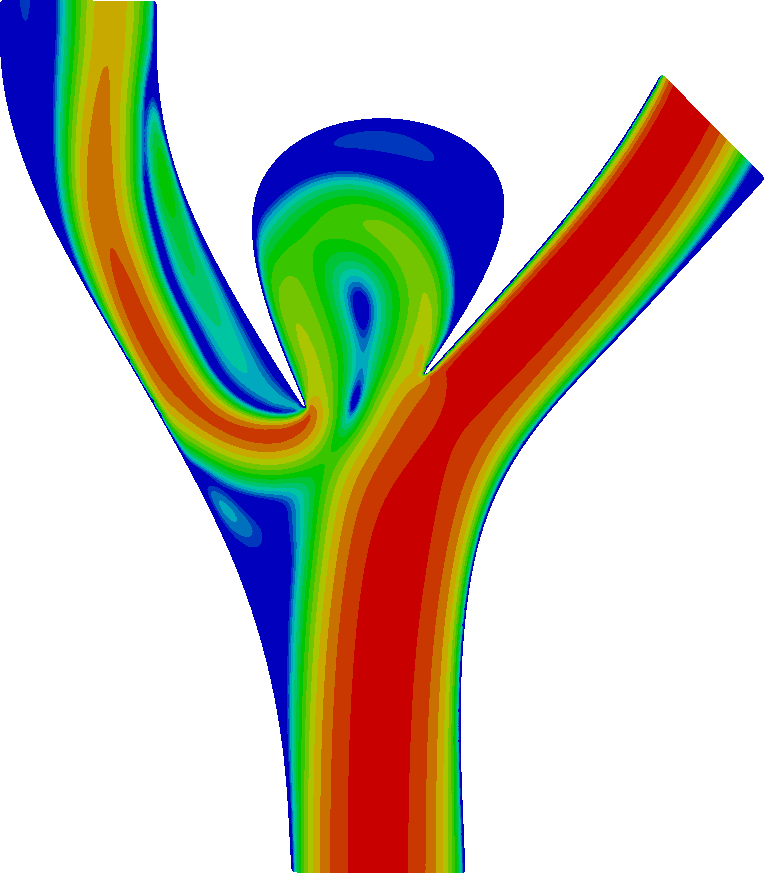}\\
    \rotatebox{90}{\qquad \quad Assimilation}\quad
    \includegraphics[width=0.29\textwidth]{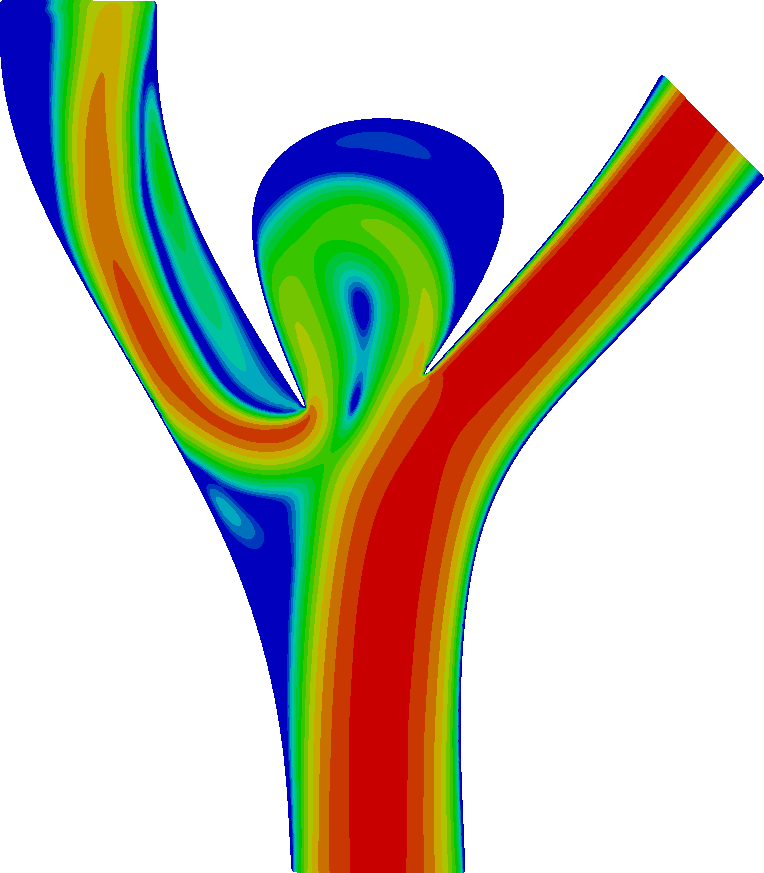}
    \includegraphics[width=0.29\textwidth]{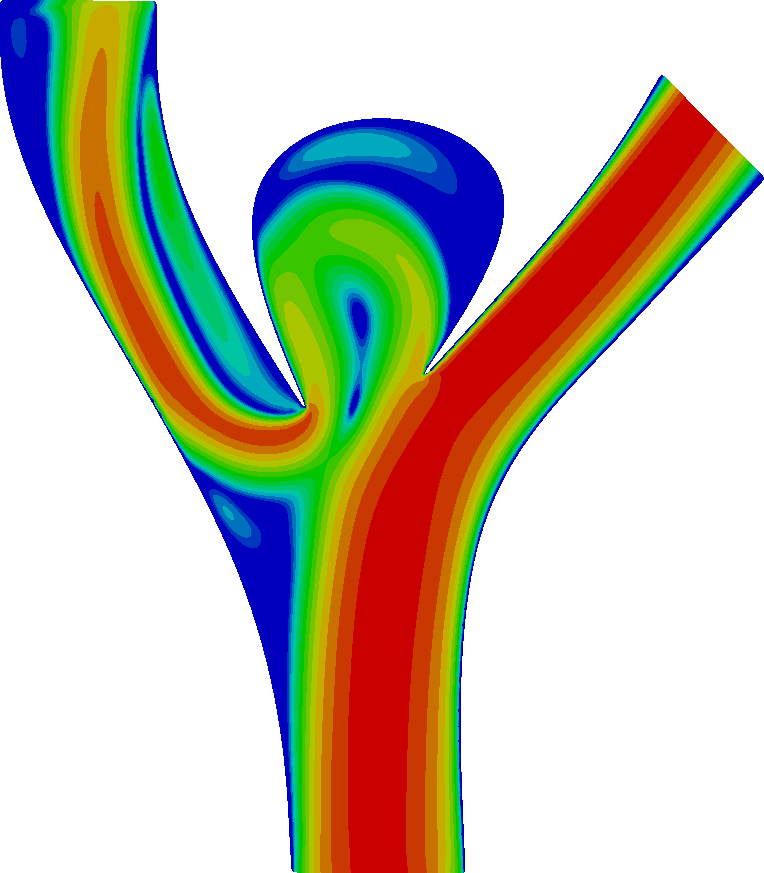}
    \includegraphics[width=0.29\textwidth]{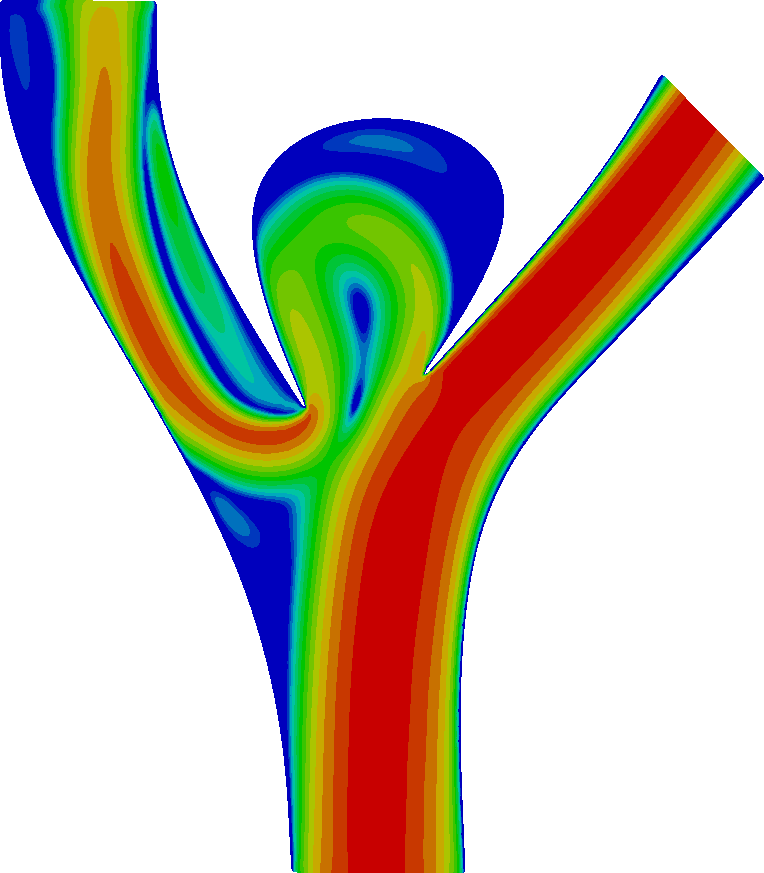}\\
    \begin{center}
        \includegraphics[width=0.3\textwidth]{results_aneurysm/scale}
    \end{center}
    \rotatebox{90}{Velocity norm in ${\Omega_{\text{ane}}}$}\quad
    \includegraphics[width=0.29\textwidth]{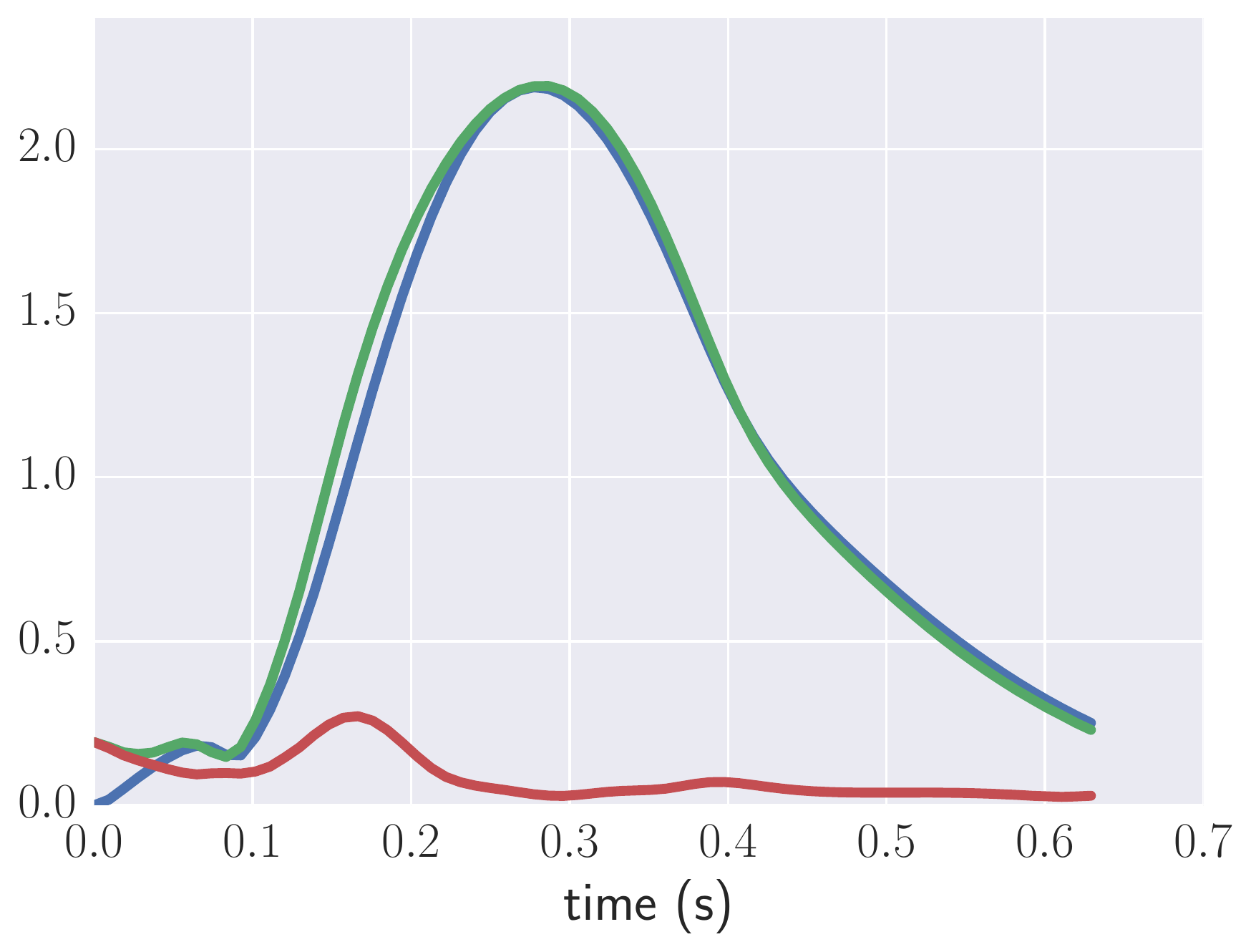}
    \includegraphics[width=0.29\textwidth]{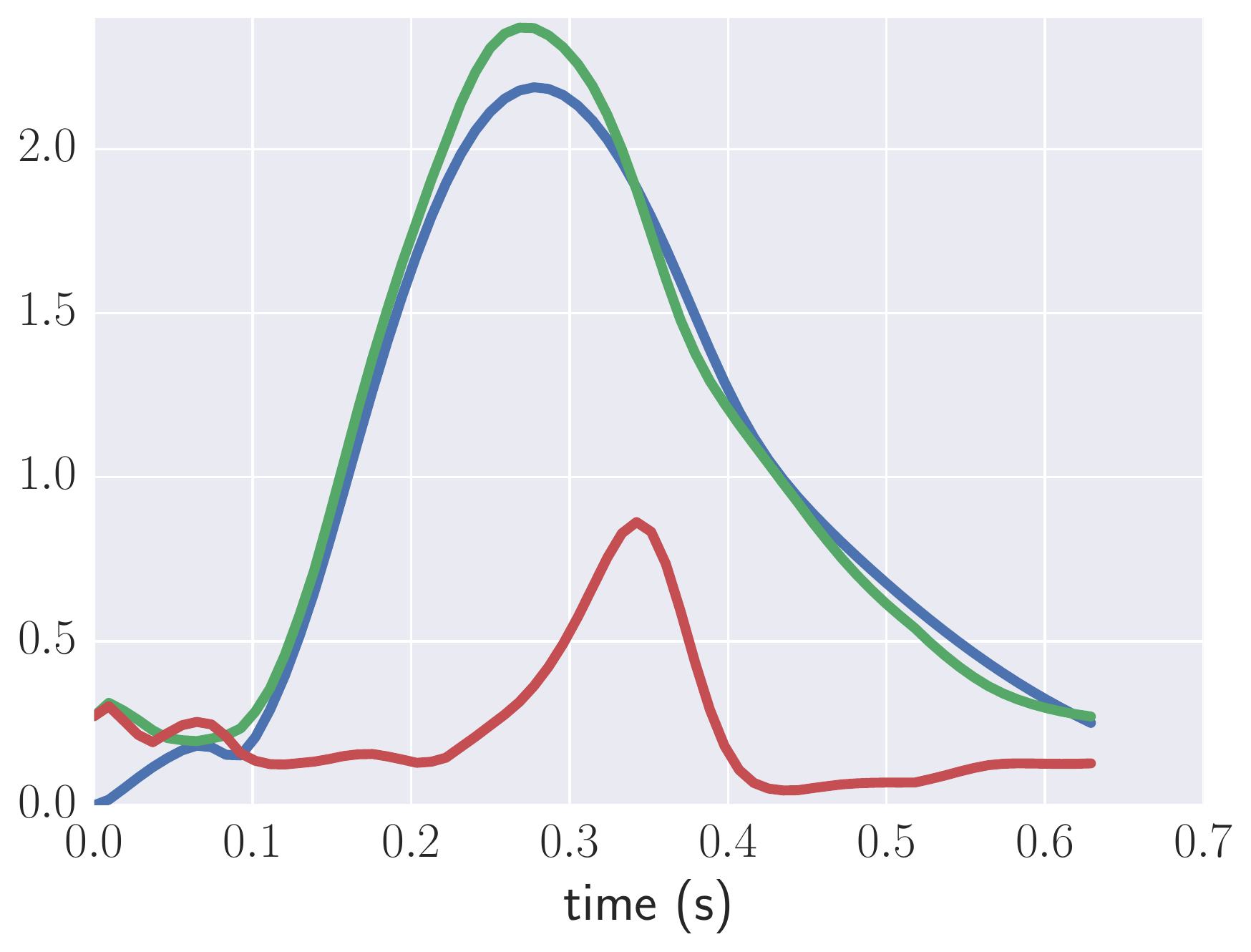}
    \includegraphics[width=0.29\textwidth]{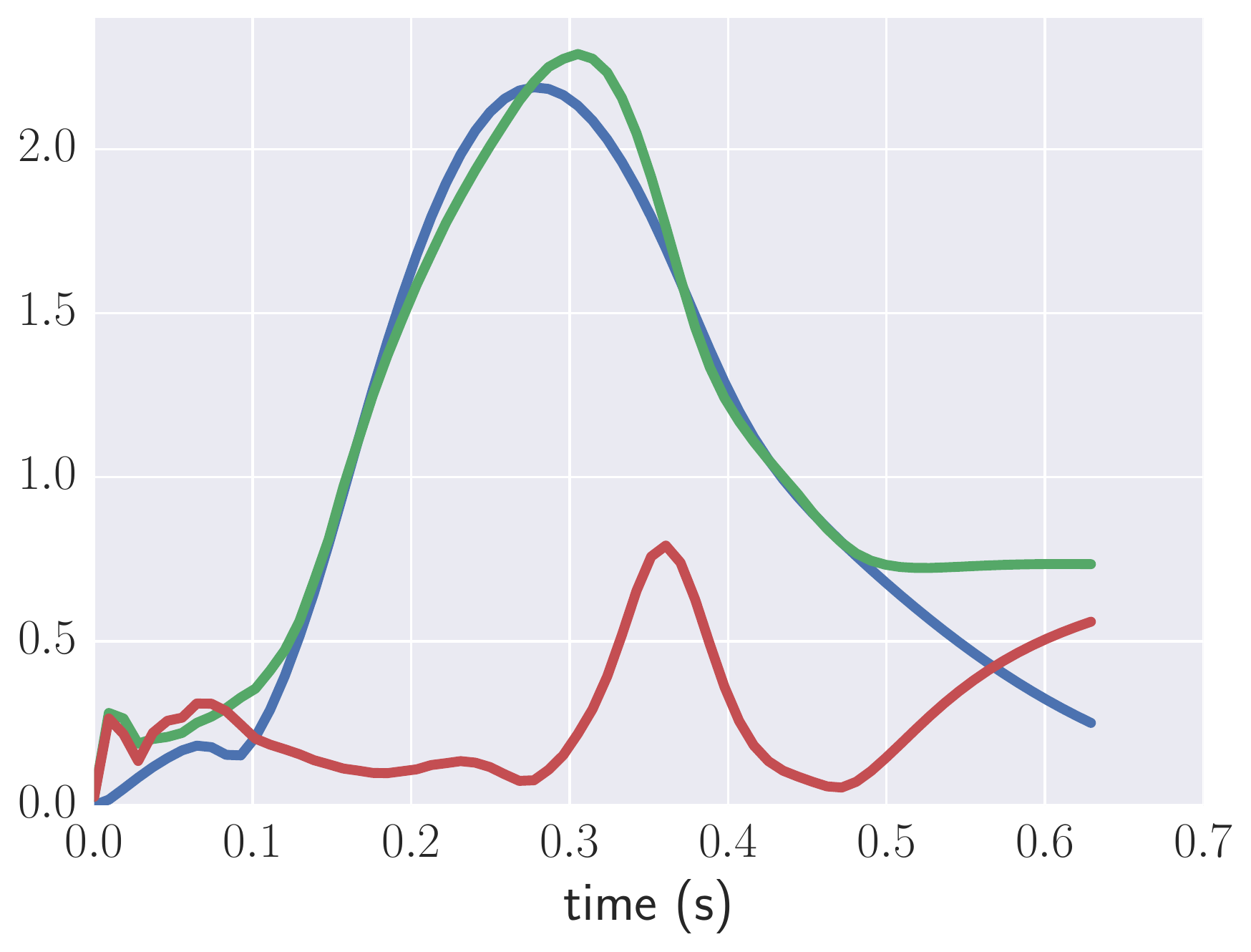}\\
    \rotatebox{90}{WSS norm on $\Gamma_{\text{ane}}$}\quad
    \includegraphics[width=0.29\textwidth]{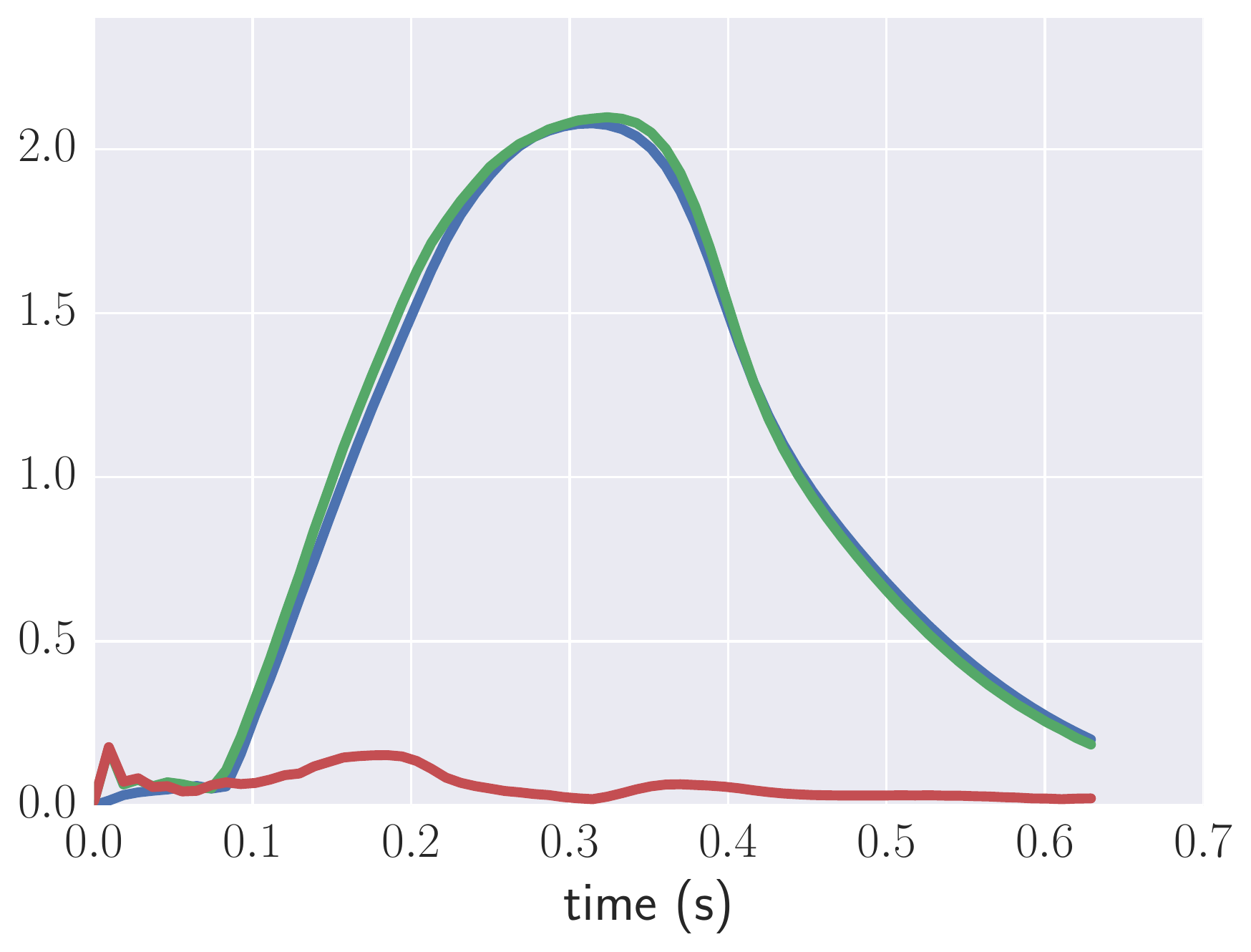}
    \includegraphics[width=0.29\textwidth]{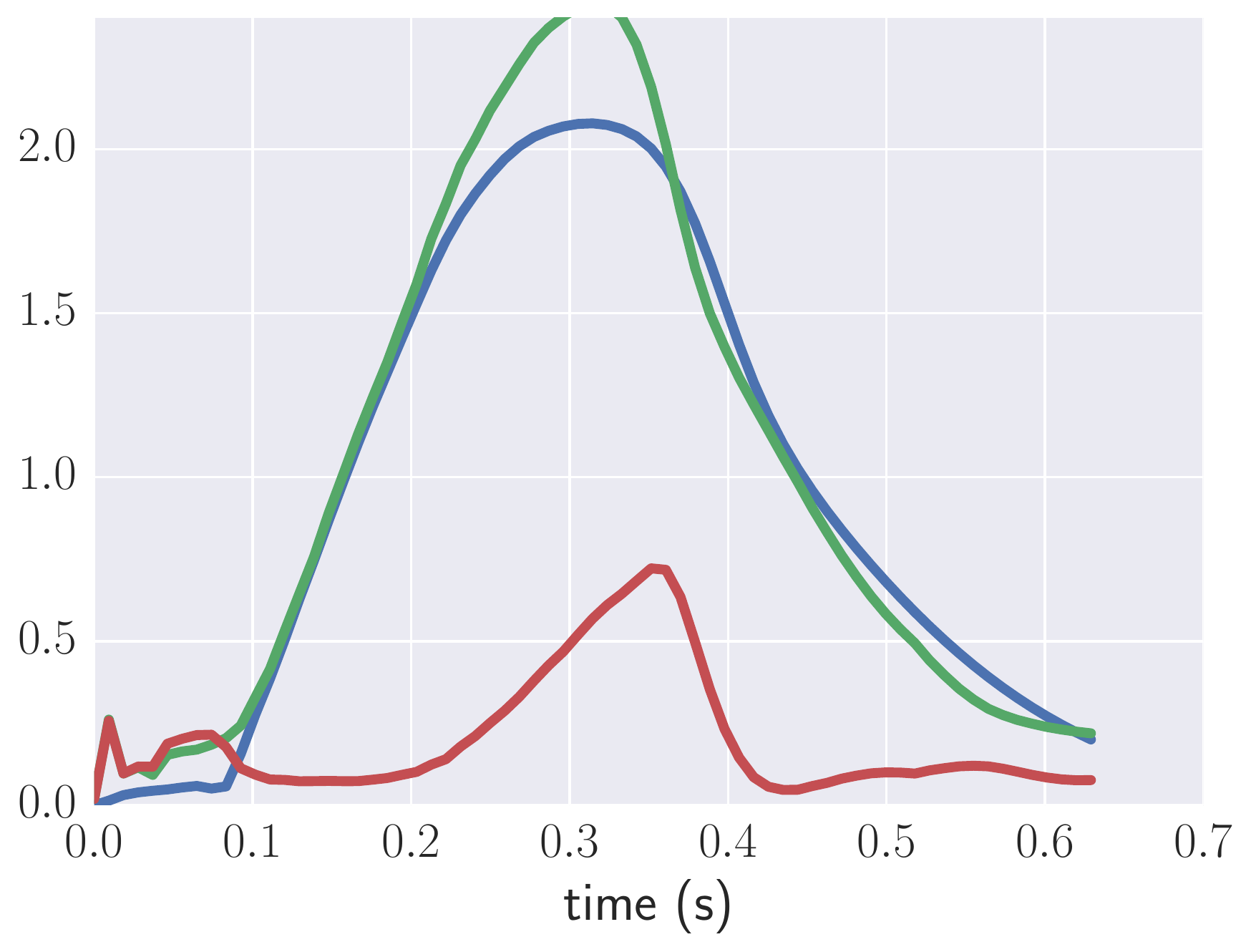}
    \includegraphics[width=0.29\textwidth]{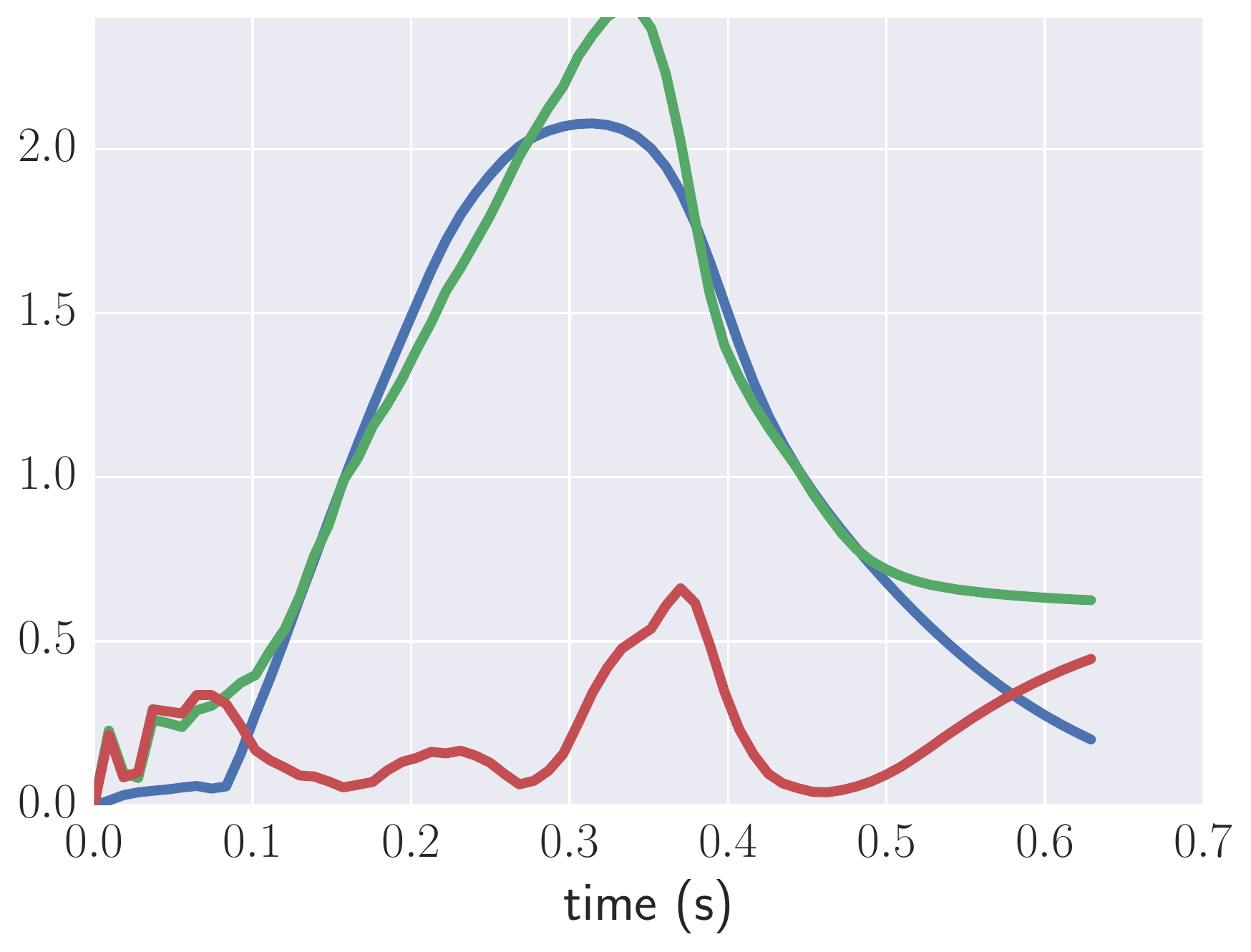}\\
    \includegraphics[width=0.59\textwidth]{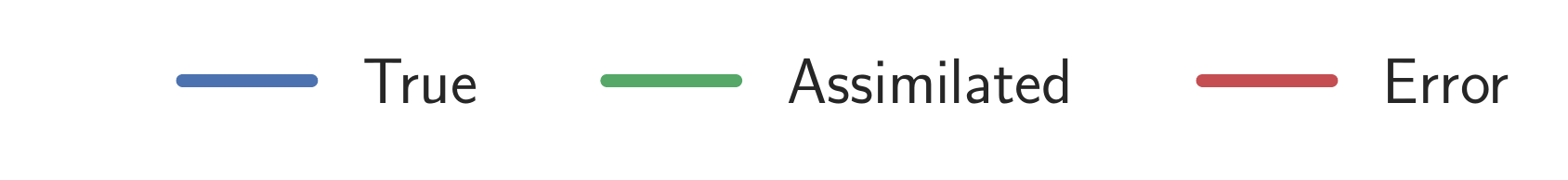}
    \caption{Results using the \textbf{instantaneous observation operator} with varying number of observations.
        The base setup (figure~\ref{fig:results_aneurysm_inst_noise}, left column) uses $N=16$ observations.
        The snapshots on the top two rows are taken at $t=0.296$s.}
    \label{fig:results_aneurysm_inst_numobs}
\end{figure}

\begin{figure}
    \centering
        \footnotesize
    \hspace{1cm}
    \parbox[b][8.5em][t]{0.29\textwidth}{
        \emph{Base setup with\\
        $N=32$}
        \vspace{0.5em}\\
        \input{results_aneurysm/nsassimilation_more_obs/source/results_aneurysm/averaged/assimilated_H1H1_0_noise/metrics}
    }
    \parbox[b][8.5em][t]{0.29\textwidth}{
        \emph{Base setup with\\
        $N=8$}
        \vspace{0.5em}\\
        \input{results_aneurysm/nsassimilation_fewer_obs/source/results_aneurysm/averaged/assimilated_H1H1_0_noise/metrics}
    }
    \parbox[b][8.5em][t]{0.29\textwidth}{
        \emph{Base setup with\\
        $N=4$}
        \vspace{0.5em}\\
        \input{results_aneurysm/nsassimilation_fewerfewerobs/source/results_aneurysm/averaged/assimilated_H1H1_0_noise/metrics}
    }
    \\
    \rotatebox{90}{\qquad \quad Observation}\quad
    \includegraphics[width=0.29\textwidth]{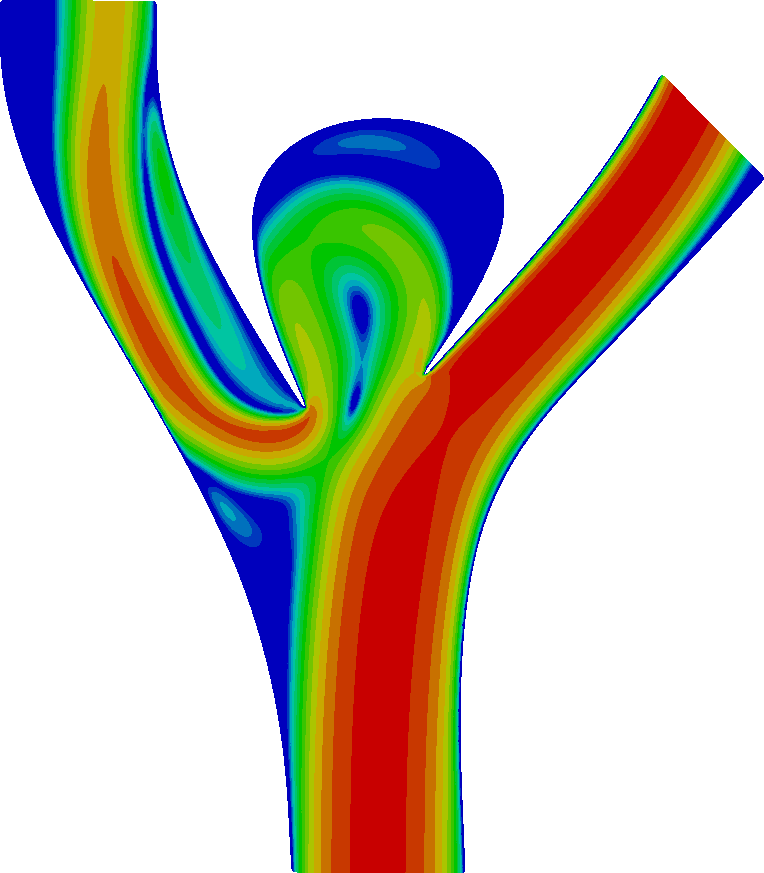}
    \includegraphics[width=0.29\textwidth]{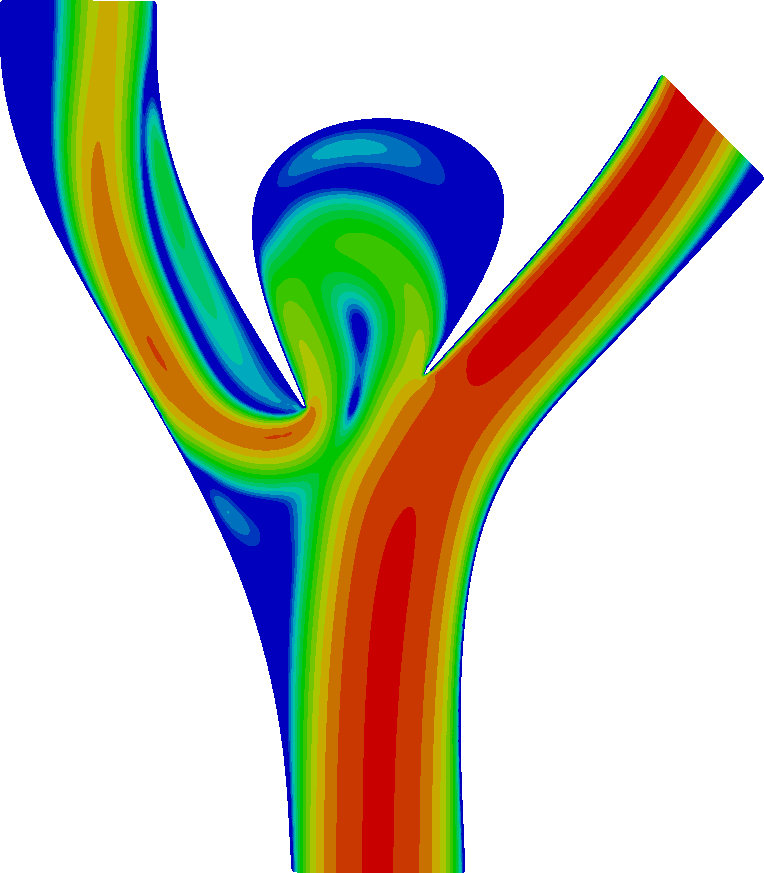}
    \includegraphics[width=0.29\textwidth]{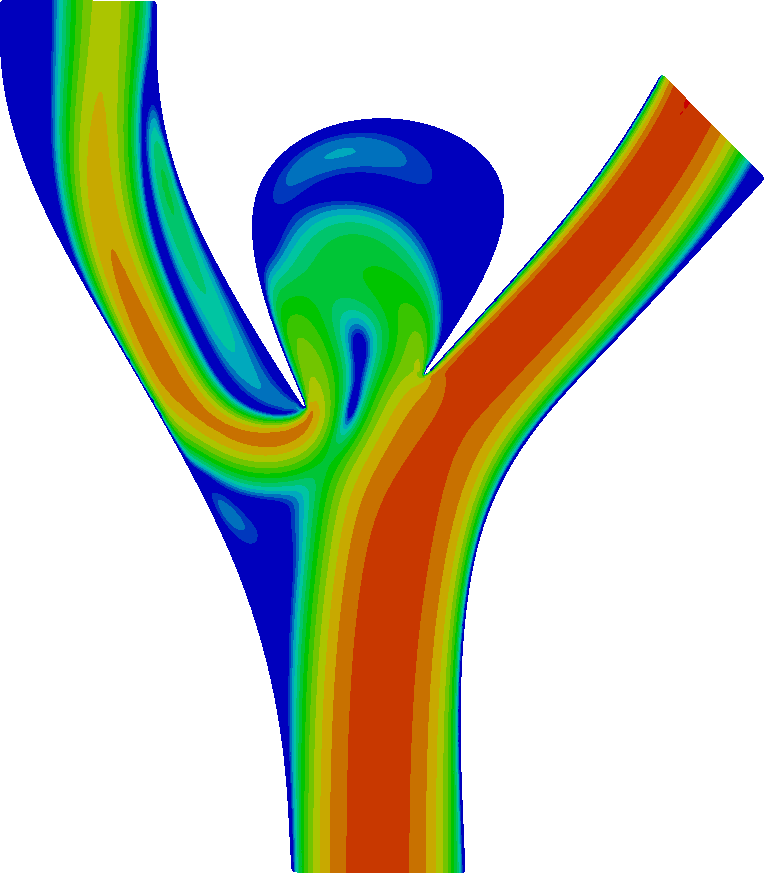}\\
    \rotatebox{90}{\qquad \quad Assimilation}\quad
    \includegraphics[width=0.29\textwidth]{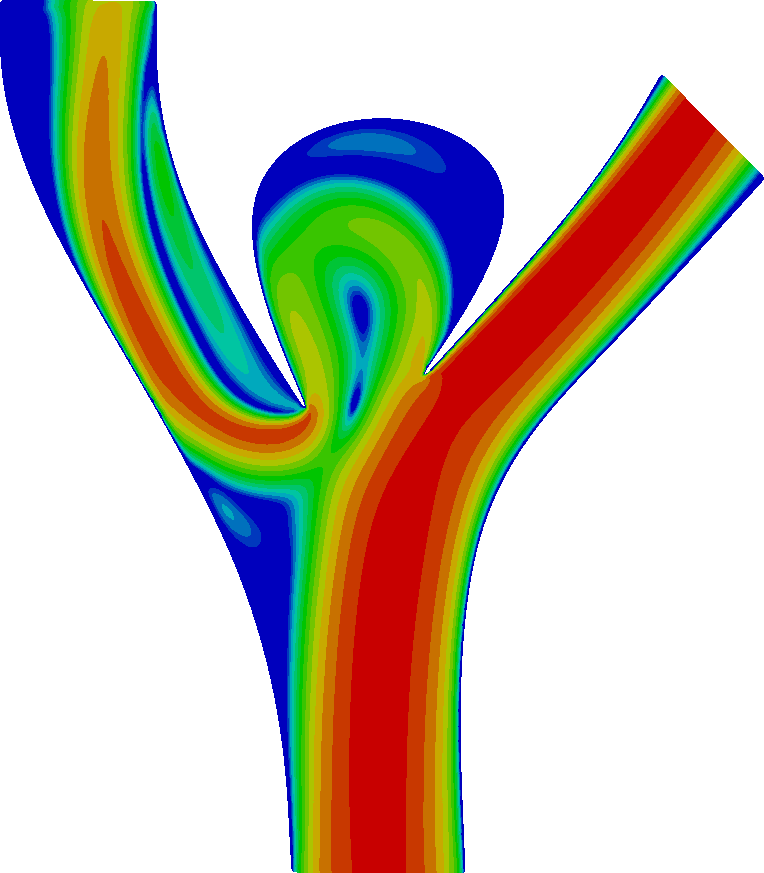}
    \includegraphics[width=0.29\textwidth]{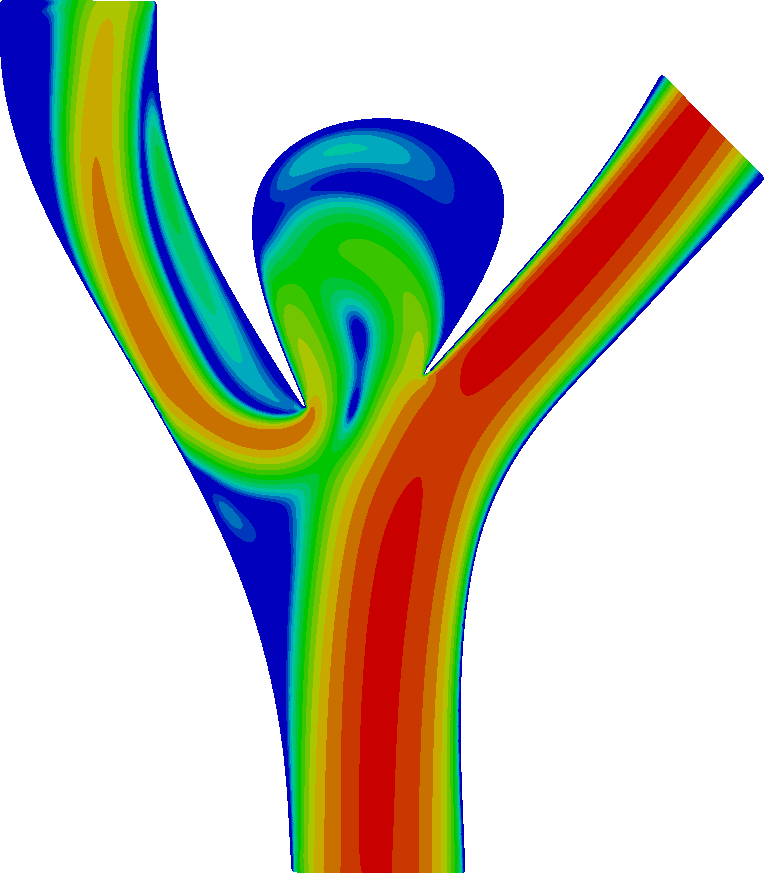}
    \includegraphics[width=0.29\textwidth]{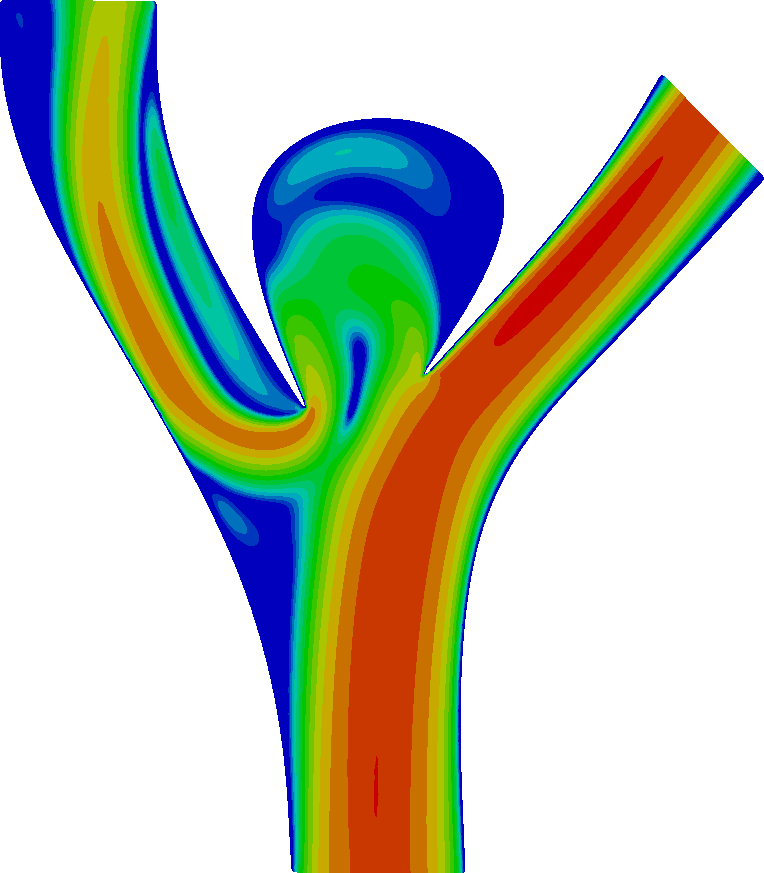}\\
    \begin{center}
        \includegraphics[width=0.3\textwidth]{results_aneurysm/scale}
    \end{center}
    \rotatebox{90}{Velocity norm in ${\Omega_{\text{ane}}}$}\quad
    \includegraphics[width=0.29\textwidth]{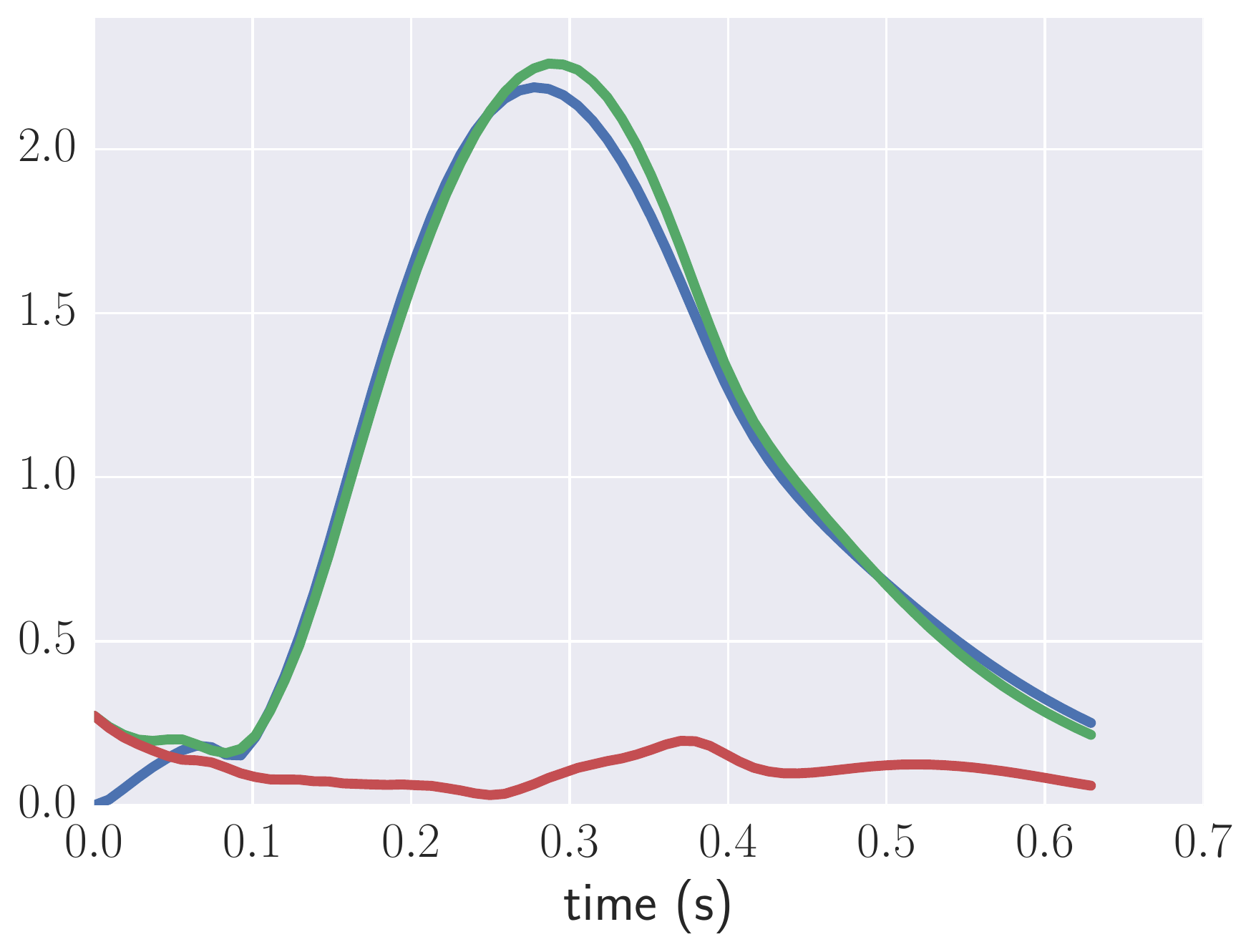}
    \includegraphics[width=0.29\textwidth]{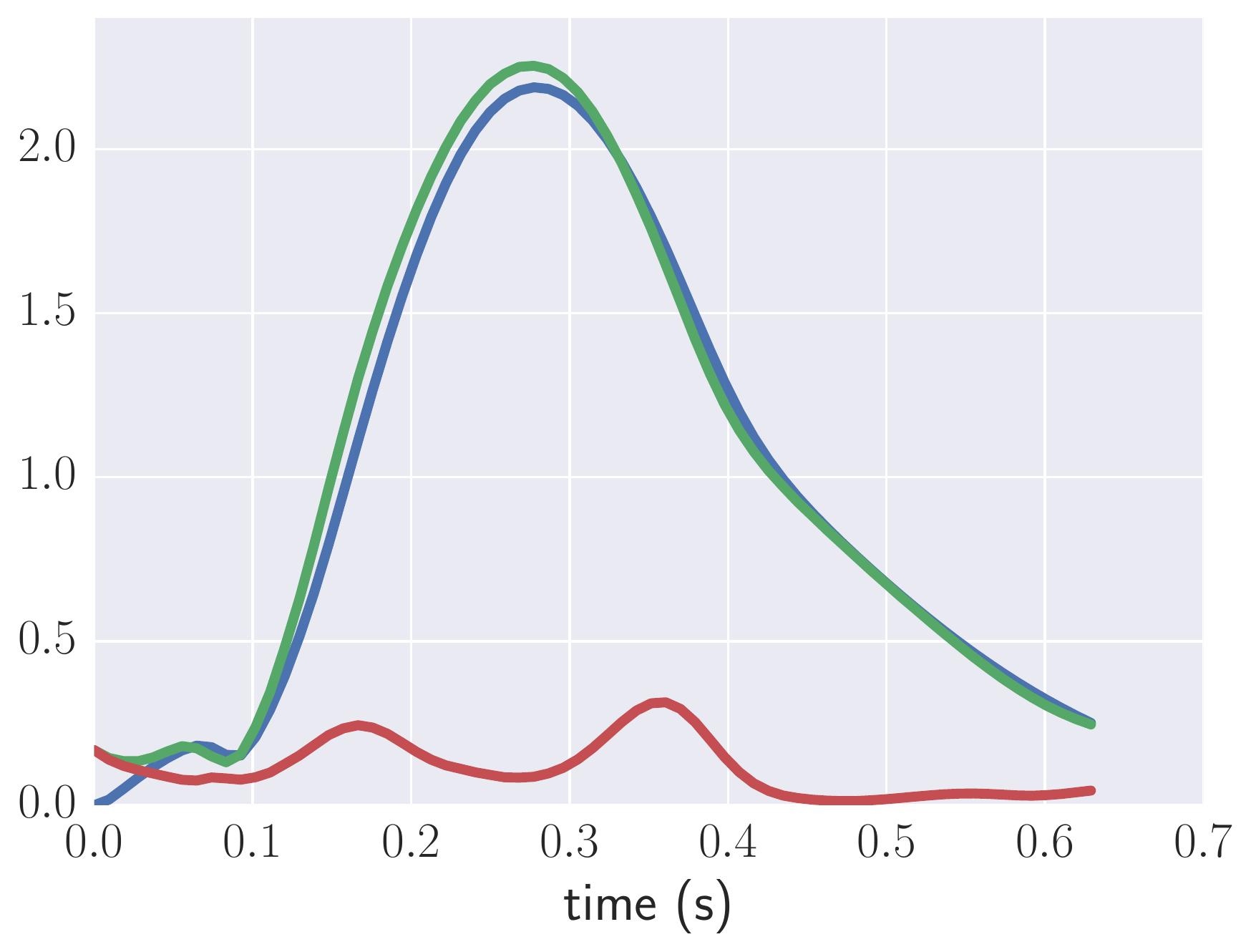}
    \includegraphics[width=0.29\textwidth]{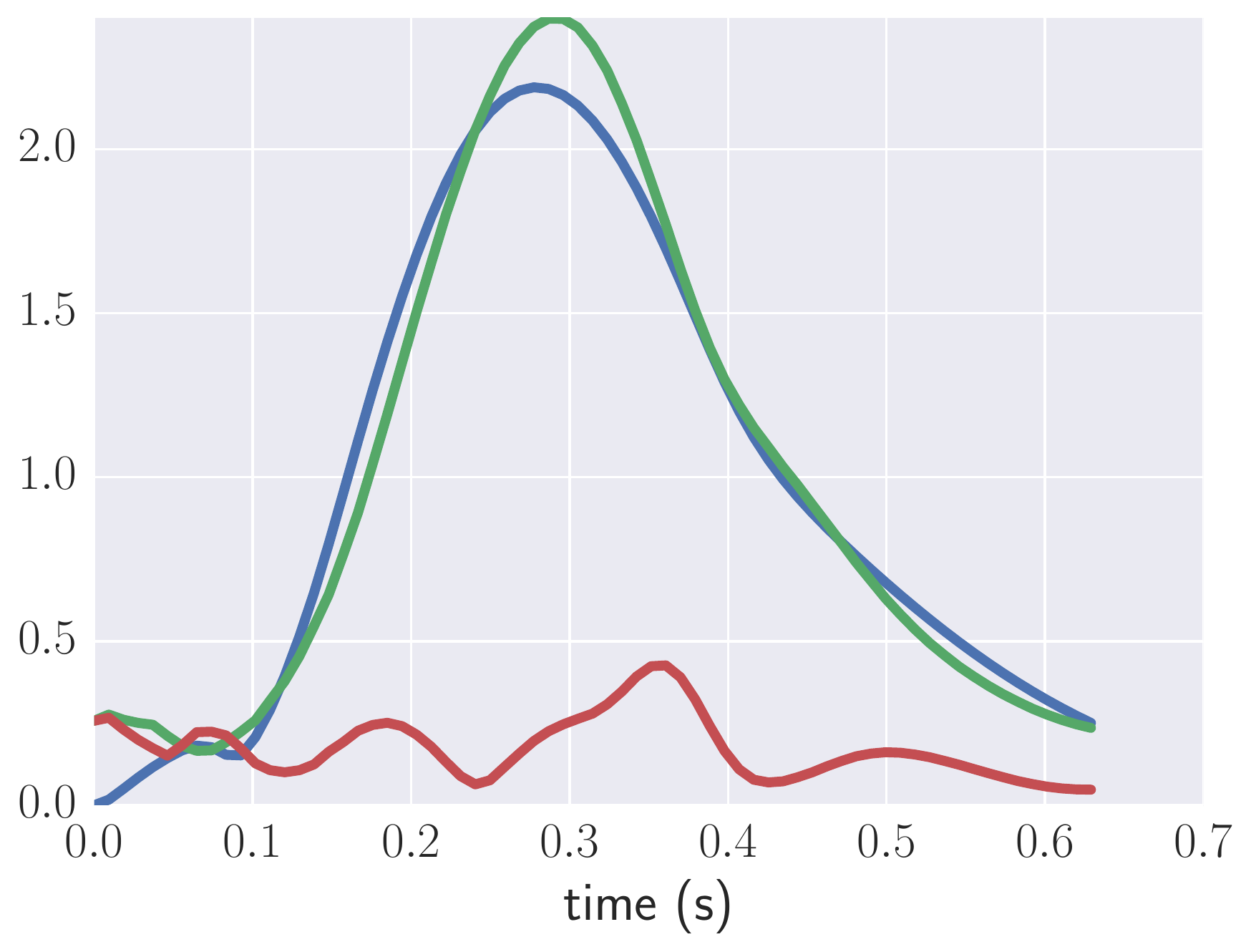}\\
    \rotatebox{90}{WSS norm on $\Gamma_{\text{ane}}$}\quad
    \includegraphics[width=0.29\textwidth]{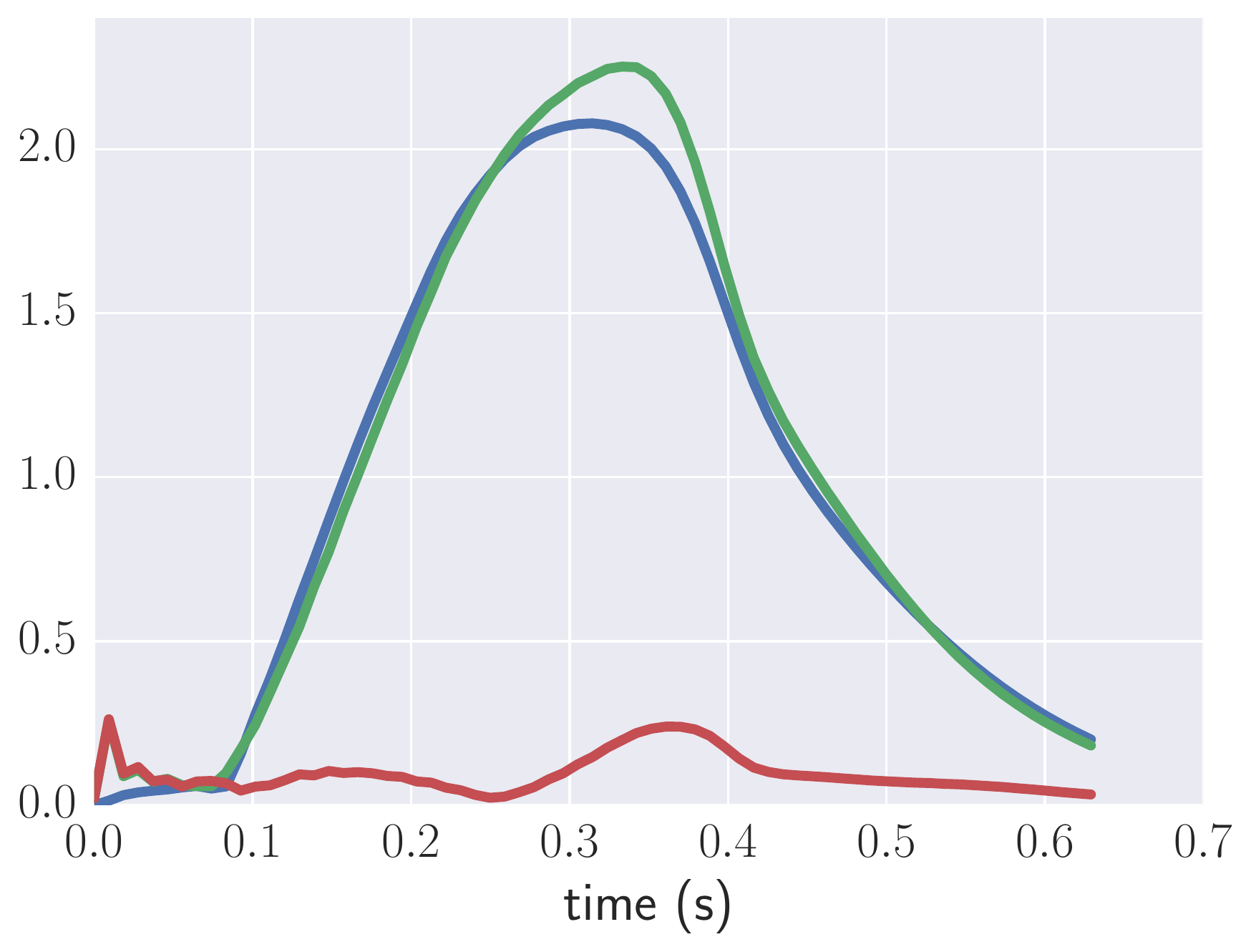}
    \includegraphics[width=0.29\textwidth]{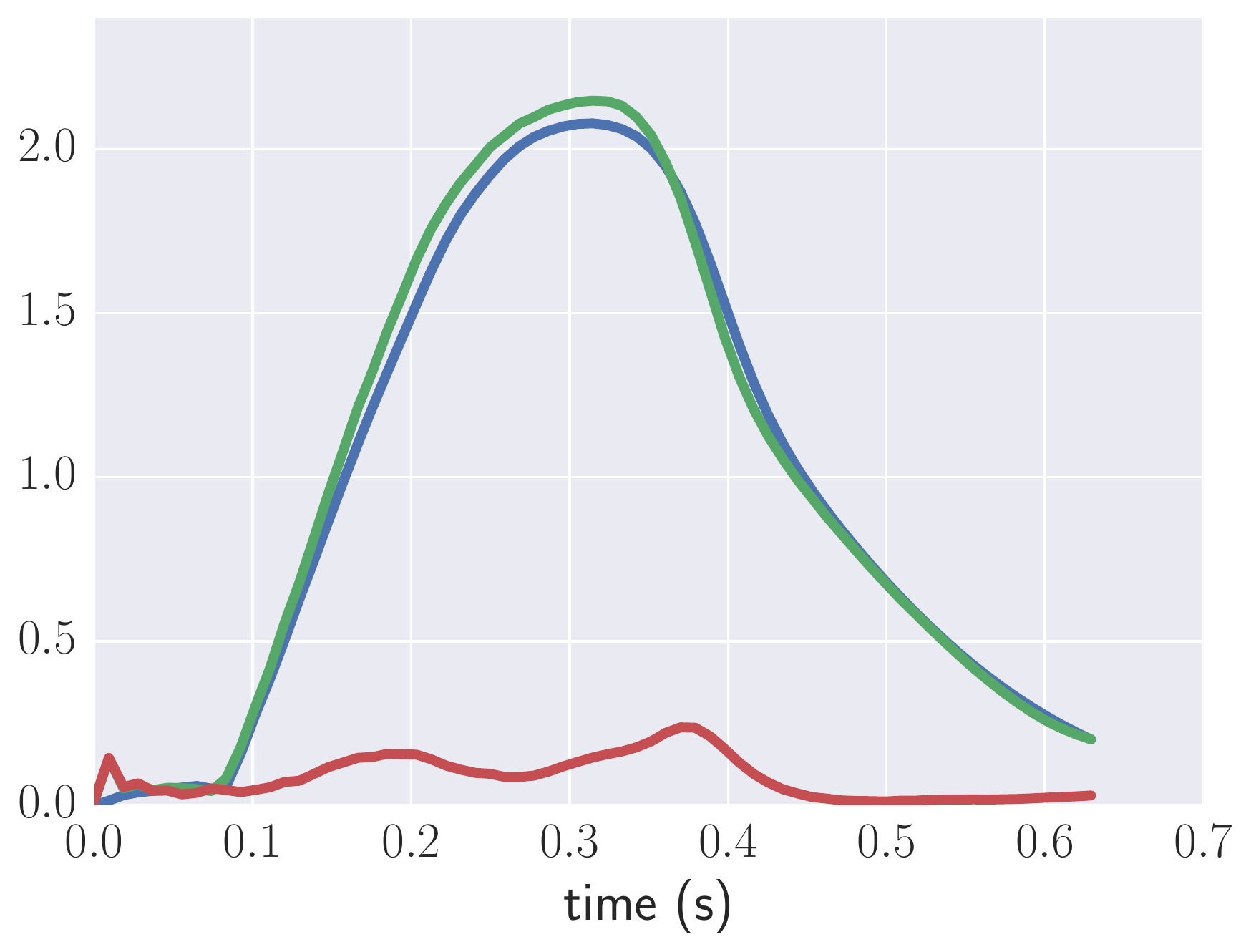}
    \includegraphics[width=0.29\textwidth]{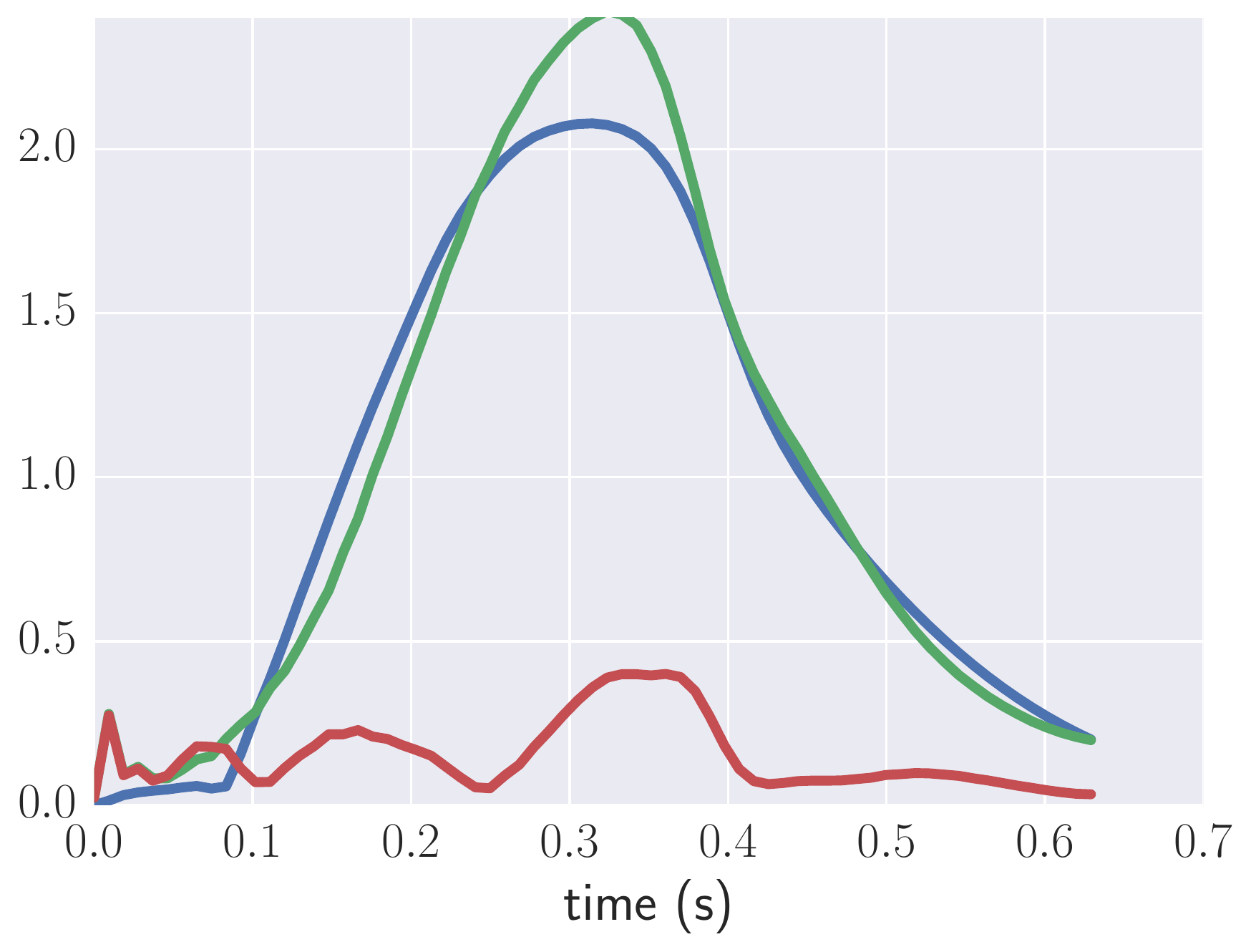}\\
    \includegraphics[width=0.59\textwidth]{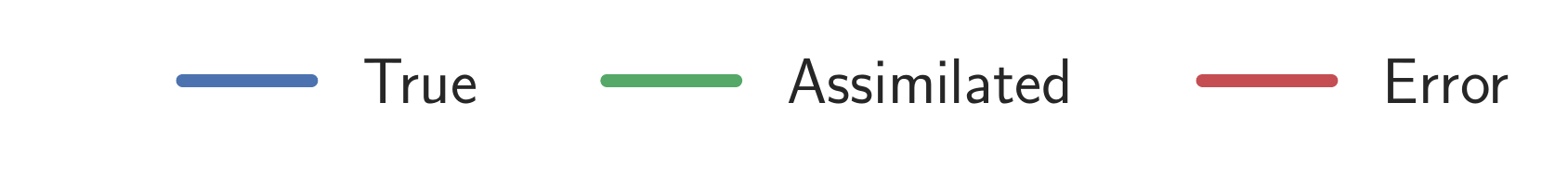}
    \caption{Results using the \textbf{time-averaging observation operator} with varying number of observations.
        The base setup (figure~\ref{fig:results_aneurysm_avg_noise}, left column) uses $N=16$ observations.
        The snapshots on the top two rows are taken at $t=0.296$s.}
    \label{fig:results_aneurysm_avg_numobs}
\end{figure}

\textbf{Data sparsity (figures \ref{fig:results_aneurysm_inst_numobs} and \ref{fig:results_aneurysm_avg_numobs}).}
Another important question is how many observations ($N$ in \eqref{eq:goal_quantity}) are required to accurately reconstruct
the blood flow. To address this question, the data assimilation was repeated with varying number of observations $N$.
The base setup (left column in figures \ref{fig:results_aneurysm_inst_noise} and \ref{fig:results_aneurysm_avg_noise}) used $N=16$ and the results for $N=4, 8$ and $32$ are shown in
figures
\ref{fig:results_aneurysm_inst_numobs} and
\ref{fig:results_aneurysm_avg_numobs}
for $\mathcal{T}^{\text{inst}}$ and $\mathcal{T}^{\text{avg}}$, respectively.

With $4$ observations the quality of the reconstruction suffers visibly,
mostly at the beginning and the end of the simulation times. The time-averaging observation
operator yields good results already with $8$ observations, while the instantaneous observation operator
requires $16$ observations to yield an accurate reconstruction. The differences between
$16$ to $32$ observations are minimal for both observation operators.

\begin{figure}
    \centering
        \footnotesize
    \parbox[b][9.5em][t]{0.29\textwidth}{
        \emph{Base setup with\\
        swapped Dirichlet control
        \\and $T^{\text{inst}}$}
        \vspace{0.5em}\\
        \input{results_aneurysm/nsassimilation_fewerfewerobs/source/results_aneurysm/averaged/assimilated_H1H1_0_noise/metrics}
    }
    \parbox[b][9.5em][t]{0.29\textwidth}{
        \emph{Base setup with\\
        swapped Dirichlet \\
        control and $T^{\text{avg}}$}
        \vspace{0.5em}\\
        \input{results_aneurysm/nsassimilation_fewer_obs/source/results_aneurysm/averaged/assimilated_H1H1_0_noise/metrics}
    }
    \\
    \rotatebox{90}{\qquad \quad Observation}\quad
    \includegraphics[width=0.29\textwidth]{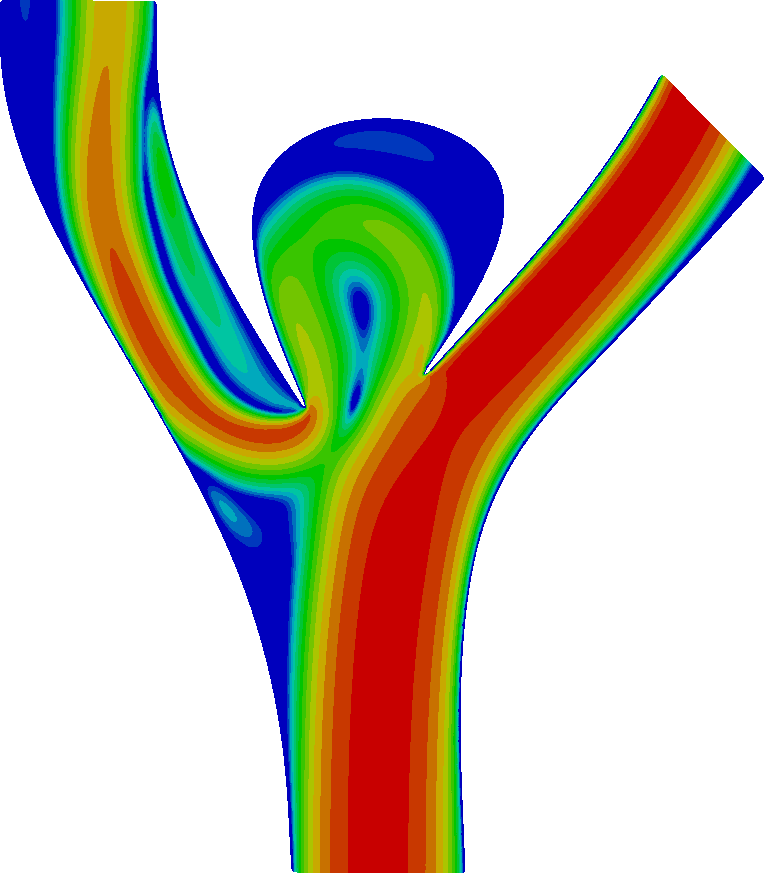}
    \includegraphics[width=0.29\textwidth]{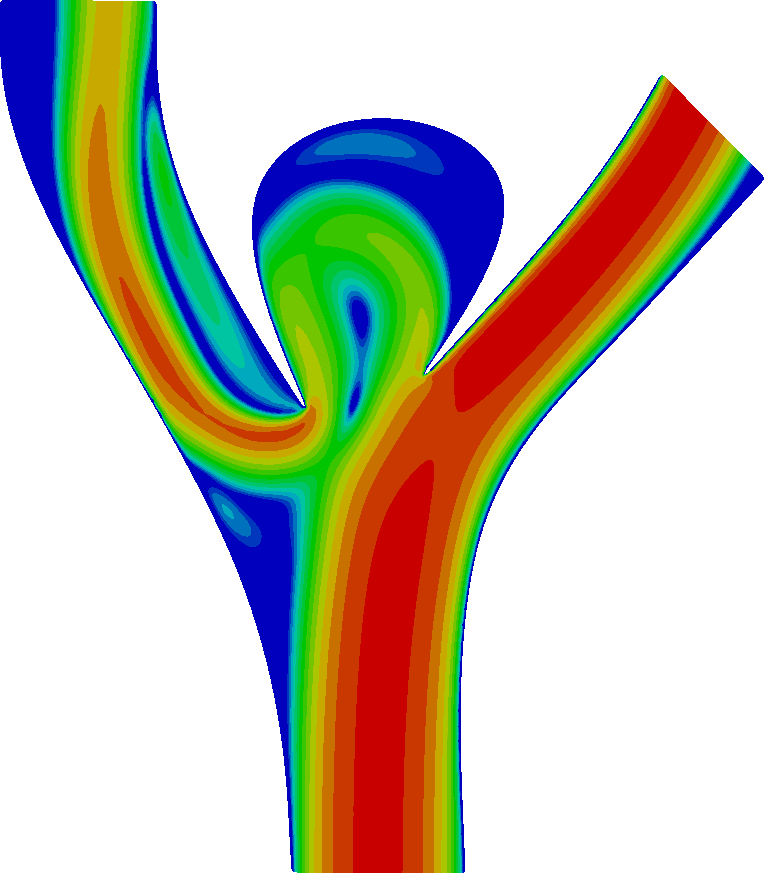}\\
    \rotatebox{90}{\qquad \quad Assimilation}\quad
    \includegraphics[width=0.29\textwidth]{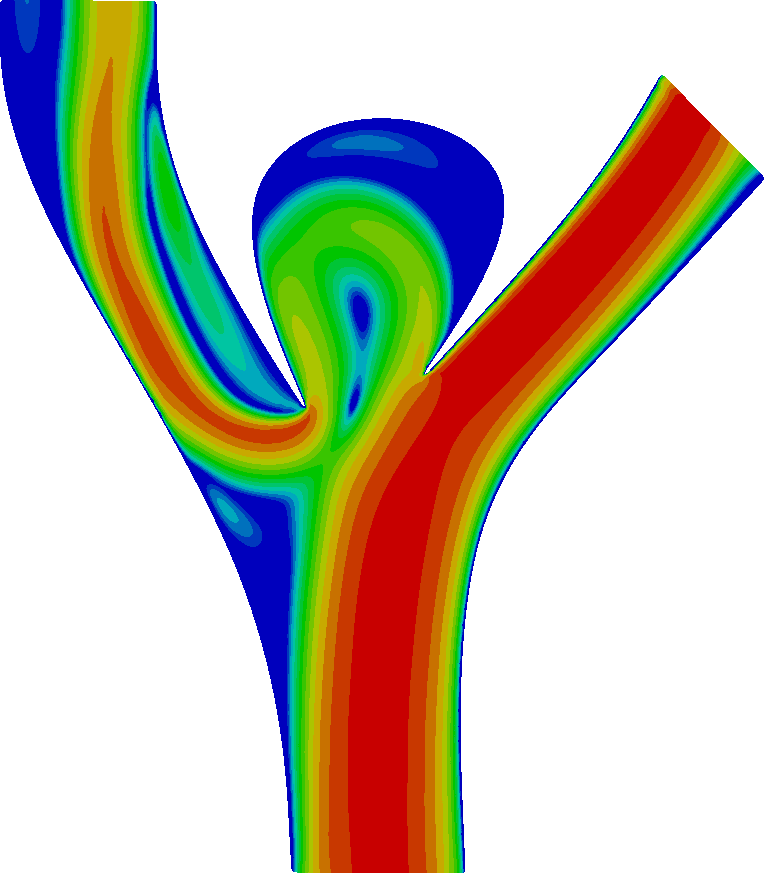}
    \includegraphics[width=0.29\textwidth]{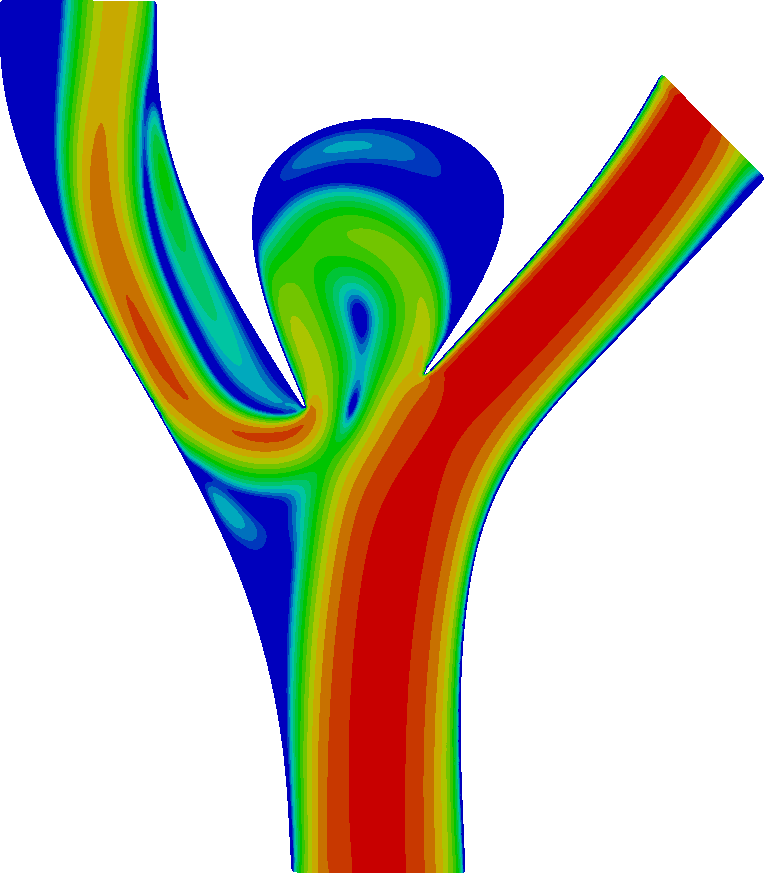}\\
    \begin{center}
        \includegraphics[width=0.3\textwidth]{results_aneurysm/scale}
    \end{center}
    \rotatebox{90}{Velocity norm in ${\Omega_{\text{ane}}}$}\quad
    \includegraphics[width=0.29\textwidth]{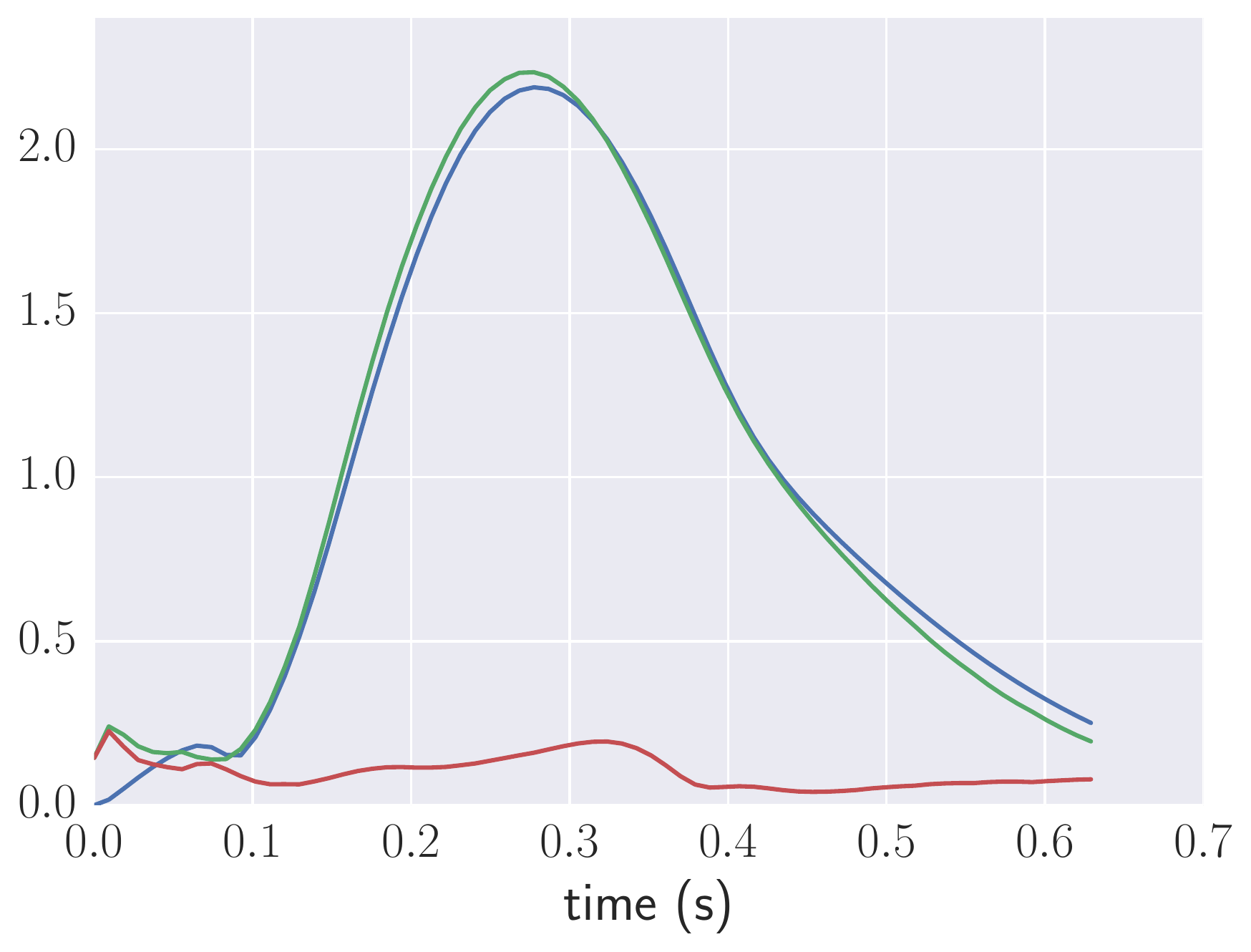}
    \includegraphics[width=0.29\textwidth]{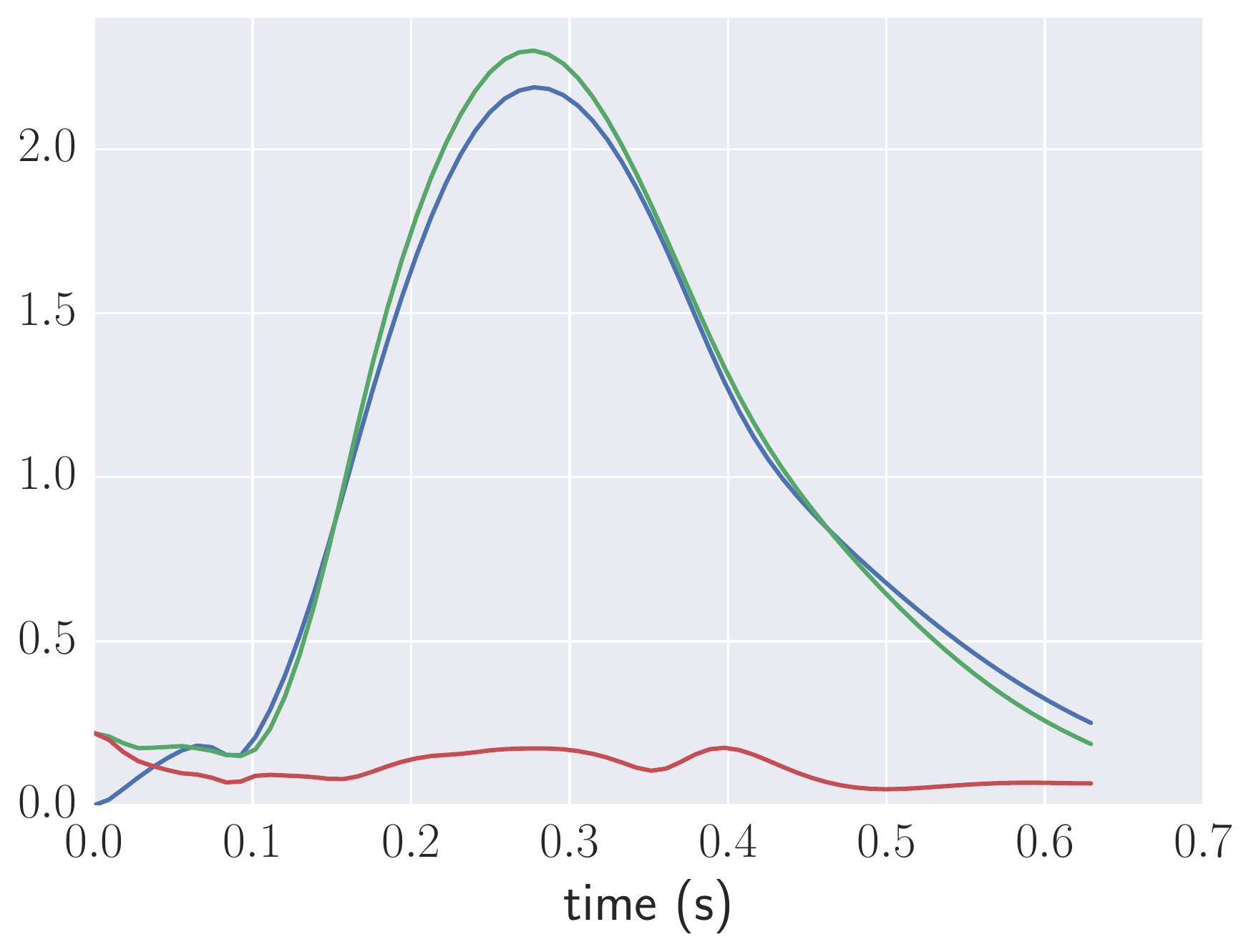} \\
    \rotatebox{90}{WSS norm on $\Gamma_{\text{ane}}$}\quad
    \includegraphics[width=0.29\textwidth]{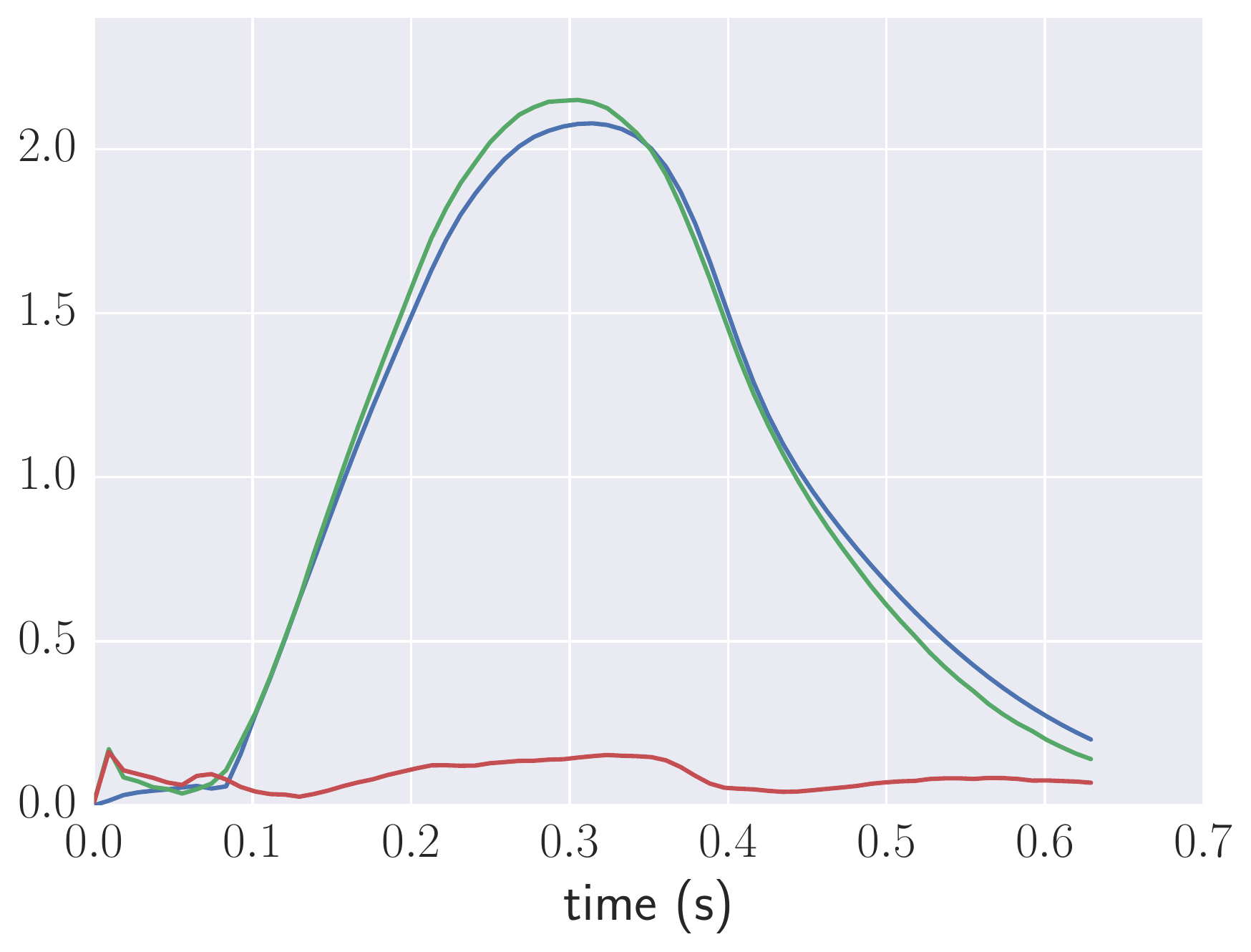}
    \includegraphics[width=0.29\textwidth]{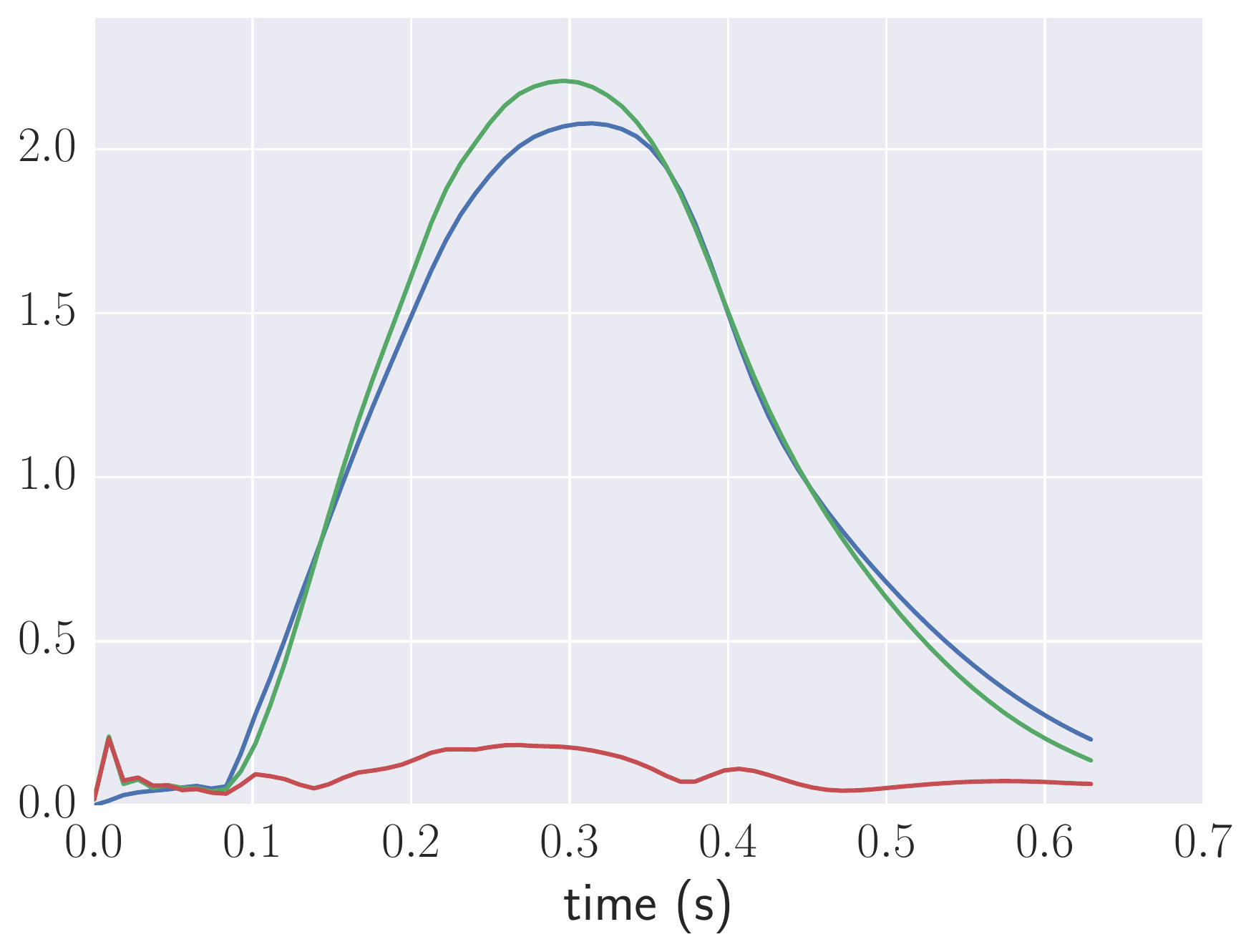} \\
    \includegraphics[width=0.59\textwidth]{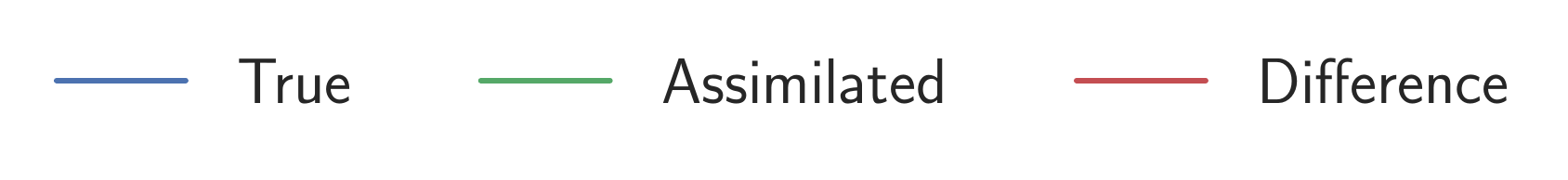}
    \caption{Assimilation results where the controlled outlets are swapped from the base cases
        (figures~\ref{fig:results_aneurysm_inst_noise} and \ref{fig:results_aneurysm_avg_noise}, left columns).
        The snapshots on the top two rows are taken at $t=0.296$s.}
    \label{fig:results_aneurysm_swapped}
\end{figure}

\textbf{Choice of controlled outflow boundary (figure \ref{fig:results_aneurysm_swapped}).}
In the problem definition \eqref{eq:ns_bc_dirichlet}, we made a choice to control the outflow on $\Gamma_{\text{out}_1}$, and to
enforce a no-stress condition on $\Gamma_{\text{out}_2}$.
It is therefore natural to check if the reconstruction works well also if $\Gamma_{\text{out}_2}$ is controlled and a non-stress condition is applied on $\Gamma_{\text{out}_1}$.
The results for this setup are shown in figure
\ref{fig:results_aneurysm_swapped}.
The reconstruction is similarly good as in the
base setup, indicating that the assimilation is not impacted significantly
by the choice of the controlled boundaries.
Nevertheless, this choice might be more significant for other setups, in particular
if one of the outflows is in close proximity to the aneurysm.

\subsection{Flow reconstruction in an aneurysm from 4D MRA measurements}\label{sec:example_3d}

In this experiment, the variational data assimilation was applied to reconstruct the blood flow
conditions in an artificially introduced aneurysm. The measurements were done using 4D PC-MRA, and are described in detail in \cite{jiang2011}. The geometry was reconstructed from an image obtained by time-averaging the observations, using VMTK (\url{www.vmtk.org}). The observations were then linearly interpolated onto the resulting mesh nodes.

The numerical settings are listed in
table~\ref{tab:3d_setup}. The assimilation terminated after 30 optimisation
iterations, when the relative change in one optimisation iteration dropped
below $0.09 \%$.
\begin{table}
    \centering
    \begin{tabular}{ l c c }
        Parameter & Symbol & Value\\
        \hline
        Viscosity & $\nu$ & $7.5$\\
        Model timestep  & $\Delta t$ & $0.004625$ s \\
        End time  & $T$ & $0.629$ s \\
        Time discretisation & $\theta$ & $1.0$ \\
        Spatial discretisation &  & $P1$-$P1$ \\
        Dimension of spatial discretisation &  & 184,464 \\
        Nitsche coefficient & $\sigma$ & $100.0$ \\
        \hline
        Number of observations & N & $16$ \\
        Regularisation parameter & $\alpha = \gamma$ & $10^{-5}$ \\
        \hline
    \end{tabular}
    \caption{The numerical settings for the reconstruction of blood flow from 4D MRA measurements.
    The first parameters specify the model setup, while the final two parameters configure the data assimilation.}
    \label{tab:3d_setup}
\end{table}

For comparison, we additionally performed a
high-resolution, low-viscosity flow simulation with ``common'' choices for
the initial and boundary conditions, but without the data assimilation procedure described in this paper.
To avoid spurious effects near the boundary for this simulation,
the segmented geometry needed to be extended by artificial straight arteries on the in-
and outlets. The high-resolution setup had 20 million DOFs with a Taylor-Hood pressure-corrector scheme and a timestep of $5.9 \cdot 10^{-4}$ s.
Womersley boundary conditions were used on the inflow and outflows. The inflow flux $Q_{\Gamma_{\text{in}}}$ was interpolated from averaged observations as
$(Q_{\Gamma_\text{in}}-Q_{\Gamma_{\text{out}_1}}-Q_{\Gamma_{\text{out}_2}})/3$, the outflow flux on $\Gamma_{\text{out}_1}$ was
averaged as $Q_{\Gamma_{\text{out}_1}}/(Q_{\Gamma_{\text{out}_1}}+Q_{\Gamma_{\text{out}_2}})$, and a traction free condition was applied on $\Gamma_{\text{out}_2}$.

\begin{figure}
    \centering
        \footnotesize
    \rotatebox{90}{\hspace{2.5cm} Observation}\quad
    \includegraphics[width=0.45\textwidth]{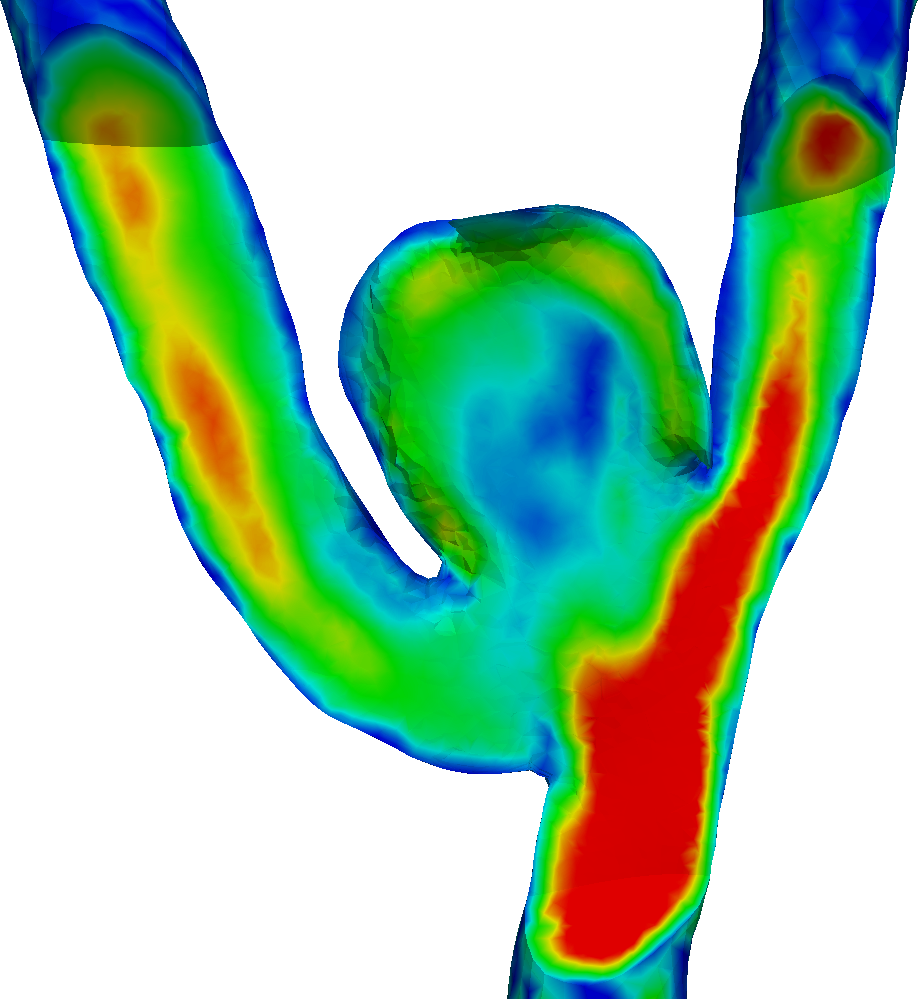}
    \includegraphics[width=0.45\textwidth]{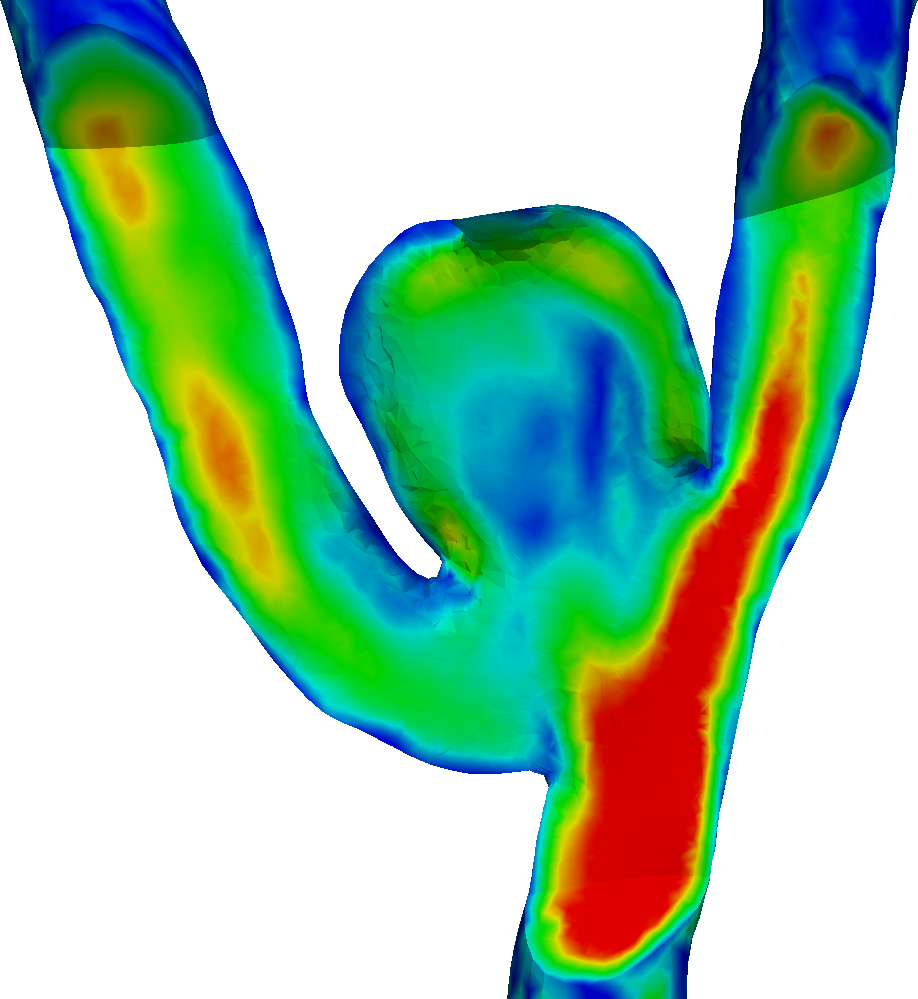}\\
    \rotatebox{90}{\hspace{2.5cm} Assimilation}\quad
    \includegraphics[width=0.45\textwidth]{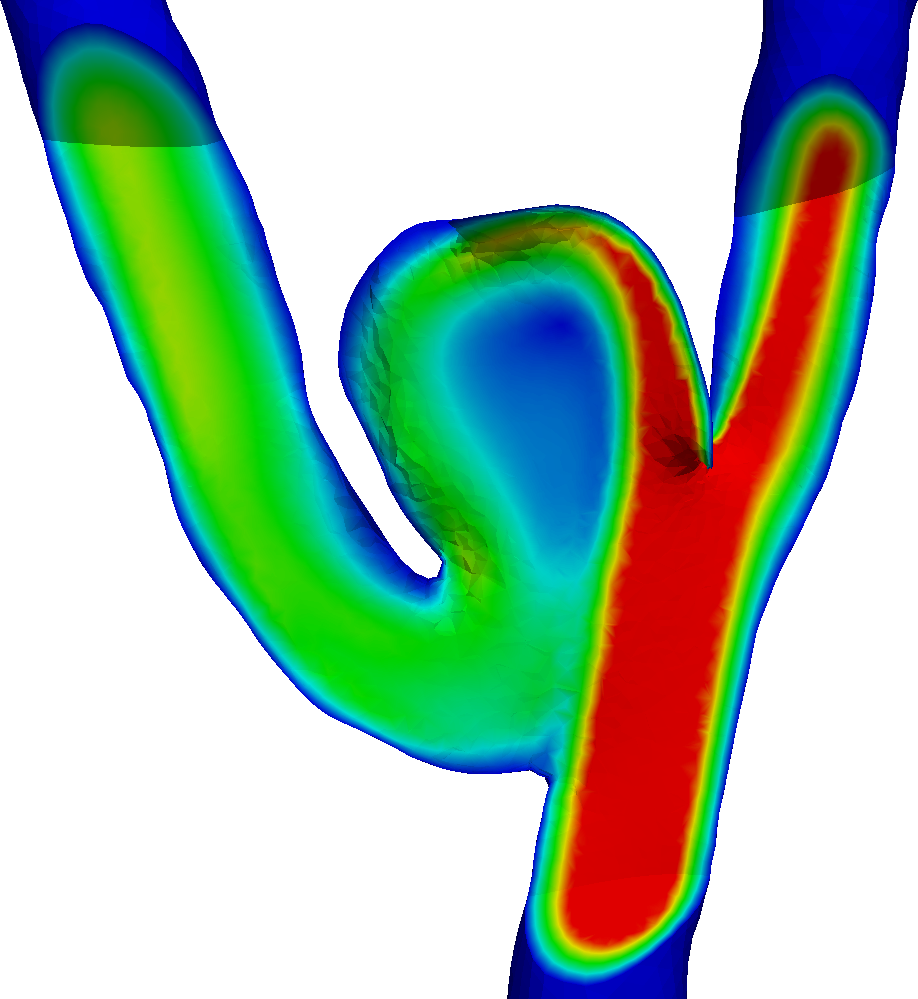}
    \includegraphics[width=0.45\textwidth]{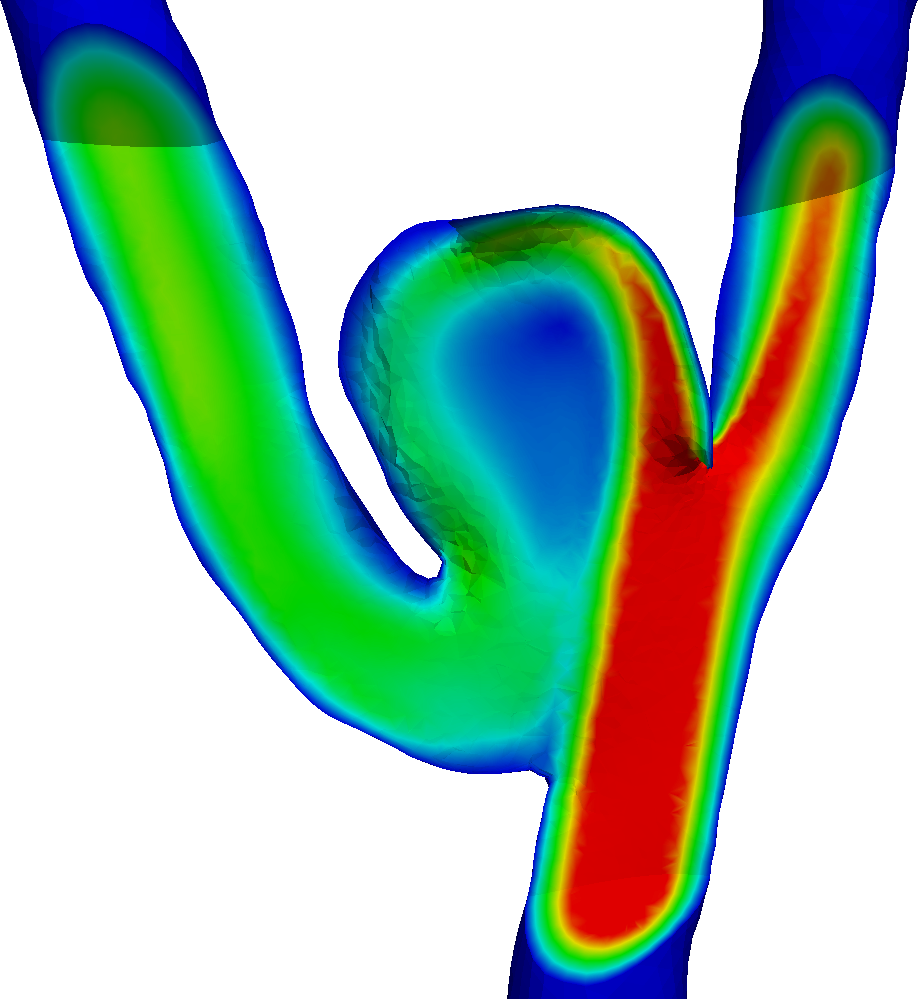}\\
    \rotatebox{90}{\hspace{2.5cm} CFD}\quad
    \includegraphics[width=0.45\textwidth]{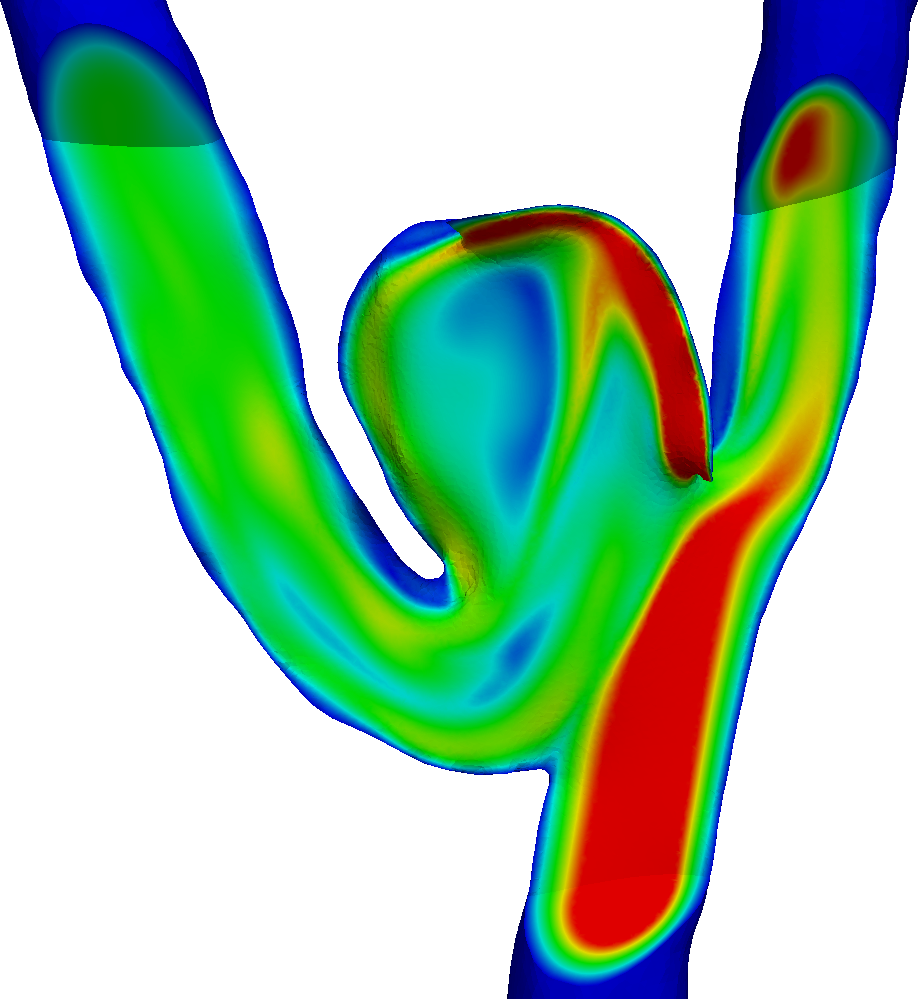}
    \includegraphics[width=0.45\textwidth]{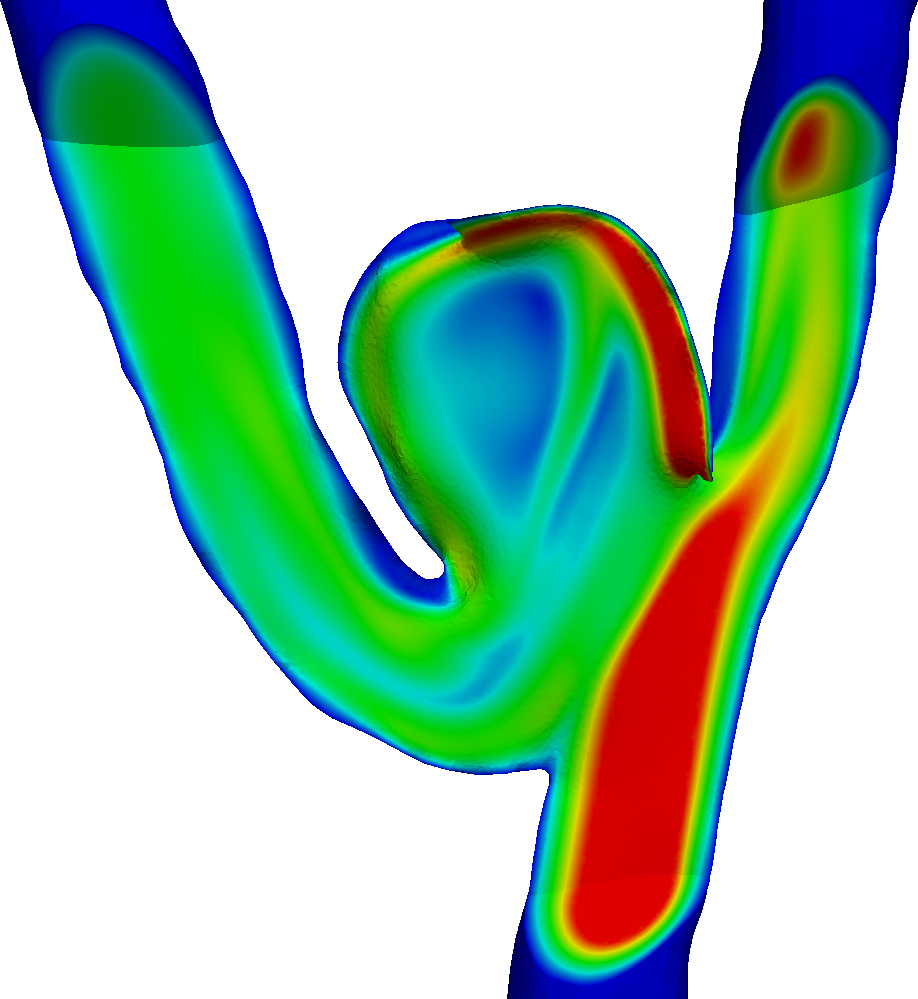}\\
    \rotatebox{90}{\qquad \qquad}\quad
    \parbox[b][1em][t]{0.45\textwidth}{
        \begin{center}
            Instantaneous observations
        \end{center}
    }
    \parbox[b][1em][t]{0.45\textwidth}{
        \begin{center}
            Averaged observations
        \end{center}
    }
    \vspace{0.5cm}
    \begin{center}
        \includegraphics[width=0.4\textwidth]{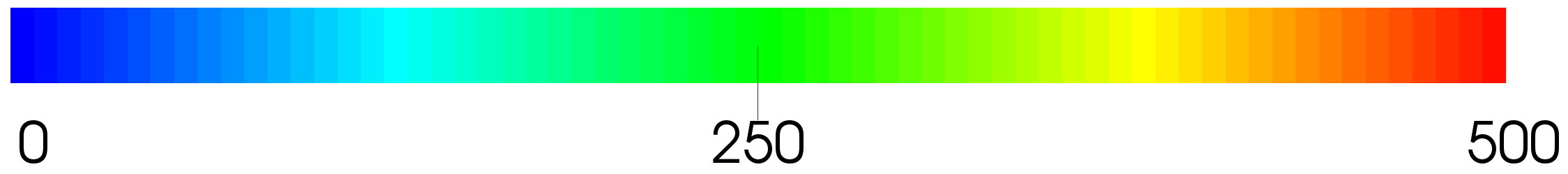}
    \end{center}
    \caption{Flow speed visualised through a slice of the 3D dog vessel. The snapshots are taken at time $t=0.296$s.}
    \label{fig:dog_results_new}
\end{figure}

The results for the data assimilation approach and the high-resolution solver are shown in figure~\ref{fig:dog_results_new}.
The figure shows the data assimilation with the instantaneous observation
operator - the results for the averaged observation operator look similar.
The high-resolution solution has transient to turbulent behaviour in the aneurysm,
while the assimilated solution is laminar. Visually, the assimilated solution
fits better to the observed velocity, both in the vessel and the aneurysm areas.

\section{Conclusion}
This paper presented the application of variational data assimilation to
reconstruct transient blood flow from observations such as MRI images.  This
technique is well known in other scientific fields such as ocean science and
meteorology, but has thus far not been applied to 3D transient blood flow
reconstruction. Mathematically, the data assimilation problem is an
optimisation problem constrained by the Navier-Stokes equations. We derived the
reduced formulation and described the numerical solution with a focus on
retaining the function spaces in the optimisation to obtain mesh-independent
iteration numbers in the optimisation step.

The data assimilation was applied to two examples: first, the reconstruction of
blood flow in an idealised blood vessel with known solution. This example was
used to demonstrate that the proposed method is robust against user parameters and
noisy observations. The second example was based on real 4D MRI measurements in
a three-dimensional domain, and the result compared to a high-resolution blood flow simulation.

Even though the considered blood flow model and observation operators are
simplified, the presented framework extends naturally to more complex setups.
Possible extension is to take into account the movement of the vessel wall,
non-Newtonian effects or a more realistic observation operator that
reimplements an existing measurement device.  Furthermore, the reconstruction
controls could be extended, for example to also reconstruct the
vessel geometry along with the initial and boundary conditions.

The data assimilation procedure introduces an additional computational burden on the
flow reconstruction process - for the discussed examples the data assimilation
is typically around 50 times more computationally expensive than a single flow
simulation.  To keep the computational time feasible, the mesh and time
resolutions had to be reduced compared to a single flow simulation study.  A
simple solution would be to first perform a data assimilation on a coarse
setup, and then apply the reconstructed initial and boundary conditions on a
high-resolution simulation.

 \section*{Acknowledgments}
 This research was supported by The Research Council of Norway through a
 Centres of Excellence grant to the Center for Biomedical Computing at Simula
 Research Laboratory, project number 179578, and a FRIPRO grant, project number 251237.
 Computations were performed on the Abel supercomputing cluster at the University of Oslo via
 NOTUR projects NN9279K and NN9316K.
 In addition, the authors would
 like to thank Jingfeng Jiang and Charles Strother for providing the 4D phase-contrast MRA data, and Lorenz John for comments on the regularisation.

 \bibliographystyle{siamplain}
 \bibliography{literature}
 \end{document}